\documentclass[leqno]{amsart}
\usepackage[T1]{fontenc}
\usepackage[utf8]{inputenc}
\usepackage{lmodern}
\usepackage{microtype}
\usepackage{amsmath, amsthm, amscd, amsfonts, amssymb, mathtools, graphicx, color}
\usepackage{mathrsfs}
\usepackage{enumitem}
\usepackage{etoolbox} 
\usepackage{tikz-cd}
\usepackage{quiver} 
\usepackage[
  bookmarksnumbered,
  plainpages,
  backref=page,
  colorlinks=false,
  pdfborder={0 0 1},
  linkbordercolor={1 0 0},
  citebordercolor={0 0.65 0},
  urlbordercolor={0 0 1}
]{hyperref}
\usetikzlibrary{shapes.geometric, arrows}

\renewcommand*{\backref}[1]{}
\renewcommand*{\backrefalt}[4]{%
  \ifcase #1\relax
  \or \space(Cited on p.~#2.)%
  \else \space(Cited on pp.~#2.)%
  \fi
}

\definecolor{MK_One_Five}{RGB}{90,180,172}
\definecolor{MK_Two_Six}{RGB}{33,102,172}
\definecolor{MK_Three_Six}{RGB}{27,120,55}
\definecolor{MK_Three_One}{RGB}{118,42,131}

\makeatletter
\renewcommand{\tocsection}[3]{%
  \indentlabel{\@ifnotempty{#2}{\bfseries\ignorespaces#1 #2\quad}}\bfseries#3}
\renewcommand{\tocsubsection}[3]{%
  \indentlabel{\@ifnotempty{#2}{\ignorespaces#1 #2\quad}}#3}

\newcommand\@dotsep{4.5}
\def\@tocline#1#2#3#4#5#6#7{\relax
  \ifnum #1>\c@tocdepth 
  \else
    \par \addpenalty\@secpenalty\addvspace{#2}%
    \begingroup \hyphenpenalty\@M
    \@ifempty{#4}{%
      \@tempdima\csname r@tocindent\number#1\endcsname\relax
    }{%
      \@tempdima#4\relax
    }%
    \parindent\z@ \leftskip#3\relax \advance\leftskip\@tempdima\relax
    \rightskip\@pnumwidth plus1em \parfillskip-\@pnumwidth
    #5\leavevmode\hskip-\@tempdima{#6}\nobreak
    \leaders\hbox{$\m@th\mkern \@dotsep mu\hbox{.}\mkern \@dotsep mu$}\hfill
    \nobreak
    \hbox to\@pnumwidth{\@tocpagenum{\ifnum#1=1\bfseries\fi#7}}\par
    \nobreak
    \endgroup
  \fi}
\AtBeginDocument{%
\expandafter\renewcommand\csname r@tocindent0\endcsname{0pt}
}
\def\l@subsection{\@tocline{2}{0pt}{2.5pc}{5pc}{}}
\makeatother

\newtheorem{theorem}{Theorem}[section]
\newtheorem{definition-theorem}[theorem]{Definition-Theorem}
\newtheorem{lemma}[theorem]{Lemma}
\newtheorem{proposition}[theorem]{Proposition}
\newtheorem{prop-def}[theorem]{Proposition-Definition}
\newtheorem{corollary}[theorem]{Corollary}
\newtheorem{conjecture}[theorem]{Conjecture}

\newtheorem{question}[theorem]{Question}

\theoremstyle{definition}
\newtheorem{definition}[theorem]{Definition}
\newtheorem{notation}[theorem]{Notation}
\newtheorem{example}[theorem]{Example}

\theoremstyle{remark}

\newtheorem{claim*}[theorem]{Claim}
\newtheorem{remark}[theorem]{Remark}

\newcommand{\closingmark}{\ensuremath{\lozenge}}

\makeatletter
\newcommand{\closeenv}[1]{%
  \AtBeginEnvironment{#1}{\renewcommand{\qedsymbol}{\closingmark}\pushQED{\qed}}%
  \AtEndEnvironment{#1}{\popQED}%
}
\newcommand{\closeenvs}[1]{%
  \@for\@tempa:=#1\do{\expandafter\closeenv\expandafter{\@tempa}}%
}
\makeatother

\closeenvs{definition,notation,example,exercise,conclusion,criterion,summary,%
axiom,problem,claim,claim*,remark,assumption}

\numberwithin{equation}{section}

\newcommand{\ZZ}{\mathbb{Z}}
\newcommand{\QQ}{\mathbb{Q}}
\newcommand{\RR}{\mathbb{R}}
\newcommand{\CC}{\mathbb{C}}
\newcommand{\HH}{\mathbb{H}}
\newcommand{\Gr}{\mathrm{Gr}}
\newcommand{\GL}{\mathrm{GL}}
\newcommand{\GSp}{\mathrm{GSp}}
\newcommand{\M}{\mathrm{M}} 
\newcommand{\Hom}{\mathrm{Hom}}
\newcommand{\Ext}{\mathrm{Ext}}
\newcommand{\End}{\mathrm{End}}
\newcommand{\Aut}{\mathrm{Aut}}
\newcommand{\pic}{\mathrm{Pic}}
\newcommand{\NS}{\mathrm{NS}} 
\newcommand{\pr}{\mathrm{pr}}
\newcommand{\id}{\mathrm{id}}

\newcommand{\spec}{\operatorname{Spec}}

\newcommand{\calA}{\mathcal{A}}
\newcommand{\calAv}{\mathcal{A}^\vee}
\newcommand{\calB}{\mathcal{B}}
\newcommand{\calY}{\mathcal{Y}}
\newcommand{\calZ}{\mathcal{Z}}
\newcommand{\calE}{\mathcal{E}}
\newcommand{\calX}{\mathcal{X}}
\newcommand{\calG}{\mathcal{G}}
\newcommand{\calH}{\mathcal{H}} 
\newcommand{\calO}{\mathcal{O}}
\newcommand{\calM}{\mathcal{M}} 
\newcommand{\calP}{\mathcal P}

\newcommand{\thetagp}{\mathscr{G}} 

\newcommand{\pl}{\mathcal{P}}          
\newcommand{\pt}{\mathcal{P}^{\times}} 
\newcommand{\lb}{\mathcal{L}}          
\newcommand{\gm}{\mathbb{G}_{\mathrm{m}}}
\newcommand{\Ga}{\mathbb{G}_a}  
\newcommand{\ttt}{\mathbb{T}}   
\newcommand{\chl}{\mathsf{X}^{\bullet}} 

\newcommand{\ag}{S_g}                   
\newcommand{\Ag}{{\mathcal{A}_g}}       
\newcommand{\Agv}{{\mathcal{A}_g^\vee}} 
\newcommand{\pgn}{\mathcal{P}^{\times}_{g,n}} 

\begin{document}
\setcounter{page}{1}

\noindent

\centerline{}

\centerline{}

\title[\scalebox{.8}{On Special Subvarieties of the Universal Semi-abelian Scheme and Pink Conjectures}]{On Special Subvarieties of the Universal Semi-abelian Scheme and Zilber--Pink Conjectures}

\author{Kaiyuan Gu$^1$ and Chenxin Huang$^2$}

\address{$^{1}$ School of Mathematical Sciences, Peking University, Beijing, China.}
\email{\textcolor{MK_One_Five}{gky@stu.pku.edu.cn}}

\address{$^{2}$ Department of Mathematics, Harvard University, Cambridge, MA 02138, USA.}
\email{\textcolor{MK_One_Five}{chuang@math.harvard.edu}}

\dedicatory{To the memory of Professor Bas Edixhoven}



\maketitle

\begin{abstract}
We classify the special subvarieties of the universal semi-abelian scheme of arbitrary toric rank, viewed as a mixed Shimura variety.

As an application, we show that the Zilber--Pink conjecture for mixed Shimura varieties implies its counterpart for semi-abelian varieties, correcting an error in an unpublished manuscript of Pink. Our classification also reveals that the relative Manin--Mumford conjecture for semi-abelian schemes must account for subvarieties of Ribet type, and we give a corrected formulation that does so.
\end{abstract}

\tableofcontents

\section{Introduction}\label{section1}

\subsection{Historical Contexts and Motivations}$\ $

We begin with the following motivating question.
\begin{question}
Let $G$ be a commutative algebraic group over an algebraically closed field $k$ of characteristic zero. If $Z\subset G$ is an irreducible subvariety containing a Zariski dense subset of torsion points of $G$, what can one say about the shape of $Z$?
\end{question}
This is a typical example of a \textbf{special point problem}. In the case of abelian varieties, the Manin--Mumford conjecture claims that such a $Z$ is the translate of an abelian subvariety by a torsion point; it was proved by Raynaud, first for curves \cite{Ray83} and then in general \cite{Ray83b}. For the torus $\gm^{n}$ the Mordell--Lang conjecture, hence a fortiori the statement about torsion points, had been proved by Laurent \cite{Lau84}. The corresponding statement for semi-abelian varieties, and more generally the Mordell--Lang conjecture over number fields, is due to Hindry \cite{Hindry88}, building on Raynaud's theorem and on Serre's results on Galois representations.

In terms of special point problems, for any semi-abelian variety $G$ we call its torsion points the \textbf{special points} and the translates of semi-abelian subvarieties by a torsion point are called the \textbf{special subvarieties}. Then Raynaud's theorem reads: If $Z\subset G$ is an irreducible subvariety containing a Zariski dense subset of special points of $G$, then $Z$ is a special subvariety.

One can consider the analogous question in the context of Shimura varieties, for instance Siegel moduli varieties, where CM points play the role of special points. The analogue of the Manin--Mumford conjecture is then the Andr\'e--Oort conjecture, proposed by Andr\'e for curves in a Shimura variety and by Oort for subvarieties of $\ag$, both with the analogy with Manin--Mumford in mind. After a long series of partial results, many of them conditional on the generalized Riemann hypothesis, the conjecture was approached through the point-counting strategy of Pila and Zannier; the case of $\ag$ was proved by Tsimerman \cite{Tsimerman}, and the general case was claimed by Pila, Shankar and Tsimerman \cite{AO}. We refer to \cite{PilaBook} for a survey of this circle of ideas, which we follow for the historical account below.

The Manin--Mumford conjecture and the Andr\'e--Oort conjecture share formal similarities. One might ask whether there is a common framework such that both conjectures can be viewed as specific cases of a more general special point problem. Given that one of these conjectures is formulated for a single abelian variety, while the other concerns, for example, the moduli space of abelian varieties, a compelling object to investigate is $\Ag$, the universal family over the moduli space of abelian varieties.

The universal abelian scheme $\Ag$ appearing here is an instance of a \emph{mixed Shimura variety}. Consequently, one can still formulate the Andr\'e--Oort conjecture in this setting. By results of Gao \cite{GaoAL, GaomixAO}, the Andr\'e--Oort conjecture for pure Shimura varieties actually implies its analogue for mixed Shimura varieties. This raises the question of whether one can transfer the Andr\'e--Oort statement to the setting of abelian varieties and thereby obtain a family version of the Manin--Mumford conjecture. If a family of abelian varieties contains a Zariski dense set of special points in the Andr\'e--Oort sense, then its base must contain a Zariski dense set of CM points, making the family a Shimura family of abelian varieties. Therefore, for a family that is not Shimura, this approach does not yield any information. Consequently, at this stage, we still lack a full correspondence.

It turns out that a relative version of the Manin--Mumford conjecture can be deduced from an \textbf{unlikely intersection problem}, which involves not only the inclusion of special points but also \textbf{atypical intersections} with special subvarieties. The relevant conjecture here is the Zilber--Pink conjecture. It was arrived at from different motivations and in different formulations by Zilber \cite{Zilber02}, by Bombieri, Masser and Zannier \cite{BMZ99}, and by Pink \cite{Pink05}, and it contains both Manin--Mumford and Andr\'e--Oort as special cases. In the setting of families of abelian varieties this route works, and it leads to the proof of the relative Manin--Mumford conjecture for abelian schemes in \cite{GaoHab_RMM}.

In the unpublished manuscript \cite{Pink05}, Pink proposed a version \cite[Conjecture 6.2]{Pink05} of the relative Manin--Mumford conjecture for semi-abelian schemes. He then used a construction related to the universal semi-abelian variety to deduce it from the Zilber--Pink conjecture. However, in \cite{BE20}, the authors constructed a class of counterexamples called the \emph{Ribet sections} to this version of the relative Manin--Mumford conjecture. While the authors of \cite{BE20} identified some gaps in the arguments of \cite{Pink05}, they also showed that these sections are not pathological from the mixed Shimura viewpoint: they are themselves special subvarieties. Guided by the successful prediction in the case of abelian schemes, this suggests that when working with semi-abelian schemes, there are more special subvarieties to consider than originally believed. Now the natural question is: should the examples of Bertrand--Edixhoven exhaust all the remaining special subvarieties?

In this article, we give an affirmative answer to this question. We accomplish this by thoroughly classifying all special subvarieties of the universal semi-abelian variety, viewed as mixed Shimura varieties. Thus the formulation of the relative Manin--Mumford conjecture might need to be altered to include these exceptions.

\subsection{Main Results}$\ $

Let $S$ be a scheme and let $\calA\to S$ be an abelian scheme with dual $\calAv$. Let $\pl_\calA$ be the Poincar\'e line bundle on $\calA\times_S\calAv$ and $\pt_{\calA}$ the associated $\gm$-torsor, the Poincar\'e torsor; we review its theory in Section~\ref{section2}. For $n\geq1$ put
\[
  (\calAv)^{[n]}:=\underbrace{\calAv\times_{S}\cdots\times_{S}\calAv}_{n\ \mathrm{times}},
  \qquad
  \calB_{\calA,n}:=\calA\times_{S}(\calAv)^{[n]} ,
\]
and define $\pt_{\calA,n}$ by the cartesian diagram
\[\begin{tikzcd}[column sep=6ex, row sep=5ex]
\pt_{\calA,n} \arrow[rr] \arrow[d] \arrow[drr, phantom,"\ulcorner",very near start]
&& \underbrace{\pt_{\calA}\times_{S}\cdots\times_{S}\pt_{\calA}}_{n\ \mathrm{times}} \arrow[d]\\
\calB_{\calA,n} \arrow[rr,"{(\Delta_{\calA},\,\mathrm{id})}"]
&& \underbrace{(\calA\times_{S}\calAv)\times_{S}\cdots\times_{S}(\calA\times_{S}\calAv)}_{n\ \mathrm{times}}\ \rlap{.}
\end{tikzcd}\]
Equivalently, $\pt_{\calA,n}$ is the $n$-fold fiber power of $\pt_\calA$ over $\calA$, so that it is a $\gm^{n}$-torsor over $\calB_{\calA,n}$. Viewed over $(\calAv)^{[n]}$, the scheme $\pt_{\calA,n}$ is a semi-abelian scheme of toric rank $n$, and it is the universal one: this is Proposition~\ref{prop_universal_semiabelian}. For $n=1$ it is the Poincar\'e torsor $\pt_\calA$.

We need another definition to state our main theorem, which seems strange but proves essential. This will be explained with more details in Theorem~\ref{thm_ribet_property}.

\begin{definition}[Definition~\ref{def_ribet_property}]
    For a morphism of $S$-abelian schemes $\calA'\joinrel\xrightarrow{(\alpha,\beta)} \calA\times_S \calAv$, if $\alpha^{\vee}\beta+\beta^{\vee}\alpha=0$, we say that $(\alpha,\beta)$ satisfies the \textbf{Ribet property}.
\end{definition}

For toric rank $n$ the Ribet property is imposed not once but once for each character of the toric part that is to be divided out. Let $\calA'\subseteq\calB_{\calA,n}$ be an abelian subscheme over $S$, with components $\alpha:\calA'\to\calA$ and $\beta_i:\calA'\to\calAv$, and put $\beta_\chi:=\sum_i\chi_i\beta_i$ for a character $\chi\in\chl(\gm^{n})=\ZZ^{n}$. Given a subtorus $\ttt\subseteq\gm^{n}$, we identify $\chl(\gm^{n}/\ttt)$ with the subgroup of characters of $\gm^{n}$ trivial on $\ttt$.

\begin{definition}[Definition~\ref{def_ribet_datum}]
The pair $(\calA',\ttt)$ is a \textbf{Ribet datum} if $(\alpha,\beta_\chi)$ satisfies the Ribet property for every $\chi\in\chl(\gm^{n}/\ttt)$.
\end{definition}

To a Ribet datum we attach in Section~\ref{section4} an irreducible closed subvariety $R(\calA',\ttt)\subseteq\pt_{\calA,n}$, obtained from a \textbf{Ribet multi-section} (Proposition-Definition~\ref{def_ribet_section}) of the $\gm^{n}/\ttt$-torsor over $\calA'$ deduced from $\pt_{\calA,n}$. A \textbf{subvariety of Ribet type} is a torsion translate of some $R(\calA',\ttt)$, the Ribet datum being taken over an irreducible closed subvariety of the base; see Definition~\ref{def_ribet_type}. Three features of this class are worth recording at once.
\begin{itemize}
  \item \emph{It contains the Ribet sections.} For $\ttt=1$ the subvariety $R(\calA',1)$ is the Ribet multi-section itself, and for $n=1$ this is the construction of \cite{BE20}. This is where the name comes from.
  \item \emph{It contains the subgroup schemes.} Every irreducible closed subgroup scheme of $\pt_{\calA,n}|_{T}$ that is flat over an irreducible closed subvariety $T$ of the base is of Ribet type, and so is every torsion translate of one; this is Proposition~\ref{prop_flat_subgroup}.
  \item \emph{Its fibers look like those of a subgroup scheme.} In each fiber, a subgroup scheme is a finite union of translates of one and the same subgroup, namely of its identity component; by Theorem~\ref{thm_ribet_type_fibers} the same is true of a subvariety of Ribet type. The two differ only in whether those translates are torsion, and in whether they come from a single translation over the whole base.
\end{itemize}

We can now specialize to the universal family. Let $\ag$ be the Siegel modular variety parametrizing principally polarized $g$-dimensional abelian varieties with certain level structures (or more generally of an arbitrary polarization type, see Remark~\ref{rem_poltype}). Suppose the level structure is good enough such that $\ag$ is a fine moduli space, whose specific form is unimportant in this article. Then there is a universal abelian scheme $\Ag$ over $\ag$, whose dual is denoted by $\Agv$, and we write
\[
  \pgn:=\pt_{\Ag,n},\qquad \calB_{g,n}:=\calB_{\Ag,n}
\]
for the objects just constructed in this case. Thus $\pgn$ is the universal semi-abelian scheme of toric rank $n$ over $(\Agv)^{[n]}$, and for $n=1$ it is the universal Poincar\'e torsor. It is a mixed Shimura variety, and this is the object inside which we classify the speical subvarieties.

\begin{theorem}[Theorem~\ref{thm_classification}]\label{thm_intro_classification}
Let $R\subseteq\pgn$ be an irreducible closed subvariety. Then $R$ is a special subvariety of $\pgn$, viewed as a mixed Shimura variety, if and only if $R$ is a subvariety of Ribet type, for $\calA=\Ag$ over a special subvariety of $\ag$.
\end{theorem}

The two directions of Theorem~\ref{thm_intro_classification} are proved in rather different flavors. That a subvariety of Ribet type is special, the content of Section~\ref{section5}, comes down to writing a suitable family of Hodge tensors and realizing the subvariety as the special subvariety attached to the mixed Shimura subdatum that those tensors cut out. The converse, the content of Section~\ref{section6}, is instead an interleaving of algebraic geometry with group theory. Both directions rest crucially on the algebro-geometric characterization of the Ribet property obtained in Section~\ref{section2} and on the algebro-geometric properties of subvarieties of Ribet type established in Section~\ref{section4}.

\begin{remark}\label{RmkMSVSpecial}
Since any mixed Shimura variety whose pure part is a Hodge type Shimura variety is a family of semi-abelian varieties, the classification theorem above essentially also gives a classification of special subvarieties for all such mixed Shimura varieties.
\end{remark}

It is worth comparing this with the universal family of abelian varieties over the Siegel moduli space, where a special subvariety is essentially a translation of an algebraic subgroup scheme by a torsion (multi)-section, by \cite[Proposition 1.1]{gaoziyang}. For the universal semi-abelian scheme these no longer exhaust the special subvarieties, and Theorem~\ref{thm_intro_classification} says exactly what has to be added: the Ribet multi-sections.

\subsection{Reformulation of the Relative Manin--Mumford Conjecture}$\ $

Let the base $S$ be a regular, irreducible, quasi-projective variety over $\CC$. The Relative Manin--Mumford Conjecture for an abelian scheme over $S$ was inspired by Shou-wu Zhang's ICM talk in 1998 \cite{Zhang1998} and formulated by Pink \cite{Pink05} and Zannier \cite{Zan12}. The first results are due to Masser and Zannier \cite{MZ10}, who showed that given two explicit non-torsion sections in the Legendre family of elliptic curves, only finitely many parameters make them simultaneously torsion. This is a special case of the Relative Manin--Mumford conjecture for curves in abelian surface schemes. The case of an arbitrary subvariety of an abelian surface scheme was settled by Corvaja, Masser and Zannier \cite{CMZ18}, and the conjecture in full by Gao and Habegger \cite{GaoHab_RMM}. See \cite{Zan12} and \cite[Ch.~20]{PilaBook} for a fuller account.

Pink's route to it is through unlikely intersections. He proposed a conjecture on the intersections of a subvariety of a mixed Shimura variety with the special subvarieties \cite[Conjecture 1.3]{Pink05}, and deduced from it a conjecture for semi-abelian schemes \cite[Conjectures 6.1 and 6.2]{Pink05}, which for families of abelian varieties is the Relative Manin--Mumford Conjecture. That case has recently been proved by Gao and Habegger \cite{GaoHab_RMM}. We restate their theorem in the terminology that our results suggest.

\begin{theorem}[{\cite[Theorem 1.1]{GaoHab_RMM}}]\label{thm_intro_GH}
    Let $\calA_{/S}$ be an abelian scheme of relative dimension $g$, and let $X$ be an irreducible closed subvariety of $\calA$. If $X(\CC)\cap\calA_{\mathrm{tor}}$ is Zariski dense in $X$, then either $\dim X \geq g$, or there are a finite surjective morphism $S'\to S$ and a \emph{proper} subvariety of Ribet type of $\calA\times_S S'$ containing an irreducible component of $X\times_S S'$.
\end{theorem}

For an abelian scheme no Ribet multi-section can occur, so the second alternative says nothing more than that $X$ lies in a subgroup scheme up to a finite base change: a torsion translate of an abelian subscheme over a finite covering of an irreducible closed subvariety of $S$. Calling such a subvariety ``of Ribet type'' is an abuse of terminology, for which we apologize; it is the price of having one word for the exceptional subvarieties in every toric rank.

In positive toric rank the two alternatives no longer exhaust the possibilities. Bertrand and Edixhoven \cite{BE20} exhibited a Ribet section carrying a Zariski dense set of torsion points which satisfies neither: its dimension is too small, and it is not contained in a subgroup scheme up to a finite base change. It is, however, a subvariety of Ribet type, and this is what the following reformulation accounts for.

\begin{conjecture}[Relative Manin--Mumford Conjecture for semi-abelian schemes; Conjecture~\ref{conj_re_RMM_cor}]\label{intro_RMM}
    Let $\calG_{/S}$ be a semi-abelian scheme of relative dimension $g+n$ with toric rank $n$, over an irreducible base. Let $Z$ be an irreducible closed subvariety of $\calG$. If $Z(\CC)\cap\calG_{\mathrm{tor}}$ is Zariski dense in $Z$, then either $\dim Z \geq g+n$, or there are a finite surjective morphism $S'\to S$ and a \emph{proper} subvariety of Ribet type of $\calG\times_S S'$ containing an irreducible component of $Z\times_S S'$.
\end{conjecture}

In both statements the base change is understood to split the toric part, since only then is a subvariety of Ribet type defined; such an $S'$ always exists.

The semi-abelian case of this circle of problems is not an abstraction. Two classical questions lead to it. The first is the Pell--Abel equation $A^{2}-D\,B^{2}=1$, to be solved in polynomials $A,B$ for a given squarefree $D$ of even degree: nontrivial solutions exist exactly when a prescribed divisor class has finite order in the Jacobian of the hyperelliptic curve $u^{2}=D(t)$, and it is a \emph{generalized} Jacobian that intervenes when $D$ has multiple roots. The second is the problem, going back to Abel and Liouville and revived by Davenport, of deciding for which members of a pencil of algebraic differentials the indefinite integral is elementary; by a criterion of Risch this too is a torsion condition in generalized Jacobians. Generalized Jacobians are extensions of a Jacobian by products of copies of $\gm$ and of $\Ga$, so the $\gm$-part of the question is exactly the semi-abelian relative Manin--Mumford problem treated here. Masser and Zannier carried out this analysis and obtained finiteness statements for both questions \cite{MZPell}; see \cite{ZannierICM} and \cite{Zannier8ECM} for surveys, and \cite{Known_for_semiab_surf} for the relative Manin--Mumford conjecture for semi-abelian surfaces. The generalized Jacobian is not merely a source of examples here: in \cite[\S6]{BE20} the counterexample itself is described through the generalized Jacobian of a family of elliptic curves acquiring a double point, a presentation more elementary than the one through the Poincar\'e bundle and more explicit in terms of divisors, rational functions and Weil reciprocity. The $\Ga$-part falls outside the mixed Shimura picture, and so outside the scope of the present article; for it we refer to \cite{Schmidt19}.

In Subsection~\ref{subsec_reform_RMM}, we will prove that Conjecture~\ref{intro_RMM} is a corollary of the Zilber--Pink Conjecture for mixed Shimura varieties.

\subsection{Zilber--Pink for mixed Shimura varieties implies Zilber--Pink for semi-abelian varieties}$\ $

Pink formulated his conjecture on unlikely intersections in two settings: for mixed Shimura varieties \cite[Conjecture 1.3]{Pink05}, where the intersections are with special subvarieties, and for semi-abelian varieties \cite[Conjecture 5.1]{Pink05}, where they are with algebraic subgroups. Since a semi-abelian variety is a fiber of a mixed Shimura variety, one expects the first conjecture to imply the second.

That this is so was announced by Edixhoven in \cite[Remark 1.3]{BE20}, and he communicated a proof to us \cite{EdixhovenPrivate}. In Subsection~\ref{subsec_ZPimpliesZPSA} we give a proof as well.

The two treatments share their starting point, a structural result: for any special subvariety $R$ of the universal semi-abelian scheme $\pgn\rightarrow(\Agv)^{[n]}$, every fiber of $R$ is a finite union of translates of subgroups of the corresponding fiber of $\pgn$. They differ in how it is obtained. We deduce it, as Theorem~\ref{thm_ribet_type_fibers}, from the geometric characterization of the special subvarieties of $\pgn$ in Theorem~\ref{thm_intro_classification}; Edixhoven obtains a finer version, \cite[Proposition 3.9]{EdixhovenPrivate}, by a direct group-theoretic computation generalizing \cite{BE20}, without entering into the geometry of those special subvarieties. That group-theoretic route has an earlier incarnation in the thesis \cite{Lopuhaa14}, which classifies the special algebraic subgroups of the ambient group and, for $g=1$, reads off from that classification the special subvarieties of the Poincar\'e torsor; the symmetry condition it produces is the Ribet property in disguise. See Remark~\ref{rem_group_theoretic_route}.

\subsection{Structure of the Article}$\ $

The structure of this article is as follows. Sections~\ref{section2}, \ref{section3} and \ref{section4} are pure algebraic geometry; no Hodge theory occurs in them. In Section~\ref{section2} we set up the Poincar\'e biextension and its fiber powers, and prove the criterion --- the Ribet property --- for the Poincar\'e bundle to become algebraically trivial on an abelian subscheme. Section~\ref{section3} treats torsion translations, through the theta group and its geometric Heisenberg group, which is what makes the notion of a torsion translate of a subvariety of the torsor precise. In Section~\ref{section4} we introduce Ribet data and construct the subvarieties of Ribet type, generalizing the Ribet sections of \cite{BE20}, and we describe their fibers.

Sections~\ref{section5} and \ref{section6} prove the two directions of Theorem~\ref{thm_intro_classification}. Section~\ref{section5} shows that subvarieties of Ribet type are special; it is the Hodge-theoretic part of the article, can be skipped on a first reading, and its computations largely follow the treatment of \cite{BE20}. Section~\ref{section6} proves the converse, that every special subvariety of $\pgn$ is of Ribet type.

Finally, in Section~\ref{section7} we give a more plausible reformulation of the Relative Manin--Mumford Conjecture and collect the implications of the conjectures involved. We will also prove that Zilber--Pink for mixed Shimura varieties implies Zilber--Pink for semi-abelian varieties. Both are based on our main result Theorem~\ref{thm_intro_classification}.

\subsection*{Acknowledgements}
This work is one of the projects in the summer school ``Algebra and Number Theory'' held by Peking University and the Academy of Mathematics and Systems Science, CAS. We want to express our gratitude to Professor Shou-wu Zhang and Professor Liang Xiao. Without their efforts, the summer school could not be held with such a success. We are also grateful to our mentor, Professor Ziyang Gao, for suggesting this question and for continuous guidance.
We thank Laura DeMarco, Gabriel Dill, Tom Han, Wenbin Luo, Liang Xiao and Xinyi Yuan for discussions and for their interest in this work. We also thank our group members Houjun Chen and Ziheng Zhang for their participation and discussion.
We are grateful to an anonymous referee for comments on an earlier version of this article.
CH would like to thank Peking University for its great hospitality, and the organizers of the Workshop on Arakelov Geometry and Diophantine Geometry, held at Brown University in August 2026, for their hospitality; part of the revision was carried out there.

\section{The Ribet Property}\label{section2}
This section and the two that follow are general algebraic geometry: no Hodge theory or Shimura varieties occur before Section~\ref{section5}. Here we set up the Poincar\'e biextension and its fiber powers, and prove the criterion --- an equation between $\alpha$ and $\beta$ --- for the Poincar\'e bundle to become algebraically trivial when pulled back to an abelian subscheme of $\calA\times_S\calAv$. That criterion, which we call the Ribet property, is what the construction of Section~\ref{section4} rests on, and it is what the classification of Section~\ref{section6} produces.

\subsection{The Poincar\'e Biextension}\label{subsec_poincare_biextension}
In this subsection we review the theory of Poincar\'e biextensions, together with the universal property of the Poincar\'e torsor and of its fiber powers. Our basic reference is \cite[Ch.~I]{faltings-chai}, whose notation differs from ours in two places: the dual abelian scheme is written $A^t$ there and $\calAv$ here, and the homomorphism attached to a line bundle is written $\lambda(\lb)$ there and $\phi_\lb$ here, following \cite{BE20}.

Let $S$ be a scheme, and let $\calA$ be an abelian scheme over $S$ with unit section $0_\calA$. The dual abelian scheme is denoted by $\calAv$ with unit section $0_{\calAv}$. Let $\pl_\calA$ be the Poincar\'e line bundle over $\calA\times_{S}\calAv$, whose formation comes with trivializations along $0_\calA\times_S \calAv$ and $\calA\times_S 0_{\calAv}$, compatible over $0_\calA\times_S 0_{\calAv}$.

\subsubsection{The Poincar\'e Line Bundle and the $\phi$-Construction}

The Poincar\'e line bundle is the universal line bundle over $\calA\times_S \calAv$ in the moduli definition of $\calAv=\pic^{0}_{\calA/S}$. Here, for an $S$-scheme $T$, we write $\pic^{0}_{\calA/S}(T)$ for the group of isomorphism classes of rigidified line bundles: pairs in which $\lb$ is a line bundle on $\calA_T:=\calA\times_S T$ whose restriction to every geometric fiber of $\calA_T\to T$ is algebraically trivial, and $0_{\calA_T}^{*}\lb\stackrel{\sim}{\to}\calO_T$ is a rigidification along the unit section. Rigidifying kills exactly the line bundles pulled back from $T$, so $\pic^{0}_{\calA/S}(T)$ is equally the group of such $\lb$ modulo $\pic(T)$; this is the rigidified Picard functor $\mathbf{Pic}^{0}_{e}(\calA/S)$ of \cite[I.1.4]{faltings-chai}.

The universal property of $\pl_\calA$ takes the following form.

\begin{proposition}\label{prop_universal_poincare}
For every $S$-scheme $T$ the assignment
\begin{equation}\label{eq_universal_poincare}
  b\ \longmapsto\ \pl_b:=(\id_{\calA}\times b)^{*}\pl_\calA,
\end{equation}
where $\pl_b$ carries the rigidification induced by that of $\pl_\calA$ along $0_\calA\times_S\calAv$, is an isomorphism of groups $\calAv(T)\stackrel{\sim}{\longrightarrow}\pic^{0}_{\calA/S}(T)$, natural in $T$. Thus $\calAv$ represents the fppf sheaf $\pic^{0}_{\calA/S}$, with universal object $\pl_\calA$.
\end{proposition}
\begin{proof}
See \cite[I.1.4]{faltings-chai} and the discussion following it; representability is \cite[I.1.9]{faltings-chai}.
\end{proof}

We use \eqref{eq_universal_poincare} constantly and in both directions: it turns a point of $\calAv$ into a line bundle on $\calA$, and it turns a family of fiberwise algebraically trivial line bundles into a morphism to $\calAv$. The first application is that $\pl_\calA$ is bilinear.

\begin{lemma}[{\cite[Exercise (7.4)]{BGB_AV}}]\label{lem_bilinear}
For $m,n\in \ZZ$, the pullback of $\pl_\calA$ via $(m,n):=[m]_\calA\times[n]_{\calAv}$ on $\calA\times_S \calAv$ is given by
\[(m,n)^*\pl_\calA\cong \pl_\calA^{\otimes mn}.\]
In particular $\pl_\calA$ is symmetric: $(-1,-1)^*\pl_\calA\cong\pl_\calA$.
\end{lemma}
\begin{proof}
It suffices to show $(m,1)^*\pl_\calA\cong \pl_\calA^{\otimes m}$, because of the identity $(m,n)=(m,1)\circ(1,n)$ and the self-duality of abelian schemes. View $\calA\times_S \calAv\to \calAv$ as a family of abelian schemes over $\calAv$, with $\pl_\calA$ the universal family of algebraically trivial line bundles, so that $[\pl_\calA] \in \pic^{0}_{\calA\times_S \calAv/\calAv}(\calAv)$ by Proposition~\ref{prop_universal_poincare}. By the linearity of pullbacks on $\pic^{0}$ we have $[(m,1)^{*}\pl_\calA]=[\pl_\calA^{\otimes m}]$ in $\pic^{0}_{\calA\times_S \calAv/\calAv}(\calAv)$. As $(m,1)$ restricts to the identity on $0_\calA\times_S\calAv$, both sides are rigidified there; since $\pic^{0}$ classifies rigidified line bundles, the equality of classes is an isomorphism of line bundles.
\end{proof}

The second application is the $\phi$-construction, which we now recall together with the Mumford line bundle. Both rest on the theorem of the cube, which we quote in the form of \cite[I.1.3]{faltings-chai}: for an invertible sheaf $\lb$ on $\calA$ and $m_I:\calA^{\times_S 3}\to\calA$ the morphism adding the coordinates indexed by $I\subseteq\{1,2,3\}$, with $m_\emptyset$ the zero morphism, the sheaf $\bigotimes_{I}m_{I}^{*}\lb^{\otimes(-1)^{\mathrm{Card}(I)}}$ on $\calA^{\times_S3}$ is canonically trivial.

\begin{definition}[{\cite[I.1.3]{faltings-chai}}]\label{def_mumford_bundle}
For an invertible sheaf $\lb$ on $\calA$ and $m_I:\calA\times_S\calA\to\calA$ as above, now with $I\subseteq\{1,2\}$, set
\[
  \Lambda(\lb):=\bigotimes_{I\subseteq\{1,2\}}m_{I}^{*}\lb^{\otimes(-1)^{\mathrm{Card}(I)}}
  = m^{*}\lb\otimes \pr_1^{*}\lb^{\otimes-1}\otimes \pr_2^{*}\lb^{\otimes-1}\otimes \pi^{*}0_{\calA}^{*}\lb
\]
on $\calA\times_S\calA$, where $m,\pr_1,\pr_2$ are the group law and the two projections and $\pi:\calA\times_S\calA\to S$ is the structure morphism.
\end{definition}

The term $I=\emptyset$ is a normalization making $\Lambda(\lb)$ canonically trivial along $0_\calA\times_S\calA$ and along $\calA\times_S0_\calA$, compatibly, so that $\Lambda(\lb)$ is bi-rigidified exactly like $\pl_\calA$. By the theorem of the cube, its restriction to $\calA\times_S\{a\}$ depends additively on $a$.

For $a\in\calA(T)$ that restriction is $t_a^{*}\lb\otimes\lb^{\otimes-1}$, which lies in $\pic^{0}_{\calA/S}(T)$. By Proposition~\ref{prop_universal_poincare} there is therefore a unique homomorphism of abelian schemes
\[
  \phi_{\lb}:\calA\longrightarrow\calAv,\qquad a\longmapsto \bigl[t_a^{*}\lb\otimes\lb^{\otimes-1}\bigr],
\]
and moreover
\begin{equation}\label{eq_phi_pullback}
  (\id_{\calA},\phi_{\lb})^{*}\pl_\calA\ \cong\ \Lambda(\lb),
\end{equation}
the latter because both sides restrict to $t_a^{*}\lb\otimes\lb^{\otimes-1}$ on $\calA\times_S\{a\}$ and are trivial along $0_\calA\times_S\calA$, so that the Seesaw theorem applies. This is the $\phi$-construction; it is written $\lambda(\lb)$ in \cite[I.1]{faltings-chai}. We record below some properties of the $\phi$-construction.

\begin{lemma}[{\cite[I.1, I.1.6]{faltings-chai}}]\label{lem_phi_properties}
Let $\lb,\calM$ be invertible sheaves on $\calA$.
\begin{enumerate}[label=\textup{(\arabic*)}]
  \item $\phi_{\lb\otimes\calM}=\phi_{\lb}+\phi_{\calM}$.
  \item $\phi_{\lb}=0$ if and only if $\lb\in\pic^{0}_{\calA/S}(S)$; thus $\phi$ factors through $\NS$.
  \item $\phi_{\lb}$ is symmetric: $\phi_{\lb}^{\vee}=\phi_{\lb}$.
  \item If $\lb$ is relatively ample, then $\phi_{\lb}$ is an isogeny.
  \item For a homomorphism $f:\calB\to\calA$ of abelian schemes,
        \begin{equation}\label{eq_phi_functorial}
          \phi_{f^{*}\lb}=f^{\vee}\circ\phi_{\lb}\circ f.
        \end{equation}
\end{enumerate}
\end{lemma}

By Lemma~\ref{lem_phi_properties}(2) the $\phi$-construction forgets exactly the classes in $\pic^{0}_{\calA/S}(S)$, and is therefore injective on $\NS$. One can also go back explicitly: pulling $\pl_\calA$ back along the graph of $\phi_\lb$ returns $\lb$ up to symmetrization, and in particular preserves ampleness.

\begin{lemma}[{\cite[I.1.6]{faltings-chai}}]\label{lem_graph_of_phi}
$(\id_{\calA},\phi_{\lb})^{*}\pl_\calA\cong\lb\otimes[-1]^{*}\lb$, which is $\lb^{\otimes2}$ when $\lb$ is symmetric. In particular, if $\lb$ is relatively ample then the restriction of $\pl_\calA$ to the graph of $\phi_\lb$ is relatively ample, being a tensor product of two relatively ample sheaves.
\end{lemma}

Finally we record the effect of pulling back by $[n]$, which is again a consequence of the theorem of the cube: for an invertible sheaf $\lb$ on $\calA$,
\begin{equation}\label{eq_multiplication_by_n}
  [n]^{*}\lb\ \cong\ \lb^{\otimes\frac{n^{2}+n}{2}}\otimes\bigl([-1]^{*}\lb\bigr)^{\otimes\frac{n^{2}-n}{2}}
\end{equation}
up to an invertible sheaf pulled back from $S$; see \cite[\S6]{Mum12}. In particular $[n]^{*}\lb\cong\lb^{\otimes n^{2}}$ when $\lb$ is symmetric, and $[n]^{*}\lb\cong\lb^{\otimes n}$ when $\lb$ lies in $\pic^{0}_{\calA/S}(S)$.

\subsubsection{The Poincar\'e Torsor and its Fiber Powers}
We turn to the torsor. Recall that a \textbf{$\gm$-biextension} of $\calG_1\times_S\calG_2$, for commutative group schemes $\calG_1,\calG_2$ over $S$, is a $\gm$-torsor over $\calG_1\times_S\calG_2$ together with two compatible partial group laws, one in each variable; equivalently, its restriction to $\calG_1\times_S\{g\}$ is an extension of $\calG_1$ by $\gm$ for every $g$, functorially in $g$, and symmetrically in the other variable. See \cite[VII, VIII]{SGA7I} and \cite[I.1]{faltings-chai}.

\begin{definition}\label{def_poincare_torsor}
The \textbf{Poincar\'e torsor} $\pt_\calA\to\calA\times_S\calAv$ is the $\gm$-torsor associated with $\pl_\calA$; equivalently, it represents the functor $\underline{\mathrm{Isom}}(\calO,\pl_\calA)$,
\begin{align*}
     \mathsf{Sch}_{\calA\times_S\calAv}&\longrightarrow \mathsf{Set}\\
    [T\stackrel{f}{\rightarrow}\calA\times_S\calAv]&\longmapsto\mathrm{Isom}_{\mathrm{Mod}_{\calO_T}}(\calO_T,f^{*}\pl_\calA),
\end{align*}
and is represented by the complement of the zero section in the geometric line bundle of $\pl_\calA$.
\end{definition}

\begin{proposition}[{\cite[I.1]{faltings-chai}}]\label{prop_poincare_biextension}
The bi-rigidification of $\pl_\calA$ makes $\pt_\calA$ the canonical $\gm$-biextension of $\calA\times_S\calAv$. In particular $\pt_\calA\to\calAv$ and $\pt_\calA\to\calA$ are semi-abelian schemes of toric rank $1$: the fiber of the first over $b\in\calAv(T)$ is the extension of $\calA_T$ by $\gm$ of class $b$, and the fiber of the second over $a\in\calA(T)$ is the extension of $\calAv_T$ by $\gm$ of class $a$, under biduality.
\end{proposition}

The second of these two group laws is universal. Recall that for a torus $\ttt$ over $S$ one writes $\chl(\ttt)=\Hom_S(\ttt,\gm)$ for its character group \cite[I.2.1]{faltings-chai}.

\begin{theorem}[Weil--Barsotti, {\cite[I.2]{faltings-chai}}]\label{thm_weil_barsotti}
For an abelian scheme $\calA$ and a torus $\ttt$ over $S$ there is a canonical isomorphism $\Ext^{1}_{S}(\calA,\ttt)\cong\Hom_{\ZZ}(\chl(\ttt),\calAv)$. In particular $\Ext^{1}_{S}(\calA,\gm)\cong\calAv(S)$, and for a split torus $\ttt=\gm^{n}$,
\begin{equation}\label{eq_weil_barsotti_split}
  \Ext^{1}_{S}(\calA,\gm^{n})\ \cong\ \calAv(S)^{n}=(\calAv)^{[n]}(S),
\end{equation}
where $(\calAv)^{[n]}$ is the $n$-fold fiber power of $\calAv$ over $S$; the class of an extension is the tuple of the classes of its pushouts along the $n$ projections $\gm^{n}\to\gm$.
\end{theorem}

For $n=1$ the isomorphism \eqref{eq_weil_barsotti_split} is realized by $\pt_\calA$: by Proposition~\ref{prop_poincare_biextension} the fiber of $\pt_\calA\to\calAv$ over $b$ is the extension of class $b$, so $\pt_\calA\to\calAv$ is the universal extension of $\calA$ by $\gm$. Its fiber powers do the same for higher toric rank, which is why $(\calAv)^{[n]}$ will serve as a moduli space in the rest of the paper.

\begin{definition}\label{def_fiber_power_torsor}
For $n\geq1$ set
\[
  \pt_{\calA,n}:=\underbrace{\pt_\calA\times_{\calA}\cdots\times_{\calA}\pt_\calA}_{n\ \text{times}}
  \ \longrightarrow\
  \calB_{\calA,n}:=\calA\times_S(\calAv)^{[n]},
\]
the fiber power being taken over $\calA$, so that $\pt_{\calA,n}$ is a $\gm^{n}$-torsor over $\calB_{\calA,n}$. We write $\pi:\pt_{\calA,n}\to\calB_{\calA,n}$ for the structure morphism.
\end{definition}

The moduli interpretation of $\pt_{\calA,n}$ is the following.

\begin{proposition}\label{prop_universal_semiabelian}
\begin{enumerate}[label=\textup{(\arabic*)}]
  \item Viewed over $(\calAv)^{[n]}$, the scheme $\pt_{\calA,n}$ is an extension of $\calA\times_S(\calAv)^{[n]}$ by $\gm^{n}$; in particular it is a semi-abelian scheme of toric rank $n$ over $(\calAv)^{[n]}$.
  \item It is the universal such extension: for every $S$-scheme $T$ and every $\underline b\in(\calAv)^{[n]}(T)$, the extension
  \[
    \underline{b}^{\,*}\pt_{\calA,n}:=(\id_{\calA},\underline{b})^{*}\pt_{\calA,n}
  \]
  of $\calA_T$ by $\gm^{n}$ has class $\underline b$. Thus $\underline b\mapsto\underline{b}^{\,*}\pt_{\calA,n}$ is inverse to the isomorphism \eqref{eq_weil_barsotti_split}.
  \item Consequently, if $\calG\to S$ is a semi-abelian scheme whose toric part is a split torus $\gm^{n}$ and whose abelian quotient is $\calA$, and if $\underline b\in(\calAv)^{[n]}(S)$ is its class, then $\calG\cong\underline{b}^{\,*}\pt_{\calA,n}$ as extensions of $\calA$ by $\gm^{n}$.
\end{enumerate}
\end{proposition}
\begin{proof}
(1) By Proposition~\ref{prop_poincare_biextension} the second partial group law makes $\pt_\calA\to\calAv$ an extension of $\calA\times_S\calAv$ by $\gm$. Pulling the $i$-th factor back along the $i$-th projection $\pr_i:(\calAv)^{[n]}\to\calAv$ gives an extension of $\calA\times_S(\calAv)^{[n]}$ by $\gm$ over $(\calAv)^{[n]}$, and $\pt_{\calA,n}$ is by Definition~\ref{def_fiber_power_torsor} the fiber product of these $n$ extensions over $\calA\times_S(\calAv)^{[n]}$. A fiber product of $n$ extensions by $\gm$ over a common base group scheme is an extension by $\gm^{n}$.

(2) Write $\calE_i:=b_i^{*}\pt_\calA$, an extension of $\calA_T$ by $\gm$ of class $b_i$ by Proposition~\ref{prop_poincare_biextension}. By (1), $\underline{b}^{\,*}\pt_{\calA,n}=\calE_1\times_{\calA_T}\cdots\times_{\calA_T}\calE_n$. Its pushout along the $i$-th projection $\gm^{n}\to\gm$ is its quotient by $\prod_{j\neq i}\gm$, and dividing out the $j$-th copy of $\gm$ for every $j\neq i$ forgets all coordinates but the $i$-th; so that pushout is $\calE_i$, of class $b_i$. By the last assertion of Theorem~\ref{thm_weil_barsotti} the class of $\underline{b}^{\,*}\pt_{\calA,n}$ is therefore $\underline b$.

(3) By (2) the extensions $\calG$ and $\underline{b}^{\,*}\pt_{\calA,n}$ have the same class in $\Ext^{1}_{S}(\calA,\gm^{n})$, hence are isomorphic.
\end{proof}

\subsection{The Ribet Property}\label{subsec_ribet_property}
Given a morphism of $S$-abelian schemes $\calA'\stackrel{(\alpha,\beta)}{\longrightarrow} \calA\times_S \calAv$, the Ribet property characterizes when the pullback $(\alpha,\beta)^*\pl_\calA$ lies in $\pic^{0}_{\calA'/S}(S)$, that is, when it is algebraically trivial on the fibers of $\calA'\to S$. We give such a characterization in Theorem~\ref{thm_ribet_property} below, in terms of an equation between $\alpha$ and $\beta$ (and their duals).

\begin{theorem}\label{thm_ribet_property}
Let $S$ be a connected, reduced noetherian scheme, and let $\calA$ be an abelian scheme over $S$. For a morphism of $S$-abelian schemes $\calA'\stackrel{(\alpha,\beta)}{\longrightarrow} \calA\times_S \calAv$, the following statements are equivalent:

(1) $(\alpha,\beta)^*\pl_\calA\in \pic^{0}_{\calA'/S}(S)$;

(2) $(\alpha,\beta)^*\pl_\calA\in \pic^{0}_{\calA'/S}(S)[2]$, in other words, $((\alpha,\beta)^*\pl_\calA)^{\otimes 2}$ is pulled back from $S$;

(3) $\alpha^{\vee}\beta+\beta^{\vee}\alpha=0$ holds in $\Hom_S(\calA',\calA'^\vee)$.
\end{theorem}

We extract the equivalent condition (3) from Theorem~\ref{thm_ribet_property} as a definition, which we call the \textbf{Ribet property}.

\begin{definition}\label{def_ribet_property}
  Let $S$ be a scheme, and let $\calA$ be an abelian scheme over $S$. We say that a morphism of $S$-abelian schemes $\calA'\stackrel{(\alpha,\beta)}{\longrightarrow} \calA\times_S \calAv$ satisfies the \textbf{Ribet property} if 
  \[\alpha^{\vee}\beta+\beta^{\vee}\alpha=0\]
  holds in $\Hom_S(\calA',\calA'^\vee)$. 
\end{definition}

\begin{example}\label{ex_ribet_graph}
  Let $f:\calA\rightarrow \calAv$ be a homomorphism of $S$-abelian schemes. Then the inclusion of the graph of $f-f^{\vee}$ into $\calA\times_S\calAv$ satisfies the Ribet property. For instance, if $\calA$ is a CM elliptic curve over a field, one can choose $f$ so that $f-f^{\vee}$ is non-zero. This example is the case studied in \cite{BE20}. In this case, $(\id_\calA, f-f^{\vee})^*\pl_\calA$ is pulled back from $S$, without having to take the second tensor power. In general, one can only guarantee that $((\alpha,\beta)^*\pl_\calA)^{\otimes 2}$ is pulled back from $S$.
\end{example}

Before going into the proof of Theorem~\ref{thm_ribet_property}, we discuss a structural consequence of the property just defined. 

The Ribet property can be viewed as an isotropy condition for the Poincar\'e bundle, and this is the form in which it first appeared. An abelian subvariety of $\calA\times_S\calAv$ to which $\pl_\calA$ restricts algebraically trivially is called \emph{isotropic} in \cite{Bertrand96}; the name records that the first Chern class of $\pl_\calA$ is an alternating form, the condition asking the subvariety to be isotropic for it. The notion is extended to one-motives in \cite[p.~5]{Bertrand98}, and condition (i) of that definition is already Proposition-Definition~\ref{def_ribet_section} for $n=1$ and $\ttt=1$: the restriction being algebraically trivial and $\pl_\calA$ being symmetric and rigidified, that restriction acquires a canonical splitting up to a $2$-isogeny. See also \cite[\S8.3]{Lopuhaa14}, where the isotropic subvarieties of $A\times A^{\vee}$ are classified for simple $A$. In a symplectic vector space, an isotropic subspace has dimension at most half of the ambient dimension. The analogous dimension bound associated with the Ribet property is as follows.

\begin{proposition}\label{prop_dimension_bound}
Let $A$ be an abelian variety over a field $k$ of dimension $g\geq 1$, and let $A'\subset A\times A^\vee$ be an abelian subvariety whose inclusion satisfies the Ribet property. Then the dimension bound
\[\dim A'\leq g\]
holds.
\end{proposition}

\begin{proof}
We argue by contradiction, assuming $\dim A'>g$. Choose an ample line bundle $\lb$ on $A$ and let $\phi_\lb:A\rightarrow A^\vee$ be the associated polarization, an isogeny by Lemma~\ref{lem_phi_properties}(4). Let $\Gamma(\phi_\lb)\subset A\times A^{\vee}$ be the graph of $\phi_\lb$. Finally, let $A''=(\Gamma(\phi_\lb)\cap A')^{\circ}$ be the identity component of $\Gamma(\phi_\lb)\cap A'$, which is a non-trivial abelian subvariety of $A\times A^\vee$ because $\dim\Gamma(\phi_\lb)+\dim A'>2g=\dim(A\times A^{\vee})$.

Consider the restriction of the Poincar\'e line bundle $\pl_A$ to $A''$. On one hand, $\pl_A|_{A''}$ is in $\pic^{0}(A'')$: by Theorem~\ref{thm_ribet_property} the Ribet property gives $\pl_A|_{A'}\in\pic^{0}(A')$, and restriction along $A''\subset A'$ preserves $\pic^{0}$. On the other hand, $\pl_A|_{A''}$ is ample since $\pl_A|_{\Gamma(\phi_{\lb})}$ is ample by Lemma~\ref{lem_graph_of_phi}. On a positive-dimensional abelian variety no line bundle is both ample and algebraically trivial, a contradiction.
\end{proof}

The remaining part of this subsection is devoted to the proof of Theorem~\ref{thm_ribet_property}. Besides the preliminaries of Subsection~\ref{subsec_poincare_biextension}, we need one auxiliary fact.

\begin{lemma}\label{lem_square_pic0}
Let $\calA'_{/S}$ be an abelian scheme with $S$ connected and reduced, and let $\calM$ be a line bundle on $\calA'$. For any nonzero integer $n$, $\calM^{\otimes n}\in\pic^{0}_{\calA'/S}(S)$ if and only if $\calM\in\pic^{0}_{\calA'/S}(S)$.
\end{lemma}
\begin{proof}
The ``only if'' direction is immediate since $\pic^{0}_{\calA'/S}(S)$ is a subgroup of $\pic(\calA')$. For the ``if'' direction, Lemma~\ref{lem_phi_properties}(1) gives $\phi_{\calM^{\otimes n}}=n\phi_\calM$, so if $\calM^{\otimes n}\in\pic^{0}_{\calA'/S}(S)$ then $n\phi_\calM=0$ in $\Hom_S(\calA',\calA'^\vee)$ by Lemma~\ref{lem_phi_properties}(2). As $\calA'\to S$ has connected fibers, $\Hom_S(\calA',\calA'^\vee)$ is a torsion-free $\ZZ$-module, so $n\phi_\calM=0$ forces $\phi_\calM=0$, i.e.\ $\calM\in\pic^{0}_{\calA'/S}(S)$ by Lemma~\ref{lem_phi_properties}(2) again.
\end{proof}

Throughout the proof we abbreviate $\pl_\calA|_{\calA'}:=(\alpha,\beta)^*\pl_\calA$.

\begin{proof}[Proof of Theorem~\ref{thm_ribet_property}]

\smallskip
\noindent\textit{Step 1: (2)$\Rightarrow$(1), for general $S$.} This is exactly
Lemma~\ref{lem_square_pic0} with $n=2$, applied to $\calM=\pl_\calA|_{\calA'}$.

\smallskip
\noindent\textit{Step 2: the case $S=\spec(k)$, $k$ algebraically closed --- here we establish the
full equivalence (1)$\Leftrightarrow$(3), together with (3)$\Rightarrow$(2).}

Choose an ample line bundle $\lb$ on $\calA$ and an isogeny $\psi:\calAv\to\calA$ with
$\phi_\lb\circ\psi=[N]_{\calAv}$ for a positive integer $N$, and set $\gamma:=\psi\circ\beta$. Let
\[
  \calM_{(\alpha,\gamma)}(\lb):=(\alpha,\gamma)^{*}\Lambda(\lb)
  =(\alpha+\gamma)^*\lb\otimes\alpha^{*}\lb^{\otimes-1}\otimes\gamma^{*}\lb^{\otimes-1},
\]
where $(\alpha,\gamma):\calA'\to\calA\times_S\calA$ and the second equality is
Definition~\ref{def_mumford_bundle}, the term indexed by $I=\emptyset$ being trivial over a field.
We first check that $\alpha^\vee\beta+\beta^\vee\alpha=0$ is equivalent to
$\calM_{(\alpha,\gamma)}(\lb)\in\pic^{0}(\calA')$: indeed, by Lemma~\ref{lem_phi_properties}(1) and
\eqref{eq_phi_functorial},
\begin{align*}
    \phi_{\calM_{(\alpha,\gamma)}(\lb)}&=\phi_{(\alpha+\gamma)^*\lb}-\phi_{\alpha^*\lb}-\phi_{\gamma^*\lb}\\
    &=(\alpha+\gamma)^\vee\circ\phi_\lb\circ (\alpha+\gamma)-\alpha^\vee\circ\phi_\lb\circ\alpha-\gamma^\vee\circ\phi_\lb\circ\gamma
    \\&=\alpha^\vee\circ\phi_\lb\circ\gamma+\gamma^\vee\circ\phi_\lb\circ\alpha
    \\&=\alpha^\vee\circ[N]_{\calAv}\circ\beta+\beta^\vee\circ[N]_{\calA}\circ\alpha
    \\
    &=(\alpha^\vee\beta+\beta^\vee\alpha)\circ [N]_{\calA'},
\end{align*}
where the fourth equality uses
\[
  \phi_\lb\gamma=\phi_\lb\psi\beta=[N]_{\calAv}\beta
  \qquad\text{and, dually,}\qquad
  \gamma^\vee\phi_\lb=\beta^\vee\psi^\vee\phi_\lb=\beta^\vee[N]_{\calA},
\]
the latter because $\phi_\lb^\vee=\phi_\lb$ by Lemma~\ref{lem_phi_properties}(3). Since $\Hom_S(\calA',\calA'^\vee)$ is
torsion free, the equivalence follows from Lemma~\ref{lem_phi_properties}(2).

Next we identify $\calM_{(\alpha,\gamma)}(\lb)$. By \eqref{eq_phi_pullback},
\[
  \calM_{(\alpha,\gamma)}(\lb)=(\alpha,\gamma)^{*}(\id_\calA,\phi_\lb)^{*}\pl_\calA
  =(\alpha,\phi_\lb\gamma)^{*}\pl_\calA=(\alpha,N\beta)^{*}\pl_\calA,
\]
and $(\alpha,N\beta)=(1,N)\circ(\alpha,\beta)$, so Lemma~\ref{lem_bilinear} gives
\begin{equation}\label{eq_M_is_PN}
  \calM_{(\alpha,\gamma)}(\lb)\cong(\alpha,\beta)^{*}\pl_\calA^{\otimes N}=(\pl_\calA|_{\calA'})^{\otimes N}.
\end{equation}
Combining the two displays with Lemma~\ref{lem_square_pic0} (with $n=N$) gives the chain of
equivalences
\begin{align*}
    (3) &\iff \calM_{(\alpha,\gamma)}(\lb)\in\pic^{0}(\calA') && \text{(computation above)}\\
    &\iff (\pl_\calA|_{\calA'})^{\otimes N}\in\pic^{0}(\calA') && \text{(by \eqref{eq_M_is_PN})}\\
    &\iff \pl_\calA|_{\calA'}\in\pic^{0}(\calA') && \text{(Lemma~\ref{lem_square_pic0}, $n=N$)}\\
    &\iff (1),
\end{align*}
which is exactly (1)$\Leftrightarrow$(3).

For (3)$\Rightarrow$(2): assume $\alpha^\vee\beta+\beta^\vee\alpha=0$, i.e.\ (3). By the equivalence
just established, (3) also gives (1), i.e.\ $\pl_\calA|_{\calA'}\in\pic^{0}(\calA')$; in particular
$\pl_\calA|_{\calA'}$ is anti-symmetric, i.e.\ $[-1]_{\calA'}^*(\pl_\calA|_{\calA'})\cong(\pl_\calA|_{\calA'})^{\otimes-1}$
(for $L\in\pic^{0}$ the seesaw theorem gives $m^{*}L\cong\pr_1^{*}L\otimes\pr_2^{*}L$, and restricting to the antidiagonal gives $[-1]^{*}L\cong L^{\otimes-1}$). Combining this with
\eqref{eq_multiplication_by_n} for $n=2$, namely $[2]_{\calA'}^*L\cong L^{\otimes3}\otimes[-1]_{\calA'}^*L$,
applied to $L=\pl_\calA|_{\calA'}$, gives
\[
  [2]_{\calA'}^*\pl_\calA|_{\calA'}\cong(\pl_\calA|_{\calA'})^{\otimes3}\otimes(\pl_\calA|_{\calA'})^{\otimes-1}=(\pl_\calA|_{\calA'})^{\otimes2}.
\]
On the other hand, since $(\alpha,\beta)\circ[2]_{\calA'}=(2,2)\circ(\alpha,\beta)$ (as $(\alpha,\beta)$
is a homomorphism), Lemma~\ref{lem_bilinear} with $(m,n)=(2,2)$ gives
\[
  [2]_{\calA'}^*\pl_\calA|_{\calA'}\cong(\alpha,\beta)^*\bigl((2,2)^*\pl_\calA\bigr)\cong(\alpha,\beta)^*\pl_\calA^{\otimes4}=(\pl_\calA|_{\calA'})^{\otimes4},
\]
and this holds regardless of (3). Combining the two displayed isomorphisms gives
$(\pl_\calA|_{\calA'})^{\otimes4}\cong(\pl_\calA|_{\calA'})^{\otimes2}$, and cancelling shows
$(\pl_\calA|_{\calA'})^{\otimes2}$ is trivial, i.e.\ (2) holds.

This establishes (1)$\Leftrightarrow$(3)$\Rightarrow$(2) --- and hence, together with Step 1, the full
equivalence (1)$\Leftrightarrow$(2)$\Leftrightarrow$(3) --- when $S=\spec(k)$.

\smallskip
\noindent\textit{Step 3: the case of a general connected reduced noetherian base $S$.}

\noindent\emph{(3)$\Rightarrow$(2).} Assume (3). Since the formation of
$(\pl_\calA|_{\calA'})^{\otimes 2}$ is compatible with base change, Step 2 shows
$(\pl_\calA|_{\calA'_{\bar s}})^{\otimes 2}$ is trivial on every geometric fiber $\calA'_{\bar s}$ of
$\calA'\to S$. By the Seesaw theorem, there is a maximal closed subscheme $K\hookrightarrow S$ over
which $(\pl_\calA|_{\calA'})^{\otimes 2}$ is trivial (i.e.\ pulled back from $K$, with $K$ universal
for this property); maximality forces every geometric point of $S$ to factor through $K$, so $K=S$
since $S$ is reduced.

\noindent\emph{(1)$\Rightarrow$(3).} Assume (1). For every geometric point $\bar s\to S$,
$\pl_\calA|_{\calA'_{\bar s}}\in\pic^{0}(\calA'_{\bar s})$, so Step 2 gives
$(\alpha^\vee\beta+\beta^\vee\alpha)_{\bar s}=0$. As this holds on every (geometric) fiber of
$\calA'\to S$, the homomorphism $\alpha^\vee\beta+\beta^\vee\alpha$ is itself zero.

Combining Steps 1--3 closes the cycle (1)$\Rightarrow$(3)$\Rightarrow$(2)$\Rightarrow$(1).
\end{proof}

\section{Torsion Translations}\label{section3}

In this section we introduce the notion of torsion translation on $\pt_{\calA,n}$. After passing to suitable finite covers induced by multiplication on $\calB_{\calA,n}$, such translations are encoded by the action of the associated theta groups. This formalism will be used in Section~\ref{section4} to define subvarieties of Ribet type, in Section~\ref{section6} to remove torsion translations in the classification, and in Section~\ref{section7} to analyze their fibers.


Our reference is \cite[Ch.~I]{moret-bailly}, which treats theta groups over an arbitrary base. We follow its notation for the objects it introduces, $K(\lb)$ and $\thetagp(\lb)$, and keep that of Subsections~\ref{subsec_poincare_biextension} and~\ref{subsec_ribet_property} for everything else. 

Throughout, $\calA\to S$ is an abelian scheme as in Subsection~\ref{subsec_poincare_biextension}, and we keep the notation $\pt_{\calA,n}\to\calB_{\calA,n}=\calA\times_S(\calAv)^{[n]}$ of Definition~\ref{def_fiber_power_torsor}. We write $\pr_0:\calB_{\calA,n}\to\calA$ and $\pr_i:\calB_{\calA,n}\to\calAv$ for the projections, $1\leq i\leq n$, and
\[
  \pl^{(i)}:=(\pr_0,\pr_i)^{*}\pl_\calA ,
\]
so that $\pt_{\calA,n}$ is the fiber product over $\calB_{\calA,n}$ of the $n$ associated $\gm$-torsors. Since $\pl_\calA$ is bi-rigidified, Proposition~\ref{prop_universal_poincare} applies equally to $\calAv$ and its dual $\calA$; we therefore extend the notation $\pl_b=(\id_\calA\times b)^{*}\pl_\calA$ of \eqref{eq_universal_poincare} by setting $\pl_a:=(a\times\id_{\calAv})^{*}\pl_\calA$ for $a\in\calA(T)$, a line bundle on $\calAv_T$. Which of the two is meant is always determined by the subscript.

\subsection{The Theta Group}

We first recall the theta group. Let $\calX\to S$ be an abelian scheme and let $\lb$ be a line bundle on $\calX$, rigidified along the unit section. Write
\[
  \pi_\lb:\ \lb^{\times}\longrightarrow\calX
\]
for the associated $\gm$-torsor, as in Definition~\ref{def_poincare_torsor}, $\lb^{\times}_x$ for its fiber above a point $x$, and $\Aut_{\gm}(\lb^{\times}_T)$ for the group of automorphisms of $\lb^{\times}_T$ as a $\gm$-torsor.

\begin{definition}[{\cite[I.4.3]{moret-bailly}}]\label{def_theta_group}
The \textbf{theta group} of $\lb$ is the group functor
\begin{align*}
  \thetagp(\lb):\ \mathsf{Sch}_{S}\ &\longrightarrow\ \mathsf{Grp}\\
  T\ &\longmapsto\ \Bigl\{\,\sigma\in\Aut_{\gm}\bigl(\lb^{\times}_T\bigr)\ \Bigm|\ \pi_\lb\circ\sigma=t_x\circ\pi_\lb\ \text{ for some } x\in\calX(T)\,\Bigr\},
\end{align*}
the group law being composition of automorphisms.
\end{definition}

The point $x$ in the definition is unique, since $\pi_\lb$ is surjective, so $\sigma\mapsto x$ is a homomorphism of group functors
\[
  \rho:\ \thetagp(\lb)\longrightarrow\calX .
\] Following \cite[I.4.2]{moret-bailly} we set
\[
  K(\lb):=\ker\phi_{\lb}\ \subset\ \calX ,
\]
a closed subgroup scheme of $\calX$.

\begin{lemma}\label{lem_theta_sheaf_exact}
The kernel of $\rho$ is $\gm$ and is central in $\thetagp(\lb)$, and the image of $\rho$ is $K(\lb)$. Thus
\begin{equation}\label{eq_theta_group}
  1\longrightarrow\gm\longrightarrow\thetagp(\lb)\stackrel{\rho}{\longrightarrow}K(\lb)\longrightarrow1
\end{equation}
is an exact sequence of fppf sheaves of groups over $S$.
\end{lemma}
\begin{proof}
The kernel consists of the automorphisms of $\lb^{\times}_T$ lying over $\id_{\calX_T}$; these are the multiplications by invertible functions on $\calX_T$, and $\calO^{*}(\calX_T)=\calO^{*}(T)$ since $\calX\to S$ is proper with geometrically connected fibers and a section. So $\ker\rho=\gm$, and it is central in $\thetagp(\lb)$ because $\gm$ acts on $\lb^{\times}$ through the torsor structure.

For the image, observe that an automorphism of $\lb^{\times}_T$ lying over $t_x$ is the same thing as an isomorphism of line bundles $\lb_T\stackrel{\sim}{\to}t_x^{*}\lb_T$. Such an isomorphism can exist only if $t_x^{*}\lb_T\otimes\lb_T^{\otimes-1}$ is trivial, hence only if $\phi_\lb(x)=0$, that is, only if $x$ lies in $K(\lb)$. Conversely, if $x\in K(\lb)(T)$ then $t_x^{*}\lb_T\otimes\lb_T^{\otimes-1}$ is pulled back from $T$, so such isomorphisms exist locally on $T$ and $\rho$ is surjective as a map of fppf sheaves.
\end{proof}

We shall also use the standard geometric realization of the theta group. Restricting \eqref{eq_phi_pullback} to $K(\lb)\times_S\calX$, and using $\phi_\lb|_{K(\lb)}=0$, gives a canonical trivialization of $\Lambda(\lb)$ there. Fiberwise, this is a bilinear operation
\[
  \ast:\ \lb^{\times}_x\times\lb^{\times}_y\longrightarrow\lb^{\times}_{x+y},
  \qquad x\in K(\lb),\ y\in\calX .
\]
For $y\in K(\lb)$ it endows $\lb^{\times}|_{K(\lb)}$ with a group structure, which is also a central extension of $K(\lb)$ by $\gm$; for arbitrary $y$ it gives an action of $\lb^{\times}|_{K(\lb)}$ on $\lb^{\times}$ lifting translations of $K(\lb)$ on $\calX$. See \cite[I.3--I.4]{moret-bailly}.

\begin{lemma}[{\cite[I.4.4]{moret-bailly}}]\label{lem_theta_representable}
The map
\[
  \lb^{\times}|_{K(\lb)}\longrightarrow\thetagp(\lb),\qquad
  u\longmapsto(v\mapsto u\ast v),
\]
is an isomorphism of central extensions. In particular $\thetagp(\lb)$ is representable by a group scheme over $S$, and \eqref{eq_theta_group} is an exact sequence of group schemes over $S$.
\end{lemma}

We now specialize to $\calX=\calB_{\calA,n}$ and $\lb=\pl^{(i)}$. A translation of $\calB_{\calA,n}$ lifts to $\pt_{\calA,n}$ exactly when it lifts in every factor $(\pl^{(i)})^{\times}$, hence when it lies in $\bigcap_iK(\pl^{(i)})$. The next two lemmas compute these kernels.

\begin{lemma}\label{lem_translate_poincare}
Let $T$ be an $S$-scheme and let $c=(a,\underline b)\in\calB_{\calA,n}(T)$. Then, up to a line bundle pulled back from $T$,
\[
  t_c^{*}\pl^{(i)}\ \cong\ \pl^{(i)}\otimes\pr_0^{*}\pl_{b_i}\otimes\pr_i^{*}\pl_{a}
  \qquad(1\leq i\leq n).
\]
\end{lemma}
\begin{proof}
Since $\pl^{(i)}=(\pr_0,\pr_i)^{*}\pl_\calA$ and $(\pr_0,\pr_i)\circ t_c=t_{(a,b_i)}\circ(\pr_0,\pr_i)$, pulling back along $(\pr_0,\pr_i)$ reduces the statement to the case $n=1$, namely to
\[
  t_{(a,b)}^{*}\pl_\calA\ \cong\ \pl_\calA\otimes\pr_0^{*}\pl_{b}\otimes\pr_1^{*}\pl_{a}
\]
on $\calA_T\times_T\calAv_T$, where $\pr_0$ and $\pr_1$ are its two projections and $a\in\calA(T)$, $b\in\calAv(T)$. Now $\phi_{\pl_\calA}$ is a homomorphism of abelian schemes, this being the theorem of the square for $\pl_\calA$ on $\calA\times_S\calAv$; applied to $(a,b)=(a,0)+(0,b)$ it gives
\[
  t_{(a,b)}^{*}\pl_\calA\otimes\pl_\calA^{\otimes-1}\ \cong\
  \bigl(t_{(a,0)}^{*}\pl_\calA\otimes\pl_\calA^{\otimes-1}\bigr)\otimes
  \bigl(t_{(0,b)}^{*}\pl_\calA\otimes\pl_\calA^{\otimes-1}\bigr).
\]
It therefore suffices to treat the two one-variable cases
\begin{equation}\label{eq_partial_translation}
  t_{(a,0)}^{*}\pl_\calA\ \cong\ \pl_\calA\otimes\pr_1^{*}\pl_{a},
  \qquad
  t_{(0,b)}^{*}\pl_\calA\ \cong\ \pl_\calA\otimes\pr_0^{*}\pl_{b}.
\end{equation}
These are precisely the two partial group laws of the Poincar\'e biextension, Proposition~\ref{prop_poincare_biextension}.
\end{proof}

\begin{lemma}\label{lem_K_of_PN}
Identify $\calB_{\calA,n}^{\vee}$ with $\calAv\times_S\calA^{[n]}$, and let $N\geq1$.
\begin{enumerate}[label=\textup{(\arabic*)}]
  \item The homomorphism $\phi_{\pl^{(i)}}$ is the $i$-th transposition,
  \[
    \phi_{\pl^{(i)}}:\calB_{\calA,n}\longrightarrow\calB_{\calA,n}^{\vee},
    \qquad(a,\underline b)\longmapsto(b_i;\,0,\dots,0,a,0,\dots,0),
  \]
  with $a$ in the $i$-th slot.
  \item Consequently
  \[
    \phi_{[N]_{\calB_{\calA,n}}^{*}\pl^{(i)}}=N^{2}\,\phi_{\pl^{(i)}},
  \]
  and therefore
  \[
    K\bigl([N]_{\calB_{\calA,n}}^{*}\pl^{(i)}\bigr)
    =\bigl\{(a,\underline b)\ :\ N^{2}a=0\ \text{ and }\ N^{2}b_i=0\bigr\}.
  \]
  \item In particular
  \[
    \bigcap_{i=1}^{n}K\bigl([N]_{\calB_{\calA,n}}^{*}\pl^{(i)}\bigr)=\calB_{\calA,n}[N^{2}].
  \]
\end{enumerate}
\end{lemma}
\begin{proof}
The functor $\pic^{0}_{\calB_{\calA,n}/S}$ is the product of $\pic^{0}_{\calA/S}$ with $n$ copies of $\pic^{0}_{\calAv/S}$, pullback along $\pr_0,\pr_1,\dots,\pr_n$ giving the isomorphism. Combined with Proposition~\ref{prop_universal_poincare} for $\calA$ and for $\calAv$, this is the stated identification of $\calB_{\calA,n}^{\vee}$ with $\calAv\times_S\calA^{[n]}$: a tuple $(\beta;\alpha_1,\dots,\alpha_n)$, with $\beta\in\calAv(T)$ and $\alpha_j\in\calA(T)$, corresponds to the class of
\[
  \pr_0^{*}\pl_{\beta}\otimes\bigotimes_{j=1}^{n}\pr_j^{*}\pl_{\alpha_j}.
\]
Now by definition $\phi_{\pl^{(i)}}(c)$ is the class of $t_c^{*}\pl^{(i)}\otimes(\pl^{(i)})^{\otimes-1}$, which Lemma~\ref{lem_translate_poincare} computes to be $\pr_0^{*}\pl_{b_i}\otimes\pr_i^{*}\pl_{a}$. Under the dictionary just described this is the tuple $(b_i;0,\dots,0,a,0,\dots,0)$, which proves (1).

The dual of $[N]_{\calB_{\calA,n}}$ is multiplication by $N$ on $\calB_{\calA,n}^{\vee}$, so \eqref{eq_phi_functorial} gives
\[
  \phi_{[N]_{\calB_{\calA,n}}^{*}\pl^{(i)}}
  =[N]_{\calB_{\calA,n}}^{\vee}\circ\phi_{\pl^{(i)}}\circ[N]_{\calB_{\calA,n}}
  =N^{2}\,\phi_{\pl^{(i)}} .
\]
By (1), a point $c=(a,\underline b)$ is therefore killed by $\phi_{[N]_{\calB_{\calA,n}}^{*}\pl^{(i)}}$ if and only if $N^{2}b_i=0$ and $N^{2}a=0$; this proves (2). Intersecting over $i$ gives (3).
\end{proof}

\subsection{The Geometric Heisenberg Group}

The $n$ theta groups of Lemma~\ref{lem_K_of_PN} act on the $n$ factors of $[N]_{\calB_{\calA,n}}^{*}\pt_{\calA,n}$ separately. Requiring them to lie over one and the same translation of $\calB_{\calA,n}$ assembles them into a single group.

\begin{definition}\label{def_heisenberg}
Let $N\geq1$ be an integer. The \textbf{geometric Heisenberg group} of level $N$ is
\[
  \calH_N:=\thetagp\bigl([N]_{\calB_{\calA,n}}^{*}\pl^{(1)}\bigr)\times_{\calB_{\calA,n}}\cdots\times_{\calB_{\calA,n}}\thetagp\bigl([N]_{\calB_{\calA,n}}^{*}\pl^{(n)}\bigr),
\]
the theta groups being those of Definition~\ref{def_theta_group} and the fiber products being taken along the maps $\rho$ of \eqref{eq_theta_group}. Concretely, $\calH_N(T)$ is the group of pairs $(c,\tau)$ with $c\in\calB_{\calA,n}(T)$ and $\tau$ an automorphism of the $\gm^{n}$-torsor $([N]_{\calB_{\calA,n}}^{*}\pt_{\calA,n})_T$ lying over the translation $t_c$ of $(\calB_{\calA,n})_T$; that is, $\calH_N$ is the group of automorphisms of the $\gm^{n}$-torsor $[N]_{\calB_{\calA,n}}^{*}\pt_{\calA,n}$ that lie over a translation of $\calB_{\calA,n}$.
\end{definition}

\begin{proposition}\label{prop_heisenberg_exact}
$\calH_N$ is representable by a group scheme over $S$, and there is an exact sequence of group schemes over $S$
\begin{equation}\label{eq_heisenberg_exact}
  1\longrightarrow\gm^{n}\longrightarrow\calH_N\stackrel{\rho}{\longrightarrow}\calB_{\calA,n}[N^{2}]\longrightarrow1,
\end{equation}
in which $\gm^{n}$ is central and $\rho(c,\tau)=c$.
\end{proposition}
\begin{proof}
By Lemma~\ref{lem_theta_representable} each $\thetagp([N]_{\calB_{\calA,n}}^{*}\pl^{(i)})$ is a group scheme and \eqref{eq_theta_group} is an exact sequence of group schemes; take the fiber product of these $n$ sequences over $\calB_{\calA,n}$. The kernel is $\gm^{n}$ and the image is $\bigcap_{i}K([N]_{\calB_{\calA,n}}^{*}\pl^{(i)})$, which is $\calB_{\calA,n}[N^{2}]$ by Lemma~\ref{lem_K_of_PN}(3). Centrality of $\gm^{n}$ is inherited factor by factor.
\end{proof}

\begin{remark}\label{rem_heisenberg_name}
For $n=1$, $\calH_N$ is the usual theta group of $[N]_{\calB_{\calA,1}}^{*}\pl^{(1)}$, with commutator pairing as in \cite[I.3.3.6, V.0]{moret-bailly}; for general $n$ it is the analogue for $\gm^{n}$-extensions.
\end{remark}

\subsection{Torsion Translations}

\begin{definition}\label{def_torsion_multisection}
Let $\calG_{/S}$ be a commutative group scheme. A \textbf{torsion section} is a torsion element $c\in\calG(S)$. A \textbf{torsion multi-section} is a torsion element $c\in\calG(S')$ for some finite surjective morphism $S'\to S$.
\end{definition}

\begin{lemma}\label{lem_torsion_lift}
Let $N\geq1$. Every torsion section $c\in\calB_{\calA,n}[N^{2}](S)$ lifts to a torsion element $\tau\in\calH_N(S')$ for some finite surjective $S'\to S$.
\end{lemma}
\begin{proof}
By Proposition~\ref{prop_heisenberg_exact} the fiber $\rho^{-1}(c)$ is a $\gm^{n}$-torsor over $S$, so it has a point $\tau$ over some finite surjective $S'\to S$. Since $c$ is killed by $N^{2}$ and $\gm^{n}$ is central, $\tau^{N^{2}}=\zeta\in\gm^{n}(S')$. Extracting an $N^{2}$-th root $\eta$ of $\zeta^{-1}$, which again requires only a finite surjective base change, we get $(\eta\tau)^{N^{2}}=1$, and $\eta\tau$ still lies over $c$.
\end{proof}

The level $N$ is part of the datum of a torsion translation, not a quantity fixed once and for all: by Lemma~\ref{lem_torsion_lift} it makes available the translations in $\calB_{\calA,n}[N^{2}]$, so every torsion section of $\calB_{\calA,n}$ in the sense of Definition~\ref{def_torsion_multisection} is reached for $N$ large enough. Given $N\geq1$ and a finite surjective $S'\to S$, write
\[
  \varpi_N:\ \bigl([N]_{\calB_{\calA,n}}^{*}\pt_{\calA,n}\bigr)_{S'}\longrightarrow\pt_{\calA,n}
\]
for the composite of the base change $S'\to S$ with the projection $[N]_{\calB_{\calA,n}}^{*}\pt_{\calA,n}\to\pt_{\calA,n}$; it is finite and surjective, being a base change of the isogeny $[N]_{\calB_{\calA,n}}$. It is along $\varpi_N$ that a subvariety is carried up to the level where $\calH_N$ acts, and then back down.

\begin{definition}\label{def_torsion_translation}
Let $Z\subset\pt_{\calA,n}$ be an irreducible closed subscheme. A \textbf{torsion translate} of $Z$ is a subvariety of the form $\varpi_N(\tau(Z_N))$, where $N$, $S'$ and $\varpi_N$ are as above, $\tau\in\calH_N(S')$ is torsion, and $Z_N$ is an irreducible component of $\varpi_N^{-1}(Z)$; it is irreducible and closed, $\tau$ being an automorphism and $\varpi_N$ being finite. We denote it by $Z+\tau$ when the level $N$ and the component $Z_N$ are understood from the context. Torsion translation of an irreducible closed subscheme of $\calB_{\calA,n}$ is defined in the same way, translating by a torsion multi-section of $\calB_{\calA,n}$.
\end{definition}

\begin{remark}\label{rem_translation_design}
\begin{enumerate}[label=\textup{(\arabic*)}]
  \item The group $\calH_N$ packages all lifts of torsion translations into one action on $[N]_{\calB_{\calA,n}}^{*}\pt_{\calA,n}$. Since $\calB_{\calA,n}[N^{2}]$ is commutative and $\gm^{n}$ acts trivially on $\calB_{\calA,n}$, \eqref{eq_heisenberg_exact} absorbs the order of the translations and the choice of lift.

  \item For the purpose of classifying subvarieties, torsion translation is invertible: $Z$ is an irreducible component of $(Z+\tau)+\tau^{-1}$. Indeed $\varpi_N$ is finite, so $\varpi_N^{-1}(Z+\tau)$ has the same dimension as $Z$, and for dimension reasons one of its irreducible components is $\tau(Z_N)$; applying $\tau^{-1}$ returns $Z_N$, whose image under $\varpi_N$ is $Z$, again for dimension reasons.\qedhere
\end{enumerate}
\end{remark}

\subsection{The Action on Fibers}

We record how torsion translation acts on the fibers of $\pt_{\calA,n}\to(\calAv)^{[n]}$, since this is the form needed in Section~\ref{section7}.

\begin{lemma}\label{lem_noconst}
Let $\pi:\calX\to T$ be an abelian scheme over a variety $T$ over $\CC$, and let $\ttt'$ be a torus over $T$. Then every $T$-morphism $\calX\to\ttt'$ factors through $\pi$, that is, is constant on the fibers of $\pi$. Consequently $\Hom_T(\calX,\ttt')=0$: every homomorphism of group schemes $\calX\to\ttt'$ over $T$ is trivial. For $T=\spec\CC$ the first assertion reads: every morphism from an abelian variety to a torus is constant.
\end{lemma}
\begin{proof}
Being proper and flat with geometrically integral fibers, $\pi$ satisfies $\pi_{*}\calO_\calX=\calO_T$. Both assertions may be checked after an \'etale base change on $T$, so $\ttt'$ may be assumed split; a $T$-morphism into a split torus is the tuple of its composites with the coordinate characters, so we may take $\ttt'=\gm$. A morphism $\calX\to\gm$ is a unit of $\Gamma(\calX,\calO_\calX)=\Gamma(T,\calO_T)$, hence the pullback of a morphism $h:T\to\gm$. If it is a homomorphism then $h=h\circ\pi\circ e_\calX=1$, where $e_\calX$ is the zero section of $\calX$.
\end{proof}

We shall use the following rigidity observation to turn equivariant maps between fibers into homomorphisms, up to translation.

\begin{lemma}\label{lem_rigidity}
Let $A$ be an abelian variety over $\CC$, and let $E_1$ and $E_2$ be semi-abelian varieties over $\CC$ with toric part $\gm^{n}$ and abelian quotient $A$, that is, extensions
\[
  1\longrightarrow\gm^{n}\longrightarrow E_i\longrightarrow A\longrightarrow1\qquad(i=1,2).
\]
Let $f:E_1\to E_2$ be a morphism of varieties which is equivariant for the two $\gm^{n}$-actions and lies over a translation $t_{a}$ of $A$. Then $f=t_{f(0_{E_1})}\circ f_0$ for a unique homomorphism of algebraic groups $f_0:E_1\to E_2$, where $0_{E_1}$ is the identity of $E_1$. If $f$ is an isomorphism then so is $f_0$, and $E_1$ and $E_2$ have the same class in $\Ext^{1}_{\CC}(A,\gm^{n})$.
\end{lemma}
\begin{proof}
Put $f_0:=t_{-f(0_{E_1})}\circ f$. Then $f_0$ is still $\gm^{n}$-equivariant, lies over $\id_A$, and sends $0_{E_1}$ to $0_{E_2}$. It remains to prove that $f_0$ is a homomorphism. Consider its defect
\[
  \kappa:E_1\times E_1\to E_2,\qquad \kappa(p_1,p_2):=f_0(p_1+p_2)-f_0(p_1)-f_0(p_2).
\]
As $f_0$ lies over $\id_A$, the image of $\kappa(p_1,p_2)$ in $A$ is zero, so $\kappa$ takes values in $\gm^{n}\subset E_2$. For $t\in\gm^{n}$, equivariance gives
\[
  \kappa(p_1+t,p_2)=f_0(p_1+p_2)+t-(f_0(p_1)+t)-f_0(p_2)=\kappa(p_1,p_2),
\]
and likewise in the second variable. Thus $\kappa$ descends to a morphism $\bar\kappa:A\times A\to\gm^{n}$.

By Lemma~\ref{lem_noconst}, $\bar\kappa$ is constant. Its value at the origin is zero, because $f_0(0_{E_1})=0_{E_2}$; hence $\kappa\equiv0$ and $f_0$ is a homomorphism. The uniqueness of $f_0$ is immediate from the normalization $f_0(0_{E_1})=0_{E_2}$. If $f$ is an isomorphism, then $f_0$ is an isomorphism of extensions, and the last assertion follows.
\end{proof}

After the finite cover in Definition~\ref{def_torsion_translation}, assume $\tau\in\calH_N(S)$ and write
\[
  \varpi_N:[N]_{\calB_{\calA,n}}^{*}\pt_{\calA,n}\longrightarrow\pt_{\calA,n}.
\]
There are two levels in what follows. The element $\tau$ acts upstairs, on $[N]_{\calB_{\calA,n}}^{*}\pt_{\calA,n}$, and $\varpi_N$ carries the result back downstairs to $\pt_{\calA,n}$.

Fix $s\in S(\CC)$ and put $A=\calA_s$. For $\underline b\in((\calAv)^{[n]})_s(\CC)$ let $E(\underline b)$ be the corresponding fiber of $\pt_{\calA,n}\to(\calAv)^{[n]}$ ``downstairs'', and set
\[
  E_N(\underline b'):=([N]_{\calB_{\calA,n}}^{*}\pt_{\calA,n})|_{A\times\{\underline b'\}},\qquad
  Z_{\underline b}:=Z\cap E(\underline b).
\]

\begin{proposition}\label{prop_fibral_transport}
Let $\tau\in\calH_N(S)$ be a torsion element, lying over $c=(a_0,\underline b_0)\in\calB_{\calA,n}[N^{2}](S)$, and let $Z\subset\pt_{\calA,n}$ be an irreducible closed subvariety. Fix $s$ and $A$ as above, and in \textup{(3)} fix also $\underline b$. Then:
\begin{enumerate}[label=\textup{(\arabic*)}]
  \item for every $\underline b'$, the fiber $E_N(\underline b')$ is a semi-abelian variety with toric part $\gm^{n}$ and abelian quotient $A$, and the restriction $\varpi_{N,\underline b'}:E_N(\underline b')\to E(N\underline b')$ of $\varpi_N$, which lies over $[N]_A$, is a finite surjective homomorphism and an isomorphism on toric parts;
  \item for every $\underline b'$, the automorphism $\tau$ restricts to an isomorphism of $\gm^{n}$-torsors $\tau_{\underline b'}:E_N(\underline b')\to E_N(\underline b'+\underline b_0)$ over $t_{a_0}$, hence by Lemma~\ref{lem_rigidity} is a homomorphism of algebraic groups followed by a translation;
  \item let $Z_N$ be an irreducible component of $\varpi_N^{-1}(Z)$ and let $Z+\tau=\varpi_N(\tau(Z_N))$ be the corresponding torsion translate. Then
  \[
    (Z+\tau)_{\underline b+N\underline b_0}
    \ =\ \bigcup_{\underline b'\,:\,N\underline b'=\underline b}
    \varpi_{N,\underline b'+\underline b_0}\Bigl(\tau_{\underline b'}\bigl((Z_N)_{\underline b'}\bigr)\Bigr),
  \]
  the map $\varpi_{N,\underline b'+\underline b_0}$ of \textup{(1)} having target $E(\underline b+N\underline b_0)$.
\end{enumerate}
\end{proposition}

The maps in the proposition are summarized by
\[
\begin{tikzcd}[row sep=5ex, column sep=10ex]
E_N(\underline b')\arrow[r,"\tau_{\underline b'}","\sim"']\arrow[d,"\varpi_{N,\underline b'}"'] & E_N(\underline b'+\underline b_0)\arrow[d,"\varpi_{N,\underline b'+\underline b_0}"]\\
E(\underline b) & E(\underline b+N\underline b_0)
\end{tikzcd}
\qquad
\begin{tikzcd}[row sep=5ex, column sep=10ex]
A\arrow[r,"t_{a_0}"]\arrow[d,"{[N]_A}"'] & A\arrow[d,"{[N]_A}"]\\
A\arrow[r,"t_{Na_0}"] & A
\end{tikzcd}
\]
where the square on the right records the induced maps on abelian quotients. The left-hand square has no bottom arrow in general; see Remark~\ref{rem_transport_directions}.

\begin{proof}
\textup{(1)} A point $(a,\underline b')$ of the upstairs base maps under $[N]_{\calB_{\calA,n}}$ to $(Na,N\underline b')$. Hence
\[
  E_N(\underline b')=[N]_A^{*}E(N\underline b').
\]
Under this identification $\varpi_{N,\underline b'}$ is the canonical map from the pullback extension to $E(N\underline b')$. It is therefore a homomorphism, is an isomorphism on toric parts, and is finite and surjective because $[N]_A$ is.

\textup{(2)} The element $\tau$ lies over the translation by $c=(a_0,\underline b_0)$ on $\calB_{\calA,n}$. Thus it sends the upstairs fiber over $A\times\{\underline b'\}$ isomorphically to the upstairs fiber over $A\times\{\underline b'+\underline b_0\}$, covers $t_{a_0}$ on $A$, and is $\gm^{n}$-equivariant. By \textup{(1)} these fibers are semi-abelian varieties with toric part $\gm^{n}$ and abelian quotient $A$, so Lemma~\ref{lem_rigidity} applies.

\textup{(3)} Decompose $Z_N$ set-theoretically according to the upstairs dual coordinate:
\[
  Z_N=\bigcup_{\underline b'}(Z_N)_{\underline b'}.
\]
If $z\in (Z_N)_{\underline b'}$, then $\tau(z)$ lies in $E_N(\underline b'+\underline b_0)$ by \textup{(2)}, and applying $\varpi_N$ sends it to the downstairs fiber
\[
  E\bigl(N(\underline b'+\underline b_0)\bigr)=E(N\underline b'+N\underline b_0).
\]
This is the fiber $E(\underline b+N\underline b_0)$ exactly when $N\underline b'=\underline b$. Taking the union over those $\underline b'$ gives the formula.
\end{proof}

The formula in Proposition~\ref{prop_fibral_transport} still involves several upstairs fibers. They all have the same extension class, so they can be canonically identified.

\begin{corollary}\label{cor_fibral_transport}
In the situation of Proposition~\ref{prop_fibral_transport}, fix $\underline b$ as in part \textup{(3)}. Then:
\begin{enumerate}[label=\textup{(\arabic*)}]
  \item every fiber $E_N(\underline b'')$ with $N^{2}\underline b''=N\underline b$ has class $N\underline b$ in $\Ext^{1}_{\CC}(A,\gm^{n})$; hence all such fibers are canonically identified by the unique isomorphisms as extensions of $A$ by $\gm^{n}$. We write $E_N$ for the common semi-abelian variety;
  \item under these identifications, all maps $\varpi_{N,\underline b'+\underline b_0}$ with $N\underline b'=\underline b$ become one $\gm^{n}$-equivariant homomorphism
  \[
    f:\ E_N\longrightarrow E(\underline b+N\underline b_0)
  \]
  lying over $[N]_A$; it is finite, surjective, and an isomorphism on toric parts;
  \item under the same identifications, each $\tau_{\underline b'}$ with $N\underline b'=\underline b$ becomes translation by a point $u_{\underline b'}\in E_N$ lying over $a_0$, and therefore
  \[
    (Z+\tau)_{\underline b+N\underline b_0}
    \ =\ \bigcup_{\underline b'\,:\,N\underline b'=\underline b}
    f\bigl((Z_N)_{\underline b'}+u_{\underline b'}\bigr).
  \]
\end{enumerate}
\end{corollary}
\begin{proof}
\textup{(1)} By Proposition~\ref{prop_fibral_transport}(1),
\[
  E_N(\underline b'')=[N]_A^{*}E(N\underline b'').
\]
Under \eqref{eq_weil_barsotti_split}, pullback by $[N]_A$ multiplies the extension class by $N$, so the class of $E_N(\underline b'')$ is $N^{2}\underline b''$. This is $N\underline b$ by hypothesis. Two extensions of $A$ by $\gm^{n}$ with the same class are isomorphic as extensions; the isomorphism is unique because two isomorphisms of extensions differ by an element of $\Hom_{\CC}(A,\gm^{n})=0$, by Lemma~\ref{lem_noconst}. This proves the canonical identifications.

\textup{(2)} If $N\underline b'=\underline b$, then $N^{2}(\underline b'+\underline b_0)=N\underline b$, since $\underline b_0$ is killed by $N^{2}$. Thus the source of $\varpi_{N,\underline b'+\underline b_0}$ is one of the fibers identified in \textup{(1)}, and its target is
\[
  E\bigl(N(\underline b'+\underline b_0)\bigr)=E(\underline b+N\underline b_0).
\]
Each $\varpi_{N,\underline b'+\underline b_0}$ is $\gm^{n}$-equivariant and lies over $[N]_A$. Such a homomorphism is unique, again because the difference of two of them is a homomorphism $A\to\gm^{n}$. Hence all these maps coincide with a single map $f$, which is finite, surjective, and an isomorphism on toric parts by Proposition~\ref{prop_fibral_transport}(1).

\textup{(3)} By Proposition~\ref{prop_fibral_transport}(2), the map $\tau_{\underline b'}$ is a $\gm^{n}$-equivariant homomorphism over $\id_A$ followed by a translation. After the identifications of \textup{(1)}, that homomorphism must be the identity, by uniqueness. Thus $\tau_{\underline b'}$ becomes translation by $u_{\underline b'}:=\tau_{\underline b'}(0)$, which lies over $a_0$. Substituting this into Proposition~\ref{prop_fibral_transport}(3) gives the displayed formula.
\end{proof}

\begin{corollary}\label{cor_transport_subgroup}
In the situation of Corollary~\ref{cor_fibral_transport}, let $H\subseteq E_N$ be a connected algebraic subgroup such that, for every $\underline b'$ with $N\underline b'=\underline b$, every irreducible component of $(Z_N)_{\underline b'}$ is a translate of $H$. Then every irreducible component of $(Z+\tau)_{\underline b+N\underline b_0}$ is a translate of $f(H)$; the toric part of $f(H)$ is that of $H$, and its abelian quotient is the image under $[N]_A$ of that of $H$.
\end{corollary}
\begin{proof}
By Corollary~\ref{cor_fibral_transport}, the fiber $(Z+\tau)_{\underline b+N\underline b_0}$ is a finite union of sets of the form $f(Y+u_{\underline b'})$, where $Y$ is an irreducible component of $(Z_N)_{\underline b'}$. By hypothesis $Y$ is a translate of $H$, so $f(Y+u_{\underline b'})$ is a translate of the connected algebraic subgroup $f(H)$. Since distinct cosets of $f(H)$ are disjoint, the irreducible components of the finite union are again such translates. The toric part and the abelian quotient of $f(H)$ are as stated because $f$ is $\gm^{n}$-equivariant, an isomorphism on toric parts, and lies over $[N]_A$.
\end{proof}

\begin{remark}\label{rem_transport_directions}
Fiberwise, $Z+\tau=\varpi_N(\tau(Z_N))$ is a correspondence rather than a morphism between downstairs fibers: one passes to the upstairs fibers $E_N(\underline b')$, applies $\tau$, and pushes down by $\varpi_N$. Note also that the shift downstairs is by $\underline b+N\underline b_0$, not by $\underline b+\underline b_0$.
\end{remark}

\section{Subvarieties of Ribet Type}\label{section4}

This section constructs subvarieties of Ribet type and describes their fibers. In toric rank $n$, the construction is determined by a subtorus $\ttt\subseteq\gm^{n}$ and Ribet conditions on the quotient by $\ttt$, without a choice of coordinates on $\gm^{n}$. The identification with special subvarieties of $\pgn$ is established in Sections~\ref{section5} and~\ref{section6}.

Throughout this section, $S$ is an irreducible variety over $\CC$, $\calA\to S$ is an abelian scheme, and all objects are taken over $S$. We denote by $\pr:\calB_{\calA,n}\to(\calAv)^{[n]}$ the projection onto the second factor, and use the same notation for its composite with $\pi:\pt_{\calA,n}\to\calB_{\calA,n}$ and for the restriction of either to any subscheme.

\subsection{Ribet Data}\label{subsubsec_ribet_data}

Let $\calA'\subseteq\calB_{\calA,n}$ be an abelian subscheme over $S$, with components
\[
  \alpha:\calA'\to\calA,\qquad \beta_i:\calA'\to\calAv\quad(1\leq i\leq n),
\]
and write
\[
  \pt_{\calA,n}|_{\calA'}:=\pt_{\calA,n}\times_{\calB_{\calA,n}}\calA' ,
\]
which is a $\gm^{n}$-torsor over $\calA'$. For a character $\chi=(\chi_1,\dots,\chi_n)$ of $\gm^{n}$, identified with an element $\chi\in\ZZ^{n}$, put
\[
  \beta_\chi:=\sum_{i=1}^{n}\chi_i\beta_i\ :\ \calA'\longrightarrow\calAv .
\]
By the bilinearity of the Poincar\'e biextension (Proposition~\ref{prop_poincare_biextension} and Lemma~\ref{lem_bilinear}), the pushout of $\pt_{\calA,n}|_{\calA'}$ along $\chi$ is the $\gm$-torsor associated with
\begin{equation}\label{eq_pushout_along_chi}
  \bigotimes_{i=1}^{n}\bigl((\alpha,\beta_i)^{*}\pl_\calA\bigr)^{\otimes\chi_i}\ \cong\ (\alpha,\beta_\chi)^{*}\pl_\calA .
\end{equation}

\begin{definition}\label{def_ribet_datum}
Let $\ttt\subseteq\gm^{n}$ be a subtorus. We identify the character group $\chl(\gm^{n}/\ttt)$ of $\gm^{n}/\ttt$ with the subgroup of $\chl(\gm^{n})=\ZZ^{n}$ consisting of characters of $\gm^{n}$ that are trivial on $\ttt$. The pair $(\calA',\ttt)$ is a \textbf{Ribet datum} if $(\alpha,\beta_\chi)$ satisfies the Ribet property for every $\chi\in\chl(\gm^{n}/\ttt)$. These conditions are called the \textbf{Ribet conditions} of $(\calA',\ttt)$.
\end{definition}

Since $\beta\mapsto\alpha^{\vee}\beta+\beta^{\vee}\alpha$ is additive in $\beta$, it is enough to verify the Ribet conditions on a basis of the free abelian group $\chl(\gm^{n}/\ttt)$.

\subsection{Ribet Sections}\label{subsubsec_ribet_sections}

Fix a Ribet datum $(\calA',\ttt)$, and write
\[
  q_\ttt:\gm^{n}\longrightarrow\gm^{n}/\ttt,\qquad
  \pt_{\calA,n}\big|_{\calA',\gm^{n}/\ttt}:=(q_\ttt)_{*}\bigl(\pt_{\calA,n}|_{\calA'}\bigr),
\]
for the quotient homomorphism and the pushout of $\pt_{\calA,n}|_{\calA'}$ along it, respectively. The latter is a $\gm^{n}/\ttt$-torsor over $\calA'$; the two indices after the restriction symbol specify its base and structure group. Consider the cartesian square
\[
\begin{tikzcd}[column sep=9ex, row sep=5ex]
{[2]_{\calA'}^{*}\bigl(\pt_{\calA,n}|_{\calA',\gm^{n}/\ttt}\bigr)} \arrow[d] \arrow[r,"\varpi"]\arrow[dr, phantom,"\ulcorner",very near start]
& \pt_{\calA,n}|_{\calA',\gm^{n}/\ttt} \arrow[d]\\
\calA' \arrow[r,"{[2]_{\calA'}}"] \arrow[u,"\tilde r",dotted,bend left]
& \calA' .
\end{tikzcd}
\]

\begin{prop-def}\label{def_ribet_section}
With the notation above:
\begin{enumerate}[label=\textup{(\arabic*)}]
  \item $\pt_{\calA,n}|_{\calA',\gm^{n}/\ttt}$ is an extension of $\calA'$ by $\gm^{n}/\ttt$; under the isomorphism
  \[
    \Ext^{1}_{S}\bigl(\calA',\gm^{n}/\ttt\bigr)\ \cong\ \Hom_{\ZZ}\bigl(\chl(\gm^{n}/\ttt),\calA'^{\vee}(S)\bigr)
  \]
  of Theorem~\ref{thm_weil_barsotti}, its class is $\chi\mapsto\bigl[(\alpha,\beta_\chi)^{*}\pl_\calA\bigr]$ and is killed by $2$;
  \item $[2]_{\calA'}^{*}\bigl(\pt_{\calA,n}|_{\calA',\gm^{n}/\ttt}\bigr)$ is therefore the trivial torsor, and there is a unique trivialization compatible with the bi-rigidification of $\pl_\calA$;
  \item this trivialization defines a section $\tilde r$ of $[2]_{\calA'}^{*}\bigl(\pt_{\calA,n}|_{\calA',\gm^{n}/\ttt}\bigr)\to\calA'$. The image of $r:=\varpi\circ\tilde r$ is an irreducible closed subscheme of $\pt_{\calA,n}|_{\calA',\gm^{n}/\ttt}$, finite and surjective over $\calA'$, called the \textbf{Ribet multi-section} of $(\calA',\ttt)$. The morphism $r$ lies over $[2]_{\calA'}$, rather than over $\id_{\calA'}$.
\end{enumerate}
\end{prop-def}
\begin{proof}
(1) By \eqref{eq_pushout_along_chi}, the pushout of $\pt_{\calA,n}|_{\calA',\gm^{n}/\ttt}$ along $\chi\in\chl(\gm^{n}/\ttt)$ is the $\gm$-torsor associated with $(\alpha,\beta_\chi)^{*}\pl_\calA$. The Ribet conditions and Theorem~\ref{thm_ribet_property}(1) imply that its class lies in $\pic^{0}_{\calA'/S}(S)=\calA'^{\vee}(S)$, where it is killed by $2$ by Theorem~\ref{thm_ribet_property}(2). A $\gm^{n}/\ttt$-torsor whose pushouts along all characters lie in $\pic^{0}$ is an extension. Theorem~\ref{thm_weil_barsotti} identifies its class with these pushouts.

(2) Under the isomorphism in (1), pullback along $[2]_{\calA'}$ acts on $\Ext^{1}_{S}(\calA',\gm^{n}/\ttt)$ by composition with $[2]_{\calA'}^{\vee}=[2]_{\calA'^{\vee}}$, hence by multiplication by $2$. Since the class is killed by $2$ by (1), the pullback is trivial. The bi-rigidification of $\pl_\calA$ gives uniqueness of the trivialization.

(3) The trivialization gives the section and its image as stated. The finiteness and surjectivity follow from those of $[2]_{\calA'}$.
\end{proof}

\begin{remark}
For $n=1$ and $\ttt=1$ this construction is the canonical splitting of \cite[p.~5, condition (i)]{Bertrand98}, spelled out for the Ribet point in \cite[\S1]{Bertrand11}: the square of the restriction of $\pl_\calA$ to the graph of an antisymmetric $\varphi$ is trivial, and a unique nonzero section is obtained after a $2$-isogeny. The base change by $[2]_{\calA'}$ is necessary in general: Theorem~\ref{thm_ribet_property} shows that the class of $\pt_{\calA,n}|_{\calA',\gm^{n}/\ttt}$ is killed by $2$, but does not assert that it vanishes. In the case $n=1$, $\ttt=1$ and $\calA'=\Gamma(f-f^{\vee})$ for a homomorphism $f:\calA\to\calAv$, the authors of \cite{BE20} prove that $(\alpha,\beta)^{*}\pl_\calA$ is itself trivial. Their Ribet section is therefore a genuine section; see also Example~\ref{ex_ribet_graph}.
\end{remark}

\subsection{Subvarieties of Ribet Type}\label{subsubsec_ribet_type}

We now associate a subvariety of the Poincar\'e torsor with each Ribet datum. The construction consists of taking the inverse image of the section on the pullback by multiplication by $2$, and then mapping it to the original torsor.

Write $p_\ttt$ for the pushout morphism along $q_\ttt$,
\[
p_\ttt:\ [2]_{\calA'}^{*}\bigl(\pt_{\calA,n}|_{\calA'}\bigr)\longrightarrow[2]_{\calA'}^{*}\bigl(\pt_{\calA,n}|_{\calA',\gm^{n}/\ttt}\bigr) .
\]
The inverse image of the section $\tilde r$ under this morphism is defined by the cartesian square
\[
\begin{tikzcd}[column sep=12ex, row sep=5ex]
\widetilde R(\calA',\ttt) \arrow[r,hook] \arrow[d] \arrow[dr,phantom,"\ulcorner",very near start]
& {[2]_{\calA'}^{*}\bigl(\pt_{\calA,n}|_{\calA'}\bigr)} \arrow[d,"p_\ttt"]\\
\calA' \arrow[r,"\tilde r"']
& {[2]_{\calA'}^{*}\bigl(\pt_{\calA,n}|_{\calA',\gm^{n}/\ttt}\bigr)}\ \rlap{.}
\end{tikzcd}
\]
Thus $\widetilde R(\calA',\ttt)$ is a $\ttt$-torsor over the abelian scheme $\calA'$, with toric part $\ttt$. In particular, it is irreducible. For $\ttt=\gm^{n}$, one has $\pt_{\calA,n}|_{\calA',\gm^{n}/\ttt}=0$, and $\widetilde R(\calA',\gm^{n})$ is the whole torsor $[2]_{\calA'}^{*}\bigl(\pt_{\calA,n}|_{\calA'}\bigr)$. For $\ttt=1$, the toric part is trivial and $\widetilde R(\calA',1)$ is the graph of $\tilde r$.

We then take the image of $\widetilde R(\calA',\ttt)$ under the composite of the two morphisms in the top row of
\[
\begin{tikzcd}[column sep=9ex, row sep=4.5ex]
{[2]_{\calA'}^{*}\bigl(\pt_{\calA,n}|_{\calA'}\bigr)} \arrow[r,"\varpi"] \arrow[d]
& \pt_{\calA,n}|_{\calA'} \arrow[r] \arrow[d]
& \pt_{\calA,n} \arrow[d]\\
\calA' \arrow[r,"{[2]_{\calA'}}"]
& \calA' \arrow[r,hook]
& \calB_{\calA,n}\ \rlap{,}
\end{tikzcd}
\]
both squares being cartesian. Here $\varpi$ is the canonical map from the pullback by $[2]_{\calA'}$, and the second morphism is induced by the inclusion $\calA'\hookrightarrow\calB_{\calA,n}$.

\begin{definition}\label{def_ribet_type}
Denote this image by $R(\calA',\ttt)\subseteq\pt_{\calA,n}$; it is irreducible and closed. As in Subsection~\ref{subsubsec_ribet_sections}, a tilde denotes an object on the pullback by $[2]_{\calA'}$, and the corresponding symbol without a tilde denotes its image in the original torsor. A \textbf{subvariety of Ribet type} of $\pt_{\calA,n}$ is a torsion translate, in the sense of Definition~\ref{def_torsion_translation}, of a subvariety $R(\calA',\ttt)$ associated with a Ribet datum $(\calA',\ttt)$ over an irreducible closed subvariety $T\subseteq S$. The construction is then carried out with $S$ replaced by $T$. We denote such a subvariety by $R$ and its base by $T$.
\end{definition}

To extend this definition to semi-abelian schemes, let $\calG\to S$ be a semi-abelian scheme with split toric part $\gm^{n}$ and abelian quotient $\calA$, and let $\underline b\in(\calAv)^{[n]}(S)$ be its class. Then
\[
  \calA\times\{\underline b\}:=\bigl\{(a,\underline b)\ :\ a\in\calA\bigr\}\ \subseteq\ \calB_{\calA,n}
\]
is a subscheme of $\calB_{\calA,n}$ over $S$. By Proposition~\ref{prop_universal_semiabelian}(3), the scheme $\calG$ is identified with the restriction $\pt_{\calA,n}|_{\calA\times\{\underline b\}}$, in the notation of Subsection~\ref{subsubsec_ribet_data}. Under this identification, subvarieties of $\calG$ are subvarieties of $\pt_{\calA,n}$ lying over $\calA\times\{\underline b\}$.

\begin{definition}\label{def_ribet_type_general}
Let $\calG\to S$ be a semi-abelian scheme with split toric part $\gm^{n}$ and abelian quotient $\calA$, identified with $\pt_{\calA,n}|_{\calA\times\{\underline b\}}$ as above. An irreducible closed subvariety of $\calG$ is \textbf{of Ribet type} if it is an irreducible component of
\[
  R\ \cap\ \pt_{\calA,n}\big|_{\calA\times\{\underline b\}}
\]
for some subvariety $R$ of Ribet type of $\pt_{\calA,n}$.
\end{definition}

\begin{remark}\label{rem_ribet_type_extremes}
The class of subvarieties of Ribet type is intended to include subgroup schemes over arbitrary irreducible closed subvarieties of the base. This inclusion is not immediate from Definition~\ref{def_ribet_type}; it is established in Proposition~\ref{prop_flat_subgroup} below. For $n=0$, that is, $\calG\to S$ is an abelian scheme, every subvariety of Ribet type is a torsion translate of an abelian subscheme over an irreducible closed subvariety of $S$.

We reserve the term \emph{special} for the mixed Shimura setting; see Sections~\ref{section5} and~\ref{section6}, where we relate these two notions.
\end{remark}

\begin{proposition}\label{prop_flat_subgroup}
Let $\calG\to S$ be a semi-abelian scheme with split toric part $\gm^{n}$ and abelian quotient $\calA$, let $T\subseteq S$ be an irreducible closed subvariety, and let $\calE\subseteq\calG|_T$ be an irreducible closed subgroup scheme, flat over $T$. Then $\calE$ is of Ribet type in the sense of Definition~\ref{def_ribet_type_general}, and the Ribet conditions occurring in it are trivial.
\end{proposition}
\begin{proof}
Since Definition~\ref{def_ribet_type} allows Ribet data over any irreducible closed subvariety of $S$, we may replace $S$ by $T$ and assume that $T=S$. Write $\underline b\in(\calAv)^{[n]}(S)$ for the class of $\calG$ and $b_\chi=\sum_i\chi_i b_i$. We identify $\calG$ with $\pt_{\calA,n}|_{\calA\times\{\underline b\}}$ and construct a subvariety $R$ of Ribet type such that $\calE=R\cap\calG$.

\noindent\textit{Step 1: the subgroup data.}
Flatness and characteristic zero imply that $\calE\to S$ is smooth. For every integer $r\geq1$, multiplication $[r]_{\calE}$ is finite and \'etale: finiteness follows from that of $[r]_{\calG}$, and \'etaleness follows because its differential is multiplication by $r$. Its image is therefore open and closed in the irreducible scheme $\calE$, hence is all of $\calE$. Thus multiplication by every positive integer is surjective on each fiber of $\calE$, so its finite group of connected components is trivial. Consequently $\calE$ is a semi-abelian scheme. Its toric rank is constant, since both its relative dimension and the degree of $[2]_{\calE}$ are constant.

Let $\ttt=(\calE\cap\gm^{n})^{\circ}$ be its toric part, a constant subtorus of $\gm^{n}$, and let $\mathcal D=\calE/\ttt$ be its abelian quotient. The induced map to $\calA$ factors as
\[
  \mathcal D\stackrel{f}{\longrightarrow}\calX
  \stackrel{\iota}{\hookrightarrow}\calA,
\]
where $\calX$ is the image of $\calE$ in $\calA$ and $f$ is an isogeny. The kernel of $f$ is $(\calE\cap\gm^{n})/\ttt$, which need not be trivial.

Put $Q=\gm^{n}/\ttt$ and $L=\chl(Q)$. The inclusion $\calE\hookrightarrow\calG$ induces a homomorphism $\mathcal D\to(q_\ttt)_*(\calG|_{\calX})$ over $f$. Hence the pullback of this extension along $f$ splits. By Theorem~\ref{thm_weil_barsotti},
\begin{equation}\label{eq_flat_subgroup_class}
  f^{\vee}\iota^{\vee}b_\chi=0
  \qquad\text{for every }\chi\in L.
\end{equation}
Choose $m\geq1$ annihilating $\ker f^{\vee}$. Then $m\iota^{\vee}b_\chi=0$ for every $\chi\in L$.

\noindent\textit{Step 2: a Ribet datum with trivial conditions.}
Choose a basis $\chi_1,\dots,\chi_k$ of $L$ and put
\[
  \calY:=\ker\!\left((\calAv)^{[n]}\longrightarrow(\calX^{\vee})^k,
  \quad\underline c\longmapsto
  (\iota^{\vee}c_{\chi_1},\dots,\iota^{\vee}c_{\chi_k})\right).
\]
This kernel is connected. Indeed, $L$ is a saturated sublattice of $\ZZ^n$, so its chosen basis extends to a basis of $\ZZ^n$. In the resulting coordinates, $\calY$ is isomorphic to
\[
  (\ker\iota^{\vee})^k\times_S(\calAv)^{[n-k]}.
\]
Here $\ker\iota^{\vee}=(\calA/\calX)^{\vee}$ is an abelian scheme. Thus $\calY$ is an abelian subscheme, and $m\underline b\in\calY(S)$. After a finite surjective base change, choose $\underline b_0\in\calY$ with $m\underline b_0=m\underline b$. Then
\[
  \underline b=\underline b_0+\sigma,\qquad m\sigma=0.
\]
All finite covers below are taken as in Definition~\ref{def_torsion_translation}; we work on an irreducible component dominating $S$ and suppress the base change from the notation.

Set $\calA'=\calX\times_S\calY\subseteq\calB_{\calA,n}$. For every $\chi\in L$, the line bundle $(\alpha,\beta_\chi)^*\pl_\calA$ restricts trivially to each $\calX\times\{\underline c\}$, since $\iota^{\vee}c_\chi=0$, and to $\{0_{\calX}\}\times\calY$ by rigidification. The seesaw theorem makes it trivial on $\calA'$. Hence $(\calA',\ttt)$ is a Ribet datum with trivial Ribet conditions.

Write $R_0=R(\calA',\ttt)$. The quotient torsor $\pt_{\calA,n}|_{\calA',Q}$ has a trivialization normalized at the origin, whose section we denote by $s_0$. Its pullback along $[2]_{\calA'}$ is the section $\tilde r$ in Proposition-Definition~\ref{def_ribet_section}. Thus $R_0=p_\ttt^{-1}(s_0(\calA'))$, where $p_\ttt:\pt_{\calA,n}|_{\calA'}\to\pt_{\calA,n}|_{\calA',Q}$ is the quotient map. In particular, $R_0$ is a $\ttt$-torsor over $\calA'$. The restriction of $s_0$ to $\{0_{\calX}\}\times\calY$ agrees with the rigidification, by Lemma~\ref{lem_noconst}. For each $\underline c\in\calY$, its restriction to $\calX\times\{\underline c\}$ is therefore the homomorphic splitting. Consequently $(R_0)_{\underline c}$ is a connected subgroup of $E(\underline c)$, extending $\calX$ by $\ttt$. This group structure is relative to $\calY$.

\noindent\textit{Step 3: a torsion lift preserving the identities.}
Choose $N$ divisible by $m$ and, after a further finite surjective base change, choose $\sigma_0$ with $N\sigma_0=\sigma$. Then $N^2\sigma_0=0$. Lemma~\ref{lem_torsion_lift} supplies a torsion lift $\tau\in\calH_N$ of $(0,\sigma_0)$.

The rigidification identifies the restriction of $[N]_{\calB_{\calA,n}}^*\pt_{\calA,n}$ to $\{0_{\calA}\}\times(\calAv)^{[n]}$ with the trivial $\gm^{n}$-torsor. On this restriction the action of $\tau$ has the form
\[
  (\underline c,z)\longmapsto
  (\underline c+\sigma_0,u(\underline c)z).
\]
The morphism $u:(\calAv)^{[n]}\to\gm^{n}$ factors through $S$ by Lemma~\ref{lem_noconst}, so we write it as $u\in\gm^{n}(S)$. If $\tau$ has order $\ell$, then $u^{\ell}=1$. Since $\gm^{n}$ is central in $\calH_N$, replacing $\tau$ by $u^{-1}\tau$ preserves its torsion property and makes its action preserve the identity in every fiber. By Proposition~\ref{prop_fibral_transport}\textup{(2)}, each resulting map
\[
  \tau_{\underline c}:E_N(\underline c)
  \longrightarrow E_N(\underline c+\sigma_0)
\]
is a homomorphism.

\noindent\textit{Step 4: the entire intersection.}
Let $Z_N$ be the restriction of $\varpi_N^{-1}(R_0)$ to the component $\calA'$ of $[N]_{\calB_{\calA,n}}^{-1}(\calA')$. It is a $\ttt$-torsor over $\calA'$, hence an irreducible component of $\varpi_N^{-1}(R_0)$. For every $\underline c\in\calY$, its fiber $(Z_N)_{\underline c}$ is the pullback of $(R_0)_{N\underline c}$ along $[N]_{\calX}$, so it is a connected subgroup of $E_N(\underline c)$ extending $\calX$ by $\ttt$. Define
\[
  R:=\varpi_N\bigl(\tau(Z_N)\bigr),
\]
where $\varpi_N$ includes the finite map back to the original base. This is a subvariety of Ribet type, lying over $\calA'+(0,\sigma)$.

Fix a geometric point of the chosen base cover and use the fiber notation of Proposition~\ref{prop_fibral_transport}. For every $\underline c\in\calY$ satisfying $N\underline c=\underline b_0$, put
\[
  H_{\underline c}:=
  \varpi_{N,\underline c+\sigma_0}
  \bigl(\tau_{\underline c}((Z_N)_{\underline c})\bigr)
  \ \subseteq\ E(\underline b).
\]
Both maps are homomorphisms, and $\varpi_{N,\underline c+\sigma_0}$ is finite and an isomorphism on toric parts. Thus $H_{\underline c}$ is a connected subgroup with toric part $\ttt$ and image $\calX$ in $\calA$.

We compare it with the fiber of $\calE$. In the quotient of $\calG|_{\calX}$ by $\ttt$, the images of these two subgroups are abelian subvarieties mapping isogenously onto $\calX$. Their sum is again an abelian subvariety. Its intersection with the torus $Q$ is finite, so its dimension is at most $\dim\calX$. Each summand already has that dimension; hence the two images coincide. Taking their preimages under the quotient by $\ttt$ gives $H_{\underline c}=\calE$ on this fiber.

Proposition~\ref{prop_fibral_transport}\textup{(3)} now gives
\[
  R\cap E(\underline b)
  =\bigcup_{\substack{\underline c\in\calY\\N\underline c=\underline b_0}}
  H_{\underline c}
  =\calE.
\]
The indexing set is nonempty because $[N]_{\calY}$ is surjective. This equality holds on every geometric fiber and for every point above the original base. Taking the finite image back to $S$ therefore yields $R\cap\calG=\calE$ as closed subvarieties, as required.
\end{proof}

\subsection{Fibers of Subvarieties of Ribet Type}

We now describe the fibers of a subvariety of Ribet type over $(\calAv)^{[n]}$. We first establish the group-theoretic lemma used in the proof.

\begin{lemma}\label{lem_elementary_fiber}
Let $A$ be an abelian variety over $\CC$, and let $E$ be a semi-abelian variety over $\CC$ with toric part $\gm^{n}$ and abelian quotient $A$, that is, an extension
\[
  1\longrightarrow\gm^{n}\longrightarrow E\stackrel{\pi_E}{\longrightarrow}A\longrightarrow1 .
\]
Let $\ttt'\subseteq\gm^{n}$ be a subtorus, let $q:\gm^{n}\to\gm^{n}/\ttt'$ be the quotient homomorphism, and let $p:E\to q_{*}E$ be the morphism associated with the pushout of $E$ along $q$. The semi-abelian variety $q_{*}E$ has toric part $\gm^{n}/\ttt'$ and abelian quotient $A$.
\begin{enumerate}[label=\textup{(\arabic*)}]
  \item For every abelian subvariety $X\subseteq A$ there is at most one algebraic subgroup of $E$ with toric part $\ttt'$ and abelian quotient $X$.
  \item Let $K\subseteq A$ be an algebraic subgroup, let $C\subseteq A$ be a coset of $K$, and let $\varsigma$ be a section of $q_{*}E\to A$ over $C$. Put
  \[
  \Sigma:=\pi_E^{-1}(C)\cap p^{-1}\bigl(\varsigma(C)\bigr).
  \]
  Then $E$ has an algebraic subgroup with toric part $\ttt'$ and abelian quotient $K^{\circ}$, and every irreducible component of $\Sigma$ is a translate of this subgroup.
\end{enumerate}
\end{lemma}
\begin{proof}
\textup{(1)} Let $H$ and $H'$ be algebraic subgroups of $E$ with toric part $\ttt'$ and abelian quotient $X$. Both lie in $\pi_E^{-1}(X)$. Their images in the pushout of $\pi_E^{-1}(X)$ along $q$ are homomorphic sections of an extension of $X$ by $\gm^{n}/\ttt'$. The difference of two such sections belongs to $\Hom_{\CC}(X,\gm^{n}/\ttt')$, which vanishes by Lemma~\ref{lem_noconst} because $X$ is an abelian variety. The images therefore coincide, and $H=H'$ since each is the preimage of its image.

\textup{(2)} Let $K^{\circ}$ be the identity component of $K$, an abelian subvariety of $A$. The coset $C$ is a finite union of cosets of $K^{\circ}$, so $\Sigma$ is the union of the corresponding subsets $\Sigma\cap\pi_E^{-1}(K^{\circ}+c_1)$. It suffices to consider one such subset. Choose $\xi_1\in q_{*}E$ with $\pi(\xi_1)=c_1$, where $\pi:q_{*}E\to A$ is the structure homomorphism. Translation by $-\xi_1$ identifies $\pi^{-1}(K^{\circ}+c_1)$ with $\pi^{-1}(K^{\circ})$ as $(\gm^{n}/\ttt')$-torsors, and maps $\varsigma$ to a section of $\pi^{-1}(K^{\circ})\to K^{\circ}$. Hence this torsor is trivial.

The semi-abelian variety $\pi^{-1}(K^{\circ})$ has toric part $\gm^{n}/\ttt'$ and abelian quotient $K^{\circ}$. By Theorem~\ref{thm_weil_barsotti}, its class in $\Ext^{1}_{\CC}(K^{\circ},\gm^{n}/\ttt')$ is the class of the underlying torsor. The extension therefore splits, and we fix an isomorphism of algebraic groups $\pi^{-1}(K^{\circ})\cong K^{\circ}\times(\gm^{n}/\ttt')$.

Under this isomorphism, the section is the graph of a morphism $K^{\circ}\to\gm^{n}/\ttt'$. This morphism is constant by Lemma~\ref{lem_noconst}, since $K^{\circ}$ is an abelian variety. Let $\Xi\subseteq q_{*}E$ be the algebraic subgroup corresponding to $K^{\circ}\times\{1\}$. It has trivial toric part and abelian quotient $K^{\circ}$, and there is a point $\eta\in q_{*}E$ such that
\[
  \varsigma(K^{\circ}+c_1)=\Xi+\eta .
\]
Put $H:=p^{-1}(\Xi)$, an algebraic subgroup of $E$, and choose $\xi\in p^{-1}(\eta)$. Write $C$ as the disjoint union of its cosets $K^{\circ}+c_j$. Since $\varsigma(C)=\bigcup_j\varsigma(K^{\circ}+c_j)$, we have
\[
  \Sigma\cap\pi_E^{-1}(K^{\circ}+c_1)
  \ =\ \bigcup_j\Bigl[\pi_E^{-1}(K^{\circ}+c_1)\cap p^{-1}\bigl(\varsigma(K^{\circ}+c_j)\bigr)\Bigr].
\]
For $j\neq1$, the $j$-th term is empty: the set $p^{-1}\bigl(\varsigma(K^{\circ}+c_j)\bigr)$ lies over $K^{\circ}+c_j$, and the cosets are disjoint. The remaining term is $\pi_E^{-1}(K^{\circ}+c_1)\cap p^{-1}(\Xi+\eta)$. Since $p^{-1}(\Xi+\eta)$ already lies over $K^{\circ}+c_1$, this gives
\[
  \Sigma\cap\pi_E^{-1}(K^{\circ}+c_1)\ =\ p^{-1}(\Xi+\eta)\ =\ H+\xi .
\]
The subgroup $H$ contains $\ker q=\ttt'$ and maps onto $K^{\circ}$. Its toric part is therefore $\ttt'$ and its abelian quotient is $K^{\circ}$. Since $\ttt'$ and $K^{\circ}$ are connected, so is $H$; hence $H+\xi$ is irreducible.
\end{proof}

Let $R$ be a subvariety of Ribet type. Replacing $S$ by the irreducible base of its Ribet datum, we assume that the datum is defined over all of $S$.

For $\underline b\in(\calAv)^{[n]}(\CC)$, let $s$ be its image in $S$, put $A:=\calA_s$, and denote by $E(\underline b)$ the corresponding fiber of $\pt_{\calA,n}\to(\calAv)^{[n]}$. The translations in the following theorem need not be torsion translations, as shown by the example in \cite{BE20}.

\begin{theorem}\label{thm_ribet_type_fibers}
Let $R\subseteq\pt_{\calA,n}$ be a subvariety of Ribet type, with Ribet datum $(\calA',\ttt)$, and let $\underline b\in(\calAv)^{[n]}(\CC)$, with image $s\in S$. Put $K:=\calA'_s\cap\bigl(\calA_s\times\{0\}\bigr)$. Then there is exactly one algebraic subgroup $H\subseteq E(\underline b)$ with toric part $\ttt$ and abelian quotient $K^{\circ}$, and every irreducible component of the fiber $R_{\underline b}$ is a translate of it.
\end{theorem}

\begin{proof}
By Definition~\ref{def_ribet_type}, write $R=R(\calA',\ttt)+\tau$ with $\tau$ torsion. After a base change, we may assume that $\tau\in\calH_N(S)$, so the morphism in Definition~\ref{def_torsion_translation} is $\varpi_N:[N]_{\calB_{\calA,n}}^{*}\pt_{\calA,n}\to\pt_{\calA,n}$. Let $R_N$ be the irreducible component of $\varpi_N^{-1}\bigl(R(\calA',\ttt)\bigr)$ used in that definition, and let $c=(a_0,\underline b_0)\in\calB_{\calA,n}[N^{2}](S)$ be the point below $\tau$.

Applying Corollary~\ref{cor_transport_subgroup} to $Z:=R(\calA',\ttt)$ and to the fiber over $\underline b-N\underline b_0$, we reduce to the following assertion:
\begin{equation}\label{eq_fiber_to_prove}
  \parbox{0.86\textwidth}{there is a connected algebraic subgroup $H_N\subseteq E_N$, with toric part $\ttt$ and abelian quotient an abelian subvariety $X\subseteq A$ satisfying $[N]_A(X)=K^{\circ}$, of which every irreducible component of every $(R_N)_{\underline b'}$ is a translate}
\end{equation}
Assuming \eqref{eq_fiber_to_prove}, Corollary~\ref{cor_transport_subgroup} implies that every irreducible component of $R_{\underline b}$ is a translate of $f(H_N)$, which has toric part $\ttt$ and abelian quotient $[N]_A(X)=K^{\circ}$; we may therefore take $H:=f(H_N)$, and uniqueness follows from Lemma~\ref{lem_elementary_fiber}\textup{(1)}.

Fix $\underline b'$ with $N\underline b'=\underline b-N\underline b_0$. This point lies over $s$, and Proposition~\ref{prop_fibral_transport}(1) identifies the fiber $E_N(\underline b')$ as a semi-abelian variety with toric part $\gm^{n}$ and abelian quotient $A$.

\smallskip
\noindent\textit{Step 1: describing $R_N$.} Since $\varpi_N$ is an isomorphism on the fibers of the two $\gm^{n}$-torsors and induces $[N]_{\calB_{\calA,n}}$ on the base, $\varpi_N^{-1}\bigl(R(\calA',\ttt)\bigr)$ lies over $[N]_{\calB_{\calA,n}}^{-1}(\calA')$. The latter is a closed subgroup scheme containing $\calA'$ and finite flat over it, hence of the same relative dimension; so its identity component is $\calA'$, and its irreducible components are the cosets $\calA'+t$ with $t$ a torsion section of $\calB_{\calA,n}$, defined over $S$ after a further finite base change. Let $\calA'+\delta$ be the coset below $R_N$, and choose, after one more base change, a section $e$ of $\calA'$ with $2e=[N]\delta$; this is possible because $[2]_{\calA'}$ is surjective.

Put $\calP_\delta:=\bigl([2]_{\calA'}+\delta\bigr)^{*}[N]_{\calB_{\calA,n}}^{*}\pt_{\calA,n}$ and let $\varpi_\delta:\calP_\delta\to[N]_{\calB_{\calA,n}}^{*}\pt_{\calA,n}$ be the canonical map from this pullback, a finite morphism lying over $[2]_{\calA'}+\delta$. Since $[N](2b+\delta)=2\phi(b)$ for $\phi(b):=[N]b+e$, the fiber of $\calP_\delta$ at $b$ is the fiber of $\pt_{\calA,n}$ at $2\phi(b)$, that is, the fiber of $[2]_{\calA'}^{*}\bigl(\pt_{\calA,n}|_{\calA'}\bigr)$ at $\phi(b)$. So $\varsigma:=\tilde r\circ\phi$ is a section of the pushout of $\calP_\delta$ along $q_\ttt$, and
\[
  \widetilde R_N:=p_\ttt^{-1}\bigl(\varsigma(\calA')\bigr)
\]
is a $\ttt$-torsor over $\calA'$, hence irreducible, with $\varpi_\delta(\widetilde R_N)$ irreducible and closed. We claim that $R_N=\varpi_\delta(\widetilde R_N)$ for a suitable choice of $e$.

Indeed, $R(\calA',\ttt)$ is by Definition~\ref{def_ribet_type} the image of $\widetilde R(\calA',\ttt)=p_\ttt^{-1}\bigl(\tilde r(\calA')\bigr)$ under the canonical map $\varpi$, which is an isomorphism on fibers over $[2]_{\calA'}$; hence
\[
  R(\calA',\ttt)_a=\bigcup_{a_1\in\calA',\ 2a_1=a}p_\ttt^{-1}\bigl(\tilde r(a_1)\bigr),
\]
a finite union of $\ttt$-cosets. As $\varpi_N$ is likewise an isomorphism on fibers, the fiber of $\varpi_N^{-1}\bigl(R(\calA',\ttt)\bigr)$ at a point $b$ of $\calA'+\delta$ is the corresponding union over the $a_1$ with $2a_1=[N]b$. Writing $b=2b_0+\delta$, such $a_1$ are exactly the $[N]b_0+e'$ with $2e'=[N]\delta$; these $e'$ form a torsor under $\calA'[2]$, so there are finitely many, and writing $\varsigma_{e'}$ for the section attached to $e'$ as $\varsigma$ was to $e$ we get
\[
  \varpi_N^{-1}\bigl(R(\calA',\ttt)\bigr)\big|_{\calA'+\delta}
  =\bigcup_{e'}\ \varpi_\delta\Bigl(p_\ttt^{-1}\bigl(\varsigma_{e'}(\calA')\bigr)\Bigr),
\]
a finite union of irreducible closed sets. The component $R_N$ is contained in one of them and, being a component, equals it; taking $e:=e'$ proves the claim. For $R(\calA',\ttt)$ itself one has $N=1$, $\delta=0$, $e=0$ and $\varsigma=\tilde r$, recovering the square of Subsection~\ref{subsubsec_ribet_type}.

\smallskip
\noindent\textit{Step 2: restricting to the fiber over $\underline b'$.} Put $A':=\calA'_s$, a subvariety of $A\times(A^{\vee})^{n}$, the fiber of $\calB_{\calA,n}$ at $s$, write its points as pairs $(a,\underline b'')$, and write $\delta_s=(\delta_0,\underline\delta)$ for the value of $\delta$ at $s$. Since $E_N(\underline b')$ lies over $s$, its preimage in $\widetilde R_N$ lies over $A'$; by Step 1 a point over $(a,\underline b'')$ has image under $\varpi_\delta$ lying over $(2a+\delta_0,\,2\underline b''+\underline\delta)$, which is in $E_N(\underline b')$ precisely when $2\underline b''+\underline\delta=\underline b'$. Therefore
\[
  (R_N)_{\underline b'}=\bigcup_{\underline b''}\,\varpi_\delta\bigl(\widetilde R_N|_{C_{\underline b''}}\bigr),
  \qquad C_{\underline b''}:=A'\cap\bigl(A\times\{\underline b''\}\bigr),
\]
a finite union over the $\underline b''\in(A^{\vee})^{n}$ with $2\underline b''+\underline\delta=\underline b'$ and $C_{\underline b''}\neq\emptyset$. It suffices to treat one term, so fix such a $\underline b''$ and put $C:=C_{\underline b''}$, a coset of $A'\cap(A\times\{0\})=K$. The second coordinate being constant on $C$, the projection $(a,\underline b'')\mapsto a$ is injective on it; we identify $C$ with its image $C_0\subseteq A$, a coset of $K$ viewed in $A$ by the same projection.

\smallskip
\noindent\textit{Step 3: identifying $\widetilde R_N|_{C}$ inside $[2]_A^{*}E_N(\underline b')$.} By Step 2, $\varpi_\delta$ maps the fiber of $\widetilde R_N$ at a point of $C$ with first coordinate $a$ into the fiber of $E_N(\underline b')$ above $2a+\delta_0$, so $\widetilde R_N|_{C}$ lies in the pullback of $E_N(\underline b')$ along $\psi:C_0\to A$, $\psi(a)=2a+\delta_0$, and the restriction of $\varpi_\delta$ is the canonical map from that pullback. Since $\psi$ is the restriction to $C_0$ of $t_{\delta_0}\circ[2]_A$, and the class of $E_N(\underline b')$ lies in $\pic^{0}(A)$, hence is invariant under translation by Lemma~\ref{lem_phi_properties}(2),
\[
  \bigl(t_{\delta_0}\circ[2]_A\bigr)^{*}E_N(\underline b')
  =[2]_A^{*}\,t_{\delta_0}^{*}E_N(\underline b')
  \ \cong\ [2]_A^{*}E_N(\underline b') .
\]
This identifies $\widetilde R_N|_{C}$ with a subset $\Sigma$ of $[2]_A^{*}E_N(\underline b')$ over $C_0$, and the description in Step 1 gives
\[
  \Sigma=\pi_2^{-1}(C_0)\cap p^{-1}\bigl(\varsigma_C(C_0)\bigr),
\]
where $\pi_2:[2]_A^{*}E_N(\underline b')\to A$ is the projection, $p$ the pushout along $q_\ttt:\gm^{n}\to\gm^{n}/\ttt$, and $\varsigma_C$ the section over $C_0$ corresponding to $\varsigma|_{C}$.

Under this identification $\varpi_\delta$ is a $\gm^{n}$-equivariant map $[2]_A^{*}E_N(\underline b')\to E_N(\underline b')$ over $t_{\delta_0}\circ[2]_A$, and so is the canonical homomorphism $\varpi_2:[2]_A^{*}E_N(\underline b')\to E_N(\underline b')$ over $[2]_A$ followed by a translation above $\delta_0$. The two differ by a morphism $A\to\gm^{n}$, which is constant by Lemma~\ref{lem_noconst}. Hence $\varpi_\delta$ is $\varpi_2$ followed by a translation.

\smallskip
\noindent\textit{Step 4: applying the group-theoretic lemma.} The description of $\Sigma$ in Step 3 satisfies the hypotheses of Lemma~\ref{lem_elementary_fiber}\textup{(2)} for the semi-abelian variety $[2]_A^{*}E_N(\underline b')$, the algebraic subgroup $K$ and its coset $C_0$, the subtorus $\ttt\subseteq\gm^{n}$, and the section $\varsigma_C$. So every irreducible component of $\Sigma$ is a translate of an algebraic subgroup with toric part $\ttt$ and abelian quotient $K^{\circ}$. By Step 3 the map $\varpi_\delta$ is $\varpi_2$ followed by a translation, so every irreducible component of $\varpi_\delta(\widetilde R_N|_{C})$ is a translate of the image of that subgroup under $\varpi_2$; and since $\varpi_2$ is a homomorphism over $[2]_A$ and an isomorphism on toric parts, that image has toric part $\ttt$ and abelian quotient $[2]_A(K^{\circ})=K^{\circ}$.

As $K$ depends on neither $\underline b''$ nor $\underline b'$, taking the union over $\underline b''$ and letting $\underline b'$ vary proves \eqref{eq_fiber_to_prove} with $X:=K^{\circ}$, for which $[N]_A(X)=K^{\circ}$.
\end{proof}

\begin{corollary}\label{cor_fiber_dim}
With the notation of Theorem~\ref{thm_ribet_type_fibers}, the dimension of $K$ does not depend on $s$, all the nonempty fibers of $R\to\pr(R)$ have dimension $\dim K+\dim\ttt$, and
\[
  \dim R=\dim\pr(R)+\dim K+\dim\ttt .
\]
\end{corollary}
\begin{proof}
The morphism $\calA'\to\pr(\calA')$ is a surjective homomorphism of abelian schemes over $S$, hence induces a surjective homomorphism of abelian varieties on each fiber. Since the dimensions of $\calA'_s$ and $\pr(\calA')_s$ are independent of $s$, the dimension of their kernel $K$ is constant. By Theorem~\ref{thm_ribet_type_fibers}, every nonempty fiber of $R\to\pr(R)$ is a finite union of translates of the subgroup $H$, whose dimension is $\dim H=\dim K+\dim\ttt$. The dimension formula follows from the surjectivity of $R\to\pr(R)$.
\end{proof}

\begin{remark}
The proof uses the Ribet property only through the existence of the section $\tilde r$, supplied by Theorem~\ref{thm_ribet_property}, and the subtorus $\ttt$ only through the quotient homomorphism $q_\ttt$. The coset $C_0$ need not be a torsion translate of an abelian subvariety, as the example in \cite{BE20} shows. Consequently, a definition restricted to torsion translates of subgroup schemes or of Ribet sections would not include all the subvarieties obtained here.
\end{remark}

\section{Subvarieties of Ribet Type Are Special}\label{section5}

In this section we review the constructions in mixed Hodge theory needed to show that the subvarieties of Ribet type of Section~\ref{section4} are special subvarieties of the universal semi-abelian scheme $\pgn$. Our treatment largely follows \cite[\S\S4--5]{BE20}, where the case $n=1$, $\ttt=1$ and $\calA'=\Gamma(f-f^{\vee})$ is carried out.
\subsection{Review of Mixed Hodge Structures}\label{subsec_MHS}

For $k\in\ZZ$, a \textbf{$\ZZ$-Hodge structure of weight $k$} is a finitely generated $\ZZ$-module $M$ together with a decomposition
\[
  M_{\CC}=\bigoplus_{p+q=k}M^{p,q}_{\CC}
\]
such that $\overline{M^{p,q}_{\CC}}=M^{q,p}_{\CC}$ for all $p,q$ with $p+q=k$. Morphisms, $\Hom$'s and tensor products of Hodge structures are the standard ones, and we do not repeat them. The \textbf{Tate structure} $\ZZ(m)$ is the Hodge structure of weight $-2m$ whose underlying module is the $\ZZ$-submodule $(2\pi i)^{m}\ZZ$ of $\CC$ and whose only nonzero component in its Hodge decomposition lies in bi-degree $(-m,-m)$.

Let $C$ be the $\CC$-linear operator on $M_{\CC}$ acting on $M^{p,q}_{\CC}$ by $i^{-p}(-i)^{-q}$; it descends to an $\RR$-linear operator on $M_{\RR}$, still written $C$. A \textbf{polarization} of $M$ is a morphism of Hodge structures $\Psi:M\otimes M\to\ZZ(-k)$ for which the $\RR$-bilinear form
\[
  M_{\RR}\times M_{\RR}\longrightarrow\RR,\qquad (x,y)\longmapsto(2\pi i)^{k}\Psi(x\otimes Cy)
\]
is symmetric and positive definite.

\begin{example}[Principally polarized abelian varieties]\label{ex_ppav}
Let $g\geq1$, let $M$ be free of rank $2g$ with a Hodge decomposition $M_{\CC}=M^{-1,0}_{\CC}\oplus M^{0,-1}_{\CC}$, and let $\Psi:M\otimes M\to\ZZ(1)$ be a polarization inducing an isomorphism $M\to M^{\vee}(1)$. From this polarized Hodge structure one constructs a principally polarized abelian variety of dimension $g$ whose $\CC$-points are $M\backslash M_{\CC}/M^{0,-1}_{\CC}$, and by Riemann's theorem this construction is an equivalence of categories.

One may choose a basis of $M$ identifying $(M,\Psi)$ with $\ZZ^{2g}$ equipped with
\[
  \Psi:\ZZ^{2g}\otimes\ZZ^{2g}\to\ZZ(1),\qquad
  x\otimes y\mapsto 2\pi i\,x^{\top}\begin{pmatrix}\begin{smallmatrix}0&-1_g\\1_g&0\end{smallmatrix}\end{pmatrix}y .
\]
The set $D_{\Psi}$ of Hodge structures of type $\{(-1,0),(0,-1)\}$ on $\ZZ^{2g}$ polarized by $\Psi$ is in bijection with the Siegel upper half space $\HH_{g}=\{\tau\in\M_{g}(\CC):\tau^{\top}=\tau,\ \mathrm{Im}(\tau)>0\}$ via
\[
  \tau\longmapsto M^{0,-1}_{\CC}:=\begin{pmatrix}\tau\\1_g\end{pmatrix}\CC^{g} .
\]

Writing $\ZZ^{2g}=\M_{2g,1}(\ZZ)$ as columns, the abelian variety attached to $\tau$ is therefore
\[
  A_{\tau}=\CC^{2g}\Big/\Bigl(\begin{pmatrix}\tau\\1_g\end{pmatrix}\CC^{g}+\M_{2g,1}(\ZZ)\Bigr)
  =\M_{g,1}(\CC)\big/(1_g\ \ {-\tau})\,\M_{2g,1}(\ZZ),
\]
the second equality being the projection $\binom{x_1}{x_2}\mapsto x_1-\tau x_2$ along $M^{0,-1}_{\CC}$. Identifying $\M_{2g,1}(\ZZ)^{\vee}$ with $\M_{1,2g}(\ZZ)$ through the pairing
\[
  \M_{2g,1}(\ZZ)\times\M_{1,2g}(\ZZ)\longrightarrow\ZZ,\qquad(x,\xi)\longmapsto\xi x ,
\]
the dual Hodge structure is the one carried by $\M_{1,2g}(\ZZ(1))$, whose $(0,-1)$ part is the annihilator $\M_{1,g}(\CC)(1_g\ \ {-\tau})$ of $M^{0,-1}_{\CC}$; the corresponding abelian variety is the dual of $A_{\tau}$, and projecting along that annihilator gives in the same way
\[
  A_{\tau}^{\vee}=\M_{1,g}(\CC)\big/\M_{1,2g}(\ZZ(1))\begin{pmatrix}\tau\\1_g\end{pmatrix}.\qedhere
\]
\end{example}

\begin{remark}\label{rem_poltype}
For a general $\Psi:M\otimes M\to\ZZ(1)$, not necessarily a principal polarization, one may choose a basis identifying $(M,\Psi)$ with $\ZZ^{2g}$ equipped with
\[
  \Psi:\ZZ^{2g}\otimes\ZZ^{2g}\to\ZZ(1),\qquad
  x\otimes y\mapsto 2\pi i\,x^{\top}\begin{pmatrix}\begin{smallmatrix}0&-\delta\\\delta&0\end{smallmatrix}\end{pmatrix}y ,
\]
where $\delta=\mathrm{diag}\{e_1,\dots,e_g\}$ with $e_1\mid e_2\mid\cdots\mid e_g$ positive integers. Example~\ref{ex_ppav} then holds with $\HH_g$ replaced by $\delta^{-1}\HH_g:=\{\tau\in\M_g(\CC):(\delta\tau)^{\top}=\delta\tau,\ \mathrm{Im}(\delta\tau)>0\}$ and with principally polarized abelian varieties replaced by polarized abelian varieties of polarization type $\delta$. All the linear-algebraic descriptions below apply after these replacements, so to keep the notation light we write the group-theoretic computations in the principally polarized case.
\end{remark}

A \textbf{$\ZZ$-mixed Hodge structure} is a finitely generated $\ZZ$-module $M$ with the following two structures.
\begin{itemize}
  \item An increasing filtration $\{W_kM\}_{k\in\ZZ}$, the \textbf{weight filtration}, with $W_kM=M_{\mathrm{tor}}$ for $k\ll0$ and $W_kM=M$ for $k\gg0$, all $M/W_kM$ being torsion free.
  \item A decreasing filtration $\{F^{p}M_{\CC}\}_{p\in\ZZ}$, the \textbf{Hodge filtration}, with $F^{p}M_{\CC}=0$ for $p\gg0$ and $F^{p}M_{\CC}=M_{\CC}$ for $p\ll0$, such that the filtration induced by $F^{\bullet}$ on $(\Gr^{W}_{k}M)_{\CC}$ is a Hodge structure of weight $k$, that is,
  \[
    (\Gr^{W}_{k}M)_{\CC}=\bigoplus_{p+q=k}(\Gr^{W}_{k}M)^{p,q}_{\CC},
    \qquad
    (\Gr^{W}_{k}M)^{p,q}_{\CC}=F^{p}(\Gr^{W}_{k}M)_{\CC}\cap\overline{F^{q}(\Gr^{W}_{k}M)_{\CC}} .
  \]
\end{itemize}

The mixed Hodge structures relevant here are those of the $1$-motives of \cite[10.1.2]{hodge3}: a \textbf{$1$-motive} over an algebraically closed field $k$ is a homomorphism $c:X\to G(k)$ from a free abelian group $X$ to the group of $k$-points of a semi-abelian variety $G$ over $k$, written $[X\stackrel{c}{\to}G]$. We find this language helpful for the moduli interpretation of $\pgn$ as a mixed Shimura variety, and it is what guides the construction of the mixed Shimura data of Subsection~\ref{subsec_partial_rank_MSD}. The remainder of this subsection gives the mixed Hodge theoretic description of $\pgn$ itself: we exhibit the mixed Hodge structures whose associated $1$-motives are the points of $\pgn$, and then read off from that description the uniformization of $\pgn$.

\begin{notation}\label{not_M_and_D}
Fix $g\geq1$ and $n\geq1$, and let
\[
  M:=\ZZ(1)^{n}\oplus\ZZ^{2g}\oplus\ZZ
\]
with standard basis $2\pi ie_1,\dots,2\pi ie_n,e_{n+1},\dots,e_{n+2g},e_{n+2g+1}$, equipped with the weight filtration
\[
  W_{-3}M=\{0\},\quad
  W_{-2}M=\ZZ(1)^{n},\quad
  W_{-1}M=\ZZ(1)^{n}\oplus\ZZ^{2g},\quad
  W_{0}M=M .
\]
Let $D$ be the set of decreasing filtrations $F^{\bullet}$ on $M_{\CC}$ such that $(M,W_{\bullet},F^{\bullet})$ is a $\ZZ$-mixed Hodge structure of type $\{(-1,-1),(-1,0),(0,-1),(0,0)\}$ and such that
\[
  \Psi:\ZZ^{2g}\otimes\ZZ^{2g}\to\ZZ(1),\qquad
  x\otimes y\mapsto 2\pi i\,x^{\top}\begin{pmatrix}\begin{smallmatrix}0&-1_g\\1_g&0\end{smallmatrix}\end{pmatrix}y
\]
is a polarization of the pure Hodge structure $\Gr^{W}_{-1}M\cong\ZZ^{2g}$, the isomorphism being the one given by the chosen basis.
\end{notation}

Every filtration in $D$ is determined by $F^{0}$, being of the form $M_{\CC}=F^{-1}M_{\CC}\supset F^{0}M_{\CC}\supset F^{1}M_{\CC}=0$, and the possible $F^{0}$ are the following.

\begin{lemma}[{\cite[Proposition 4.5]{BE20}} for $n=1$]\label{lem_D_coordinates}
There is a bijection $\HH_g\times\M_{1,g}(\CC)^{n}\times\M_{g,1}(\CC)\times\CC^{n}\to D$,
\[
  (\tau,\underline u,v,\underline w)\longmapsto
  F^{0}M_{\CC}=\begin{pmatrix}
    \underline u&\underline w\\
    \tau&v\\
    1_g&0\\
    0&1
  \end{pmatrix}\CC^{g+1} .
\]
\end{lemma}
\begin{proof}
For $n=1$ this is \cite[Proposition 4.5]{BE20}, and the argument is the same for every $n$. Let $F^{\bullet}\in D$. Since $\Gr^{W}_{-2}M=\ZZ(1)^{n}$ is of type $(-1,-1)$ and $\Gr^{W}_{0}M=\ZZ$ of type $(0,0)$, the filtration induced by $F^{0}$ is $0$ on the first and everything on the second, while on $\Gr^{W}_{-1}M\cong\ZZ^{2g}$ it is the $M^{0,-1}_{\CC}$ of a unique $\tau\in\HH_g$ by Example~\ref{ex_ppav}. In particular $\dim_{\CC}F^{0}M_{\CC}=g+1$.

Now $F^{0}M_{\CC}\cap W_{-2}M_{\CC}$ is the degree-$0$ step of the induced filtration on $W_{-2}M$, hence zero, so $F^{0}M_{\CC}\cap W_{-1}M_{\CC}$ injects into $(\Gr^{W}_{-1}M)_{\CC}$ with image contained in $\binom{\tau}{1_g}\CC^{g}$; comparing dimensions, it maps onto it. It is therefore the graph of a unique $\CC$-linear map $\binom{\tau}{1_g}\CC^{g}\to\ZZ(1)^{n}_{\CC}$, that is, the span of the first $g$ columns of the displayed matrix for a unique $\underline u\in\M_{1,g}(\CC)^{n}$.

Finally $F^{0}M_{\CC}$ meets $W_{-1}M_{\CC}$ in a subspace of codimension one and maps onto $F^{0}(\Gr^{W}_{0}M)_{\CC}=\CC$, so it is spanned by those $g$ columns together with one further vector whose last coordinate is $1$. Since the block $1_g$ occurring in the first $g$ columns is invertible, that vector may be normalized uniquely so that its $\ZZ^{2g}$-component has vanishing second half; this is the last column of the matrix, and it determines unique $\underline w\in\CC^{n}$ and $v\in\M_{g,1}(\CC)$.

Conversely, for any $(\tau,\underline u,v,\underline w)$ the displayed subspace has dimension $g+1$ and induces on the three graded pieces the filtrations just described, which are Hodge structures of the required types, with $\Psi$ a polarization on $\Gr^{W}_{-1}M$ by Example~\ref{ex_ppav}. So it lies in $D$, and the two constructions are mutually inverse.
\end{proof}

\begin{remark}\label{rem_1motive_moduli}
Lemma~\ref{lem_D_coordinates} attaches to each $(\tau,\underline u,v,\underline w)$ a filtration $F^{\bullet}\in D$, hence a mixed Hodge structure $(M,W_{\bullet},F^{\bullet})$ as in Notation~\ref{not_M_and_D}. It gives rise to the $1$-motive
\[
  \bigl[\,W_{-1}M\backslash W_{0}M\longrightarrow W_{-1}M\backslash M_{\CC}/F^{0}M_{\CC}\,\bigr]
\]
over $\CC$. Here $W_{-1}M\backslash W_{0}M\cong\ZZ$, so the $1$-motive is of the form $[\ZZ\stackrel{c}{\to}G]$ with
\[
  G:=W_{-1}M\backslash M_{\CC}/F^{0}M_{\CC}
\]
a semi-abelian variety of toric rank $n$, whose abelian quotient is
\[
  A:=(\Gr^{W}_{-1}M)\backslash(\Gr^{W}_{-1}M)_{\CC}/F^{0}(\Gr^{W}_{-1}M)_{\CC},
\]
a principally polarized abelian variety by the polarization on $\Gr^{W}_{-1}M$. Thus the four coordinates have the following meanings, each of them being a point of a universal cover rather than of the variety it describes. The matrix $\tau$ is the moduli point of $A$, and the uniformizations it determines are
\[
  \M_{g,1}(\CC)\longrightarrow A,\qquad
  \M_{1,g}(\CC)\longrightarrow A^{\vee},\qquad
  \CC\longrightarrow\gm(\CC),\quad w\longmapsto e^{w},
\]
the kernels of the first two being the period lattices attached to $\tau$ and that of the third being $\ZZ(1)$. In these terms $\underline u\in\M_{1,g}(\CC)^{n}$ is a lift of the extension class of $G$, which by Theorem~\ref{thm_weil_barsotti} is an element of $(A^{\vee})^{n}(\CC)$; the coordinate $v\in\M_{g,1}(\CC)$ is a lift of the image of $c(1)$ in $A$; and the pair $(v,\underline w)$ is a lift of $c(1)$ itself, the coordinate $\underline w\in\CC^{n}$ being the one along the fibral $\gm^{n}$. Comparing with Proposition~\ref{prop_universal_semiabelian}, whose content is exactly that $\pt_{\calA,n}$ is the universal extension, we see that $D$ parametrizes the points of the universal semi-abelian scheme: the pairs $(\tau,\underline u)$ parametrize $(\Agv)^{[n]}$, the triples $(\tau,\underline u,v)$ parametrize $\calB_{g,n}$, and $D$ itself parametrizes $\pgn$. These last three parametrizations are stated formally, as identifications of quotients of $D$ by explicit subgroups of $P$, in Proposition~\ref{prop_uniformization} and Corollary~\ref{cor_uniformization} below.
\end{remark}

The kernels of the several uniformizations of Remark~\ref{rem_1motive_moduli} are organized by group theory into a single algebraic group acting on $D$: they are the integral points of its unipotent radical, one block for each of the three directions, and the description of $\pgn$ as a quotient of $D$ in Proposition~\ref{prop_uniformization} is what this organization buys. Let $P$ be the subgroup scheme of $\GL(M)\times\GL(\ZZ(1))$ fixing
\begin{itemize}
  \item the weight filtration $W_{\bullet}$,
  \item the map $\ZZ(1)^{n}\to W_{-2}M$, $(2\pi ia_i)_{1\leq i\leq n}\mapsto\sum_i2\pi ia_ie_i$,
  \item the map $\ZZ\to\Gr^{W}_{0}M$, $a\mapsto ae_{n+2g+1}$,
  \item the polarization form $\Psi:\Gr^{W}_{-1}M\otimes\Gr^{W}_{-1}M\to\ZZ(1)$.
\end{itemize}
For $g\in\GSp(\Psi)$ we write $\mu(g)$ for the similitude factor, so that $\Psi(gv,gw)=\mu(g)\Psi(v,w)$.

\begin{lemma}\label{lem_P_points}
For a ring $C$, with all matrices written in the $\ZZ$-basis of $M$ fixed in Notation~\textup{\ref{not_M_and_D}},
\[
  P(C)=\left\{\begin{pmatrix}
    \mu(g)1_n&\underline x&\underline z\\
    0&g&y\\
    0&0&1
  \end{pmatrix}\right\},
\]
where $(g,\mu(g))\in\GSp(\Psi)(C)$, $\underline x\in\M_{n,2g}(C)$, $y\in\M_{2g,1}(C)$ and $\underline z\in C^{n}$; its unipotent radical is
\[
  P^{u}(C)=\left\{\begin{pmatrix}
    1_n&\underline x&\underline z\\
    0&1&y\\
    0&0&1
  \end{pmatrix}:\underline x\in\M_{n,2g}(C),\ y\in\M_{2g,1}(C),\ \underline z\in C^{n}\right\},
\]
and $P$ has a normal subgroup $U\cong\Ga^{n}$ with
\[
  U(C)=\left\{\begin{pmatrix}
    1_n&0&\underline z\\
    0&1&0\\
    0&0&1
  \end{pmatrix}:\underline z\in C^{n}\right\}.
\]
In the $\CC$-basis $e_1,\dots,e_{n+2g+1}$ of $M_{\CC}$, obtained from that one by dropping the factors $2\pi i$, the blocks $\underline x$ and $\underline z$ are multiplied by $2\pi i$ and the block $y$ is unchanged, as in \textup{\cite[(4.5.2)]{BE20}}.
\end{lemma}
\begin{proof}
Let $(h,\lambda)$ be a $C$-point of $\GL(M)\times\GL(\ZZ(1))$, so that $\lambda\in\gm(C)$. Preserving the weight filtration means that $h$ is block upper triangular for the decomposition $M=\ZZ(1)^{n}\oplus\ZZ^{2g}\oplus\ZZ$ of Notation~\ref{not_M_and_D}, say
\[
  h=\begin{pmatrix}A&\underline x&\underline z\\0&g&c'\\0&0&c\end{pmatrix},
  \qquad A\in\GL_n(C),\quad g\in\GL_{2g}(C),\quad c\in\gm(C),
\]
with $\underline x\in\M_{n,2g}(C)$, $c'\in\M_{2g,1}(C)$ and $\underline z\in C^{n}$ unconstrained so far. The remaining three conditions pin down $A$, $c$ and $g$ in turn.

Fixing the map $\ZZ(1)^{n}\to W_{-2}M$ means that $A$, which is how $h$ acts on the target, agrees with the action of $\lambda$ on the source; as the latter is $\lambda$ on each of the $n$ factors, this forces $A=\lambda 1_n$. Fixing the map $\ZZ\to\Gr^{W}_{0}M$ means that $c$, which is how $h$ acts on the target, is the identity on the source, that is, $c=1$. Fixing the polarization form $\Psi:\Gr^{W}_{-1}M\otimes\Gr^{W}_{-1}M\to\ZZ(1)$, on which $h$ acts through $g$ and $\lambda$ acts on the target, means $\Psi(gv,gw)=\lambda\Psi(v,w)$, that is, $g\in\GSp(\Psi)(C)$ with similitude factor $\mu(g)=\lambda$. Writing $y$ for $c'$ gives the stated description, which is manifestly functorial in $C$.

The projection $h\mapsto g$ identifies the reductive quotient of $P$ with $\GSp(\Psi)$, so $P^{u}$ is the locus $g=1$, on which also $\lambda=\mu(1)=1$; this is the second display. Finally $U$ is the subgroup $\underline x=0$, $y=0$ of $P^{u}$, and a direct computation with the above matrices gives
\[
  h\,u(\underline z_0)\,h^{-1}=u\bigl(\mu(g)\underline z_0\bigr),
  \qquad
  u(\underline z_0):=\begin{pmatrix}1_n&0&\underline z_0\\0&1&0\\0&0&1\end{pmatrix},
\]
so that $U$ is normal in $P$ and $P$ acts on $U\cong\Ga^{n}$ through the similitude $\mu$.
\end{proof}

\begin{lemma}\label{lem_Pu_law}
Writing an element of $P^{u}$ as $u=(\underline x,y,\underline z)$, multiplication in $P^{u}$ is
\[
  (\underline x,y,\underline z)\cdot(\underline x',y',\underline z')
  =(\underline x+\underline x',\ y+y',\ \underline z+\underline z'+\underline xy'),
  \qquad
  (\underline x,y,\underline z)^{-1}=(-\underline x,-y,-\underline z+\underline xy).
\]
Consequently:
\begin{enumerate}[label=\textup{(\arabic*)}]
  \item $[(\underline x,y,\underline z),(\underline x',y',\underline z')]=(0,0,\underline xy'-\underline x'y)\in U$;
  \item $(P/U)^{u}=P^{u}/U\cong\M_{n,2g}\times\M_{2g,1}$ is abelian;
  \item $1\to U\to P^{u}\to(P/U)^{u}\to1$ is exact and $U$ is central in $P^{u}$.
\end{enumerate}
\end{lemma}
\begin{proof}
Multiply the matrices of Lemma~\ref{lem_P_points}.
\end{proof}

We record the geometric dictionary, which we use repeatedly: by Corollary~\ref{cor_uniformization}, the block $y$ is the $\Ag$-direction, the block $\underline x$ the $n$ directions of $(\Agv)^{[n]}$, and the block $\underline z$ the fibral $\gm^{n}$-direction. Finally, $G_{T}$ acts on $U\cong\Ga^{n}$ through the similitude $\mu$, that is, by scalars, so \emph{every} $\QQ$-subspace of $U_{\QQ}$ is $G_{T}$-stable; this is why in Subsection~\ref{subsec_converse} the subspace $Q\cap U_{\QQ}$ may be in arbitrary position.

We refer to \cite[\S4]{BE20} for the explicit formulas for the actions of these groups on $D$. Quotienting $D$ by the action of $P(\ZZ)$ amounts to forgetting the choice of basis of $M$ made in Notation~\ref{not_M_and_D}, and the quotient is to be viewed as the moduli space of the $1$-motives just described; to obtain a fine moduli space one first quotients by the smaller group $P^{u}(\ZZ)$.

\begin{proposition}[{\cite[Proposition 4.6]{BE20}} for $n=1$]\label{prop_uniformization}
The quotient $P^{u}(\ZZ)\backslash D$ is the universal semi-abelian scheme $\pgn$ over $\HH_g$, that is, the pullback of $\pgn$ along the uniformization $\HH_g\to\ag$.
\end{proposition}
\begin{proof}
For $n=1$ this is \cite[Proposition 4.6]{BE20}. The general case follows by taking fiber powers. In the coordinates of Lemma~\ref{lem_D_coordinates}, $D$ is the $n$-fold fiber power over $\HH_g\times\M_{g,1}(\CC)$ of the corresponding space for $n=1$, the $i$-th factor contributing the pair $(\underline u_i,\underline w_i)$; by Lemma~\ref{lem_P_points} the group $P^{u}(\ZZ)$ is correspondingly the $n$-fold fiber power, over the block $y$, of the group for $n=1$, the blocks $\underline x$ and $\underline z$ splitting into their $n$ rows and coordinates. The actions match factor by factor. Since $\pgn$ is by Definition~\ref{def_fiber_power_torsor} the $n$-fold fiber power of $\pt_{\Ag}$ over $\Ag$, and since $\Ag$ over $\HH_g$ is the quotient of $\HH_g\times\M_{g,1}(\CC)$ by the block $y$, the case $n=1$ gives the general one.
\end{proof}

\begin{corollary}\label{cor_uniformization}
The quotient $P^{u}(\ZZ)\M_{n,2g}(\RR)U(\CC)\backslash D$ is the universal abelian scheme $\Ag$ over $\HH_g$, the quotient $P^{u}(\ZZ)\M_{2g,1}(\RR)U(\CC)\backslash D$ is $(\Agv)^{[n]}$ over $\HH_g$, and the quotient $P^{u}(\ZZ)U(\CC)\backslash D$ is $\calB_{g,n}=\Ag\times_{\HH_g}(\Agv)^{[n]}$ over $\HH_g$.
\end{corollary}
\begin{proof}
In the coordinates $(\tau,\underline u,v,\underline w)$ of Lemma~\ref{lem_D_coordinates}, the three unipotent blocks act by translations in the last three coordinates: $\underline z$ translates $\underline w$, $\underline x$ translates $\underline u$, and $y$ translates $v$; see \cite[\S4]{BE20} for the formulas. Dividing $D$ by $U(\CC)$ therefore forgets $\underline w$ and leaves $\HH_g\times\M_{1,g}(\CC)^{n}\times\M_{g,1}(\CC)$, whose quotient by $P^{u}(\ZZ)$ is $\calB_{g,n}$ over $\HH_g$ by Proposition~\ref{prop_uniformization}, the fibral $\gm^{n}$ having been divided out. Dividing in addition by $\M_{n,2g}(\RR)$ forgets $\underline u$ and leaves $\Ag$, and dividing instead by $\M_{2g,1}(\RR)$ forgets $v$ and leaves $(\Agv)^{[n]}$.
\end{proof}

Quotienting further by a congruence subgroup of $\GSp(\Psi)(\QQ)$ replaces $\HH_g$ by a Siegel moduli variety; the level so entering is the one suppressed throughout, as in Remark~\ref{rem_minimal_family}, and by Remark~\ref{rem_coverings} it plays no role in what follows.

In the language of \cite[Ch.~2]{PinkThesis}, all of this says that $(P,D)$ is a mixed Shimura datum whose associated mixed Shimura variety is $\pgn$. This is the point of view we take for the rest of the article.

\subsection{Mixed Shimura Data of Ribet Type}\label{subsec_partial_rank_MSD}

This subsection runs parallel to Section~\ref{section4}. There a Ribet datum over a base $S$ was defined and a subvariety of Ribet type was constructed from it; here the base is a special subvariety $T\subseteq\ag$, and we attach to each Ribet datum over $T$ a mixed Shimura datum of Ribet type. The passage proceeds in two steps, and it is the first that this subsection is about. At the level of data, we translate a Ribet datum over $T$ into linear algebra --- a lattice and a subtorus, subject to a matrix identity --- and Proposition~\ref{prop_datum_bijection} states that the translation is a bijection; the mixed Shimura datum is then constructed from the linear side in Notation~\ref{not_MSD_list}. At the level of the geometric objects they define, Subsection~\ref{subsec_comparison} identifies the mixed Shimura variety of that datum with the subvariety of Ribet type of Section~\ref{section4}. Along the way we generalize the Hodge tensor \cite[(5.1.1)]{BE20} to a family of Hodge tensors indexed by the characters trivial on $\ttt$, and compute the algebraic group of the datum explicitly.

Throughout the subsection we fix a special subvariety $T\subseteq\ag$ and write $(G_{T},\HH_{g,T})$ for its pure Shimura datum, and we keep the identifications of Example~\ref{ex_ppav}: $\M_{2g,1}(\ZZ)^{\vee}$ with $\M_{1,2g}(\ZZ)$ through the pairing $(x,\xi)\mapsto\xi x$, and $\M_{1,2g}(\ZZ(1))$ with the local system of $\Agv$.

\subsubsection*{The lattice of an abelian subscheme}

Let $\Lambda$ be a \textbf{saturated} $G_{T}$-stable $\ZZ$-submodule of $\M_{2g,1}(\ZZ)\oplus\M_{1,2g}(\ZZ(1))^{n}$ of rank $2d$. Being $G_{T}$-stable, $\Lambda$ is a sub-Hodge structure at every point of $\HH_{g,T}$, the Mumford--Tate group at such a point being contained in $G_{T}$; and since $\M_{2g,1}(\ZZ)\oplus\M_{1,2g}(\ZZ(1))^{n}$ is polarizable and $\Lambda$ is saturated, $\Lambda$ is polarizable as well. By \cite[Proposition 3.4]{gaoziyang} the assignment $\Lambda\mapsto\calA'$ is a bijection between the saturated $G_{T}$-stable submodules of $\M_{2g,1}(\ZZ)\oplus\M_{1,2g}(\ZZ(1))^{n}$ and the abelian subschemes $\calA'\subseteq\calB_{g,n}|_{T}$ over $T$, the module $\Lambda$ being the local system underlying $\calA'$, of rank twice the relative dimension $d$. This is the first half of the dictionary announced above; the second half is the subtorus, and it is a tautology.

Choosing a basis of $\Lambda$, so that $\Lambda\cong\M_{2d,1}(\ZZ)$, the inclusion $\Lambda\hookrightarrow\M_{2g,1}(\ZZ)\oplus\M_{1,2g}(\ZZ(1))^{n}$ is given by matrices
\[
  \alpha_{\ZZ}\in\M_{2g,2d}(\ZZ),\qquad \beta_{i,\ZZ}\in\M_{2d,2g}(\ZZ)\quad(1\leq i\leq n),
\]
namely by $x\mapsto\alpha_{\ZZ}x$ into $\M_{2g,1}(\ZZ)$ and by $x\mapsto2\pi i\,x^{\top}\beta_{i,\ZZ}$ into the $i$-th copy of $\M_{1,2g}(\ZZ(1))$. These matrices are the ones of the projections $\alpha:\calA'\to\Ag|_{T}$ and $\beta_i:\calA'\to\Agv|_{T}$ of Subsection~\ref{subsubsec_ribet_data}, read on the underlying local systems. For $\chi\in\ZZ^{n}$ we keep the notation $\beta_\chi=\sum_i\chi_i\beta_i$ and write $\beta_{\chi,\ZZ}=\sum_i\chi_i\beta_{i,\ZZ}$.

\subsubsection*{The subtorus and its characters}

Let $\ttt\subseteq\gm^{n}$ be a subtorus and, as in Definition~\ref{def_ribet_datum}, read
\[
  \Xi:=\chl(\gm^{n}/\ttt)\subseteq\chl(\gm^{n})=\ZZ^{n}
\]
as the saturated subgroup of characters of $\gm^{n}$ trivial on $\ttt$. The assignment $\ttt\mapsto\Xi$ is a bijection between the subtori of $\gm^{n}$ and the saturated subgroups of $\ZZ^{n}$; this is the second half of the dictionary, and the subtorus is carried along unchanged. Under $U\cong\Ga^{n}$ of Lemma~\ref{lem_P_points} the subtorus $\ttt$ corresponds in turn to the $\QQ$-subspace
\[
  U':=\{\underline z\in U:\ \chi\underline z=0\ \text{ for every }\chi\in\Xi\},
  \qquad \chi\underline z:=\textstyle\sum_i\chi_iz_i ,
\]
and $\ttt\mapsto U'$ is a bijection between the subtori of $\gm^{n}$ and the $\QQ$-subspaces of $U_\QQ$.

\subsubsection*{Ribet data in linear algebra}

The Ribet property involves the duals $\alpha^{\vee}$ and $\beta_\chi^{\vee}$, so we say what records them. By Example~\ref{ex_ppav} the local system of $(\calA')^{\vee}$ is $\M_{1,2d}(\ZZ(1))$ and that of $\Agv|_{T}=(\Ag|_{T})^{\vee}$ is $\M_{1,2g}(\ZZ(1))$, while the source of $\beta_i^{\vee}$ is $(\Agv|_{T})^{\vee}=\Ag|_{T}$, whose local system is $\M_{2g,1}(\ZZ)$. So $\alpha^{\vee}$ and $\beta_i^{\vee}$ are given on local systems by
\[
\begin{aligned}
  \M_{1,2g}(\ZZ(1))&\to\M_{1,2d}(\ZZ(1)), & \xi&\mapsto\xi\,(\alpha^{\vee})_{\ZZ},\\
  \M_{2g,1}(\ZZ)&\to\M_{1,2d}(\ZZ(1)), & y&\mapsto2\pi i\,y^{\top}(\beta_i^{\vee})_{\ZZ},
\end{aligned}
\]
and these displays define the matrices $(\alpha^{\vee})_{\ZZ},(\beta_i^{\vee})_{\ZZ}\in\M_{2g,2d}(\ZZ)$.

\begin{lemma}\label{lem_dual_matrices}
With the notation above, $(\alpha^{\vee})_{\ZZ}=\alpha_{\ZZ}$ and $(\beta_i^{\vee})_{\ZZ}=-\beta_{i,\ZZ}^{\top}$ for $1\leq i\leq n$; explicitly,
\[
\begin{aligned}
  \alpha^{\vee}&:\M_{1,2g}(\ZZ(1))\to\M_{1,2d}(\ZZ(1)), & \xi&\longmapsto \xi\alpha_{\ZZ},\\
  \beta_i^{\vee}&:\M_{2g,1}(\ZZ)\to\M_{1,2d}(\ZZ(1)), & y&\longmapsto-2\pi i\,y^{\top}\beta_{i,\ZZ}^{\top}.
\end{aligned}
\]
By linearity, $(\beta_\chi^{\vee})_{\ZZ}=-\beta_{\chi,\ZZ}^{\top}$ for every $\chi\in\ZZ^{n}$.
\end{lemma}
\begin{proof}
Fix $\tau\in\HH_{g,T}$ and take $n=1$, writing $\beta=\beta_1$; the $n$ factors of $\calB_{g,n}$ are independent of one another. Choose the basis of $\Lambda$ so that $F^{0}\Lambda_{\CC}=\binom{\tau_{\Lambda}}{1_d}\M_{d,1}(\CC)$ for some $\tau_{\Lambda}$ as in Remark~\ref{rem_poltype}, and write $\alpha_{\CC}\in\M_{g,d}(\CC)$ and $\beta_{\CC}\in\M_{d,g}(\CC)$ for the maps $s\mapsto\alpha_{\CC}s$ and $s\mapsto s^{\top}\beta_{\CC}$ induced on universal covers. Then there is a commutative square
\[\begin{tikzcd}[column sep=8ex, row sep=5ex]
  \M_{2d,1}(\ZZ)\ar[r,"{(\alpha_{\ZZ},\beta_{\ZZ})}"]\ar[d] & \M_{2g,1}(\ZZ)\oplus\M_{1,2g}(\ZZ(1))\ar[d] \\
  \M_{d,1}(\CC)\ar[r,"{(\alpha_{\CC},\beta_{\CC})}"] & \M_{g,1}(\CC)\oplus\M_{1,g}(\CC)
\end{tikzcd}\]
whose vertical maps are the uniformizations of Example~\ref{ex_ppav}, namely $s\mapsto(1_d\ \ {-\tau_{\Lambda}})s$ on the left and $(y,\xi)\mapsto\bigl((1_g\ \ {-\tau})y,\ \xi\binom{\tau}{1_g}\bigr)$ on the right. Its commutativity reads
\begin{equation}\label{eq_periods}
  \alpha_{\CC}(1_d\ \ {-\tau_{\Lambda}})=(1_g\ \ {-\tau})\,\alpha_{\ZZ},
  \qquad
  \begin{pmatrix}1_d\\-\tau_{\Lambda}^{\top}\end{pmatrix}\beta_{\CC}=2\pi i\,\beta_{\ZZ}\begin{pmatrix}\tau\\1_g\end{pmatrix} ;
\end{equation}
this is the analogue of \cite[(4.7.1)]{BE20}. We split $\alpha_{\ZZ}=\binom{\alpha_1}{\alpha_2}$ and $\beta_{\ZZ}=\binom{\beta_1}{\beta_2}$ into halves, and further $\alpha_j=(\alpha_{j1}\ \alpha_{j2})$ with $\alpha_{jk}\in\M_{g,d}(\ZZ)$ and $\beta_j=(\beta_{j1}\ \beta_{j2})$ with $\beta_{jk}\in\M_{d,g}(\ZZ)$.

Next we record the two families of extensions carried by the Poincar\'e torsor, in the normalization of \cite[(4.6.1)]{BE20}: an extension of a complex torus $X=V/L$ by $\gm$ is the cokernel of
\[
  \ZZ(1)\oplus L\longrightarrow\CC\oplus V,
  \qquad (2\pi ik,\ell)\longmapsto\bigl(2\pi ik-\theta(\ell),\,\ell\bigr),
\]
for a homomorphism $\theta\in\Hom_{\ZZ}(L,\CC)$, its \textbf{period homomorphism}; the class of the extension is the image of $\theta$ in $\Hom_{\ZZ}(L,\CC)$ modulo the $\CC$-linear functionals on $V$ and modulo $\Hom_{\ZZ}(L,\ZZ(1))$, and that quotient is $X^{\vee}$ with the uniformization of Example~\ref{ex_ppav}. The sign in the display is the one chosen in \cite[(4.6.1)]{BE20}. In the coordinates $(\tau,u,v,w)$ of Lemma~\ref{lem_D_coordinates}, \cite[(4.6.2)]{BE20} says that the fibre of $\pt_{\Ag}$ over a point $u$ of $A_{\tau}^{\vee}$ is the extension of $A_{\tau}$ by $\gm$ with
\[
  \theta_{u}\binom{m_1}{m_2}=u\,m_2,\qquad \binom{m_1}{m_2}\in\M_{2g,1}(\ZZ),
\]
and \cite[(4.6.3)]{BE20} says that its fibre over a point $v$ of $A_{\tau}$ is the extension of $A_{\tau}^{\vee}$ by $\gm$ with
\[
  \theta_{v}\bigl(2\pi i(n_1\ n_2)\bigr)=-2\pi i\,n_1v,\qquad 2\pi i(n_1\ n_2)\in\M_{1,2g}(\ZZ(1)) .
\]
The first of these, read with $\tau_{\Lambda}$ and $d$ in place of $\tau$ and $g$, gives the extension of $\calA'_{\tau}$ by $\gm$ attached to a point $u'$ of $(\calA'_{\tau})^{\vee}$, namely
\[
  \theta'_{u'}\binom{p}{q}=u'q,\qquad\binom{p}{q}\in\M_{2d,1}(\ZZ) .
\]

We can now compute, as in \cite[Proposition 4.8]{BE20}. Write $x=\binom{p}{q}$ for an element of $\M_{2d,1}(\ZZ)$, so that the period lattice of $\calA'_{\tau}$ consists of the $p-\tau_{\Lambda}q$, and recall from the previous paragraph that two period homomorphisms on it define the same extension exactly when they differ by the restriction of a $\CC$-linear functional on $\M_{d,1}(\CC)$, up to $\ZZ(1)$.

\emph{The matrix of $\alpha^{\vee}$.} Let $\xi=2\pi i(\xi_1\ \xi_2)$ be in $\M_{1,2g}(\ZZ(1))$ and let $u=\xi\binom{\tau}{1_g}=2\pi i(\xi_1\tau+\xi_2)$ be the corresponding point of $A_{\tau}^{\vee}$. The lower half of $\alpha_{\ZZ}x$ is $\alpha_2x=\alpha_{21}p+\alpha_{22}q$, so the pullback $\theta_{u}\circ\alpha_{\ZZ}$ is $x\mapsto u(\alpha_{21}p+\alpha_{22}q)$. Subtracting the $\CC$-linear functional $s\mapsto u\alpha_{21}s$, which on the period lattice takes the value $u\alpha_{21}(p-\tau_{\Lambda}q)$, leaves $u(\alpha_{22}+\alpha_{21}\tau_{\Lambda})q$; comparing with $\theta'_{u'}$ gives
\[
  \alpha^{\vee}(u)=u'=u\bigl(\alpha_{22}+\alpha_{21}\tau_{\Lambda}\bigr).
\]
On the other hand the point of $(\calA'_{\tau})^{\vee}$ attached to $\xi\alpha_{\ZZ}\in\M_{1,2d}(\ZZ(1))$ is
\[
  \xi\alpha_{\ZZ}\begin{pmatrix}\tau_{\Lambda}\\1_d\end{pmatrix}
  =2\pi i\bigl[(\xi_1\alpha_{11}+\xi_2\alpha_{21})\tau_{\Lambda}+\xi_1\alpha_{12}+\xi_2\alpha_{22}\bigr],
\]
and substituting $\alpha_{12}=-\alpha_{11}\tau_{\Lambda}+\tau(\alpha_{21}\tau_{\Lambda}+\alpha_{22})$, which is what the first identity in \eqref{eq_periods} says on the last $d$ columns, this collapses to $2\pi i(\xi_1\tau+\xi_2)(\alpha_{21}\tau_{\Lambda}+\alpha_{22})=u'$. So $\alpha^{\vee}$ carries $\xi$ to $\xi\alpha_{\ZZ}$, that is, $(\alpha^{\vee})_{\ZZ}=\alpha_{\ZZ}$.

\emph{The matrix of $\beta^{\vee}$.} Let $y=\binom{y_1}{y_2}$ be in $\M_{2g,1}(\ZZ)$ and let $v=y_1-\tau y_2$ be the corresponding point of $A_{\tau}$. The first half of $x^{\top}\beta_{\ZZ}$ is $p^{\top}\beta_{11}+q^{\top}\beta_{21}$, so the pullback $\theta_{v}\circ\beta_{\ZZ}$ is
\[
  x\longmapsto-2\pi i\bigl(p^{\top}\beta_{11}+q^{\top}\beta_{21}\bigr)v
  =-2\pi i\,v^{\top}\bigl(\beta_{11}^{\top}p+\beta_{21}^{\top}q\bigr).
\]
Subtracting the $\CC$-linear functional $s\mapsto-2\pi i\,v^{\top}\beta_{11}^{\top}s$ leaves $-2\pi i\,v^{\top}(\beta_{21}^{\top}+\beta_{11}^{\top}\tau_{\Lambda})q$; comparing with $\theta'_{u'}$ as before gives
\[
  \beta^{\vee}(v)=u'=-2\pi i\,v^{\top}\bigl(\beta_{21}^{\top}+\beta_{11}^{\top}\tau_{\Lambda}\bigr).
\]
On the other hand the point of $(\calA'_{\tau})^{\vee}$ attached to $-2\pi i\,y^{\top}\beta_{\ZZ}^{\top}\in\M_{1,2d}(\ZZ(1))$ is
\[
  -2\pi i\,y^{\top}\beta_{\ZZ}^{\top}\begin{pmatrix}\tau_{\Lambda}\\1_d\end{pmatrix}
  =-2\pi i\,y^{\top}\bigl(\beta_{1}^{\top}\tau_{\Lambda}+\beta_{2}^{\top}\bigr),
\]
whose first $g$ rows are $\beta_{11}^{\top}\tau_{\Lambda}+\beta_{21}^{\top}$ and whose last $g$ rows are $\beta_{12}^{\top}\tau_{\Lambda}+\beta_{22}^{\top}$. Transposing the second identity in \eqref{eq_periods} and using $\tau^{\top}=\tau$ gives
\[
  \beta_{12}^{\top}\tau_{\Lambda}+\beta_{22}^{\top}=-\tau\bigl(\beta_{21}^{\top}+\beta_{11}^{\top}\tau_{\Lambda}\bigr),
\]
so that $\beta_{1}^{\top}\tau_{\Lambda}+\beta_{2}^{\top}=\binom{1_g}{-\tau}(\beta_{21}^{\top}+\beta_{11}^{\top}\tau_{\Lambda})$. Since $v^{\top}=y^{\top}\binom{1_g}{-\tau}$, the two displays agree, and $(\beta^{\vee})_{\ZZ}=-\beta_{\ZZ}^{\top}$.
\end{proof}

\begin{proposition}\label{prop_ribet_symmetric}
For $\chi\in\ZZ^{n}$, the pair $(\alpha,\beta_\chi)$ satisfies the Ribet property of Definition~\ref{def_ribet_property} if and only if the matrix $\beta_{\chi,\ZZ}\alpha_{\ZZ}$ is symmetric.
\end{proposition}
\begin{proof}
The Ribet property of Definition~\ref{def_ribet_property} reads $\alpha^{\vee}\beta_\chi+\beta_\chi^{\vee}\alpha=0$, an identity between homomorphisms $\calA'\to(\calA')^{\vee}$ of abelian schemes over $T$. Such a homomorphism vanishes if and only if the induced map of the underlying local systems does, and these are $\Lambda\cong\M_{2d,1}(\ZZ)$ for $\calA'$ and $\M_{1,2d}(\ZZ(1))$ for $(\calA')^{\vee}$. On lattices, $\alpha$ and $\beta_\chi$ are given by the definition of $\alpha_{\ZZ}$ and $\beta_{\chi,\ZZ}$, and $\alpha^{\vee}$ and $\beta_\chi^{\vee}$ by Lemma~\ref{lem_dual_matrices}. Composing,
\[
  \alpha^{\vee}\beta_\chi:\ x\longmapsto2\pi i\,x^{\top}\beta_{\chi,\ZZ}\alpha_{\ZZ},
  \qquad
  \beta_\chi^{\vee}\alpha:\ x\longmapsto-2\pi i\,x^{\top}(\beta_{\chi,\ZZ}\alpha_{\ZZ})^{\top},
\]
so the left-hand side of the Ribet property is $x\mapsto2\pi i\,x^{\top}\bigl(\beta_{\chi,\ZZ}\alpha_{\ZZ}-(\beta_{\chi,\ZZ}\alpha_{\ZZ})^{\top}\bigr)$, and it vanishes if and only if $\beta_{\chi,\ZZ}\alpha_{\ZZ}$ is symmetric.
\end{proof}

Changing the basis of $\Lambda$ by $P\in\GL_{2d}(\ZZ)$ replaces $\alpha_{\ZZ}$ by $\alpha_{\ZZ}P$ and each $\beta_{i,\ZZ}$ by $P^{\top}\beta_{i,\ZZ}$, hence $\beta_{\chi,\ZZ}\alpha_{\ZZ}$ by $P^{\top}\beta_{\chi,\ZZ}\alpha_{\ZZ}P$; the symmetry of $\beta_{\chi,\ZZ}\alpha_{\ZZ}$ is therefore independent of the basis.

The two halves of the dictionary now combine into the bijection at the level of data.

\begin{definition}\label{def_linear_ribet_datum}
A \textbf{linear Ribet datum} over $T$ is a pair $(\Lambda,\ttt)$ consisting of a saturated $G_{T}$-stable submodule $\Lambda\subseteq\M_{2g,1}(\ZZ)\oplus\M_{1,2g}(\ZZ(1))^{n}$ and a subtorus $\ttt\subseteq\gm^{n}$, such that the matrix $\beta_{\chi,\ZZ}\alpha_{\ZZ}$ is symmetric for every $\chi\in\Xi$.
\end{definition}

\begin{proposition}\label{prop_datum_bijection}
The assignment $(\calA',\ttt)\mapsto(\Lambda,\ttt)$ is a bijection between the Ribet data for $\pgn|_{T}$ in the sense of Definition~\ref{def_ribet_datum} and the linear Ribet data over $T$.
\end{proposition}
\begin{proof}
On first entries the assignment is the bijection between abelian subschemes of $\calB_{g,n}|_{T}$ over $T$ and saturated $G_{T}$-stable submodules recorded above, and on second entries it is the identity, so only the matching of the conditions has to be checked. By Definition~\ref{def_ribet_datum} the Ribet conditions of $(\calA',\ttt)$ say that $(\alpha,\beta_\chi)$ has the Ribet property for every $\chi\in\Xi$, and by Proposition~\ref{prop_ribet_symmetric} this is the symmetry of $\beta_{\chi,\ZZ}\alpha_{\ZZ}$ for every $\chi\in\Xi$.
\end{proof}

The rest of the subsection constructs, from a linear Ribet datum, the mixed Shimura datum of Ribet type attached to it.

\subsubsection*{The Hodge tensors of a Ribet datum}

For $\chi\in\ZZ^{n}$ let $M_\chi:=M/\ker\bigl(\chi:W_{-2}M=\ZZ(1)^{n}\to\ZZ(1)\bigr)$ be the pushout of $M$ along $\chi$, a lattice of the shape $\ZZ(1)\oplus\ZZ^{2g}\oplus\ZZ$ with the induced weight and Hodge filtrations; on the geometric side, pushing out along $\chi$ replaces the semi-abelian variety $G$ of Remark~\ref{rem_1motive_moduli} by its quotient $\chi_{*}G$ of toric rank one, and $M_\chi$ is the lattice of \cite[\S4]{BE20}. Define
\[
  \Theta_{\Lambda,\chi}:\ \ZZ(1)\oplus\M_{2d,1}(\ZZ)\oplus\ZZ\longrightarrow M_\chi\oplus M_\chi^{\vee}(1)
\]
by
\[
  (k,w,m)\longmapsto\bigl((k,\alpha_{\ZZ}w,m),\,(-m,\,2\pi i\,w^{\top}\beta_{\chi,\ZZ},\,-k)\bigr),
\]
the second component being written in the coordinates of
\[
  M_\chi^{\vee}(1)=\Hom(M_\chi,\ZZ(1))\cong\ZZ\oplus\M_{1,2g}(\ZZ(1))\oplus\ZZ(1)
\]
dual to $M_\chi=\ZZ(1)\oplus\ZZ^{2g}\oplus\ZZ$. Give the source the weight filtration
\[
  W_{-3}=\{0\},\quad W_{-2}=\ZZ(1),\quad W_{-1}=\ZZ(1)\oplus\M_{2d,1}(\ZZ),\quad W_{0}=\ZZ(1)\oplus\M_{2d,1}(\ZZ)\oplus\ZZ .
\]
For $n=1$, $\chi=1$ and $\calA'$ the graph of $\alpha$ in $\calB_{g,1}|_{T}$, so that $\beta_{\ZZ}=1_{2g}$, the image of $\Theta_{\Lambda,1}$ is the graph of the map \cite[(5.1.1)]{BE20}.

That these maps are Hodge tensors is not an extra hypothesis: it is equivalent to $\Lambda$ being a sub-Hodge structure at every point of $\HH_{g,T}$, which follows from the $G_{T}$-stability required in Definition~\ref{def_linear_ribet_datum}.

\begin{lemma}\label{lem_tensor_is_hodge}
Let $\chi\in\ZZ^{n}$ and let $F^{\bullet}\in D$ lie over a point of $\HH_{g,T}$. Then the filtration induced by $\Theta_{\Lambda,\chi}$ makes the source a $\ZZ$-mixed Hodge structure, so that $\Theta_{\Lambda,\chi}$ is a morphism of mixed Hodge structures, that is, a Hodge tensor at $F^{\bullet}$.
\end{lemma}
\begin{proof}
The map is strict for the weight filtrations by inspection: the summand $\ZZ(1)$ goes to weight $-2$ on both sides, the summand $\ZZ$ to weight $0$ on both sides, and $\M_{2d,1}(\ZZ)$ to weight $-1$ on both sides. So the induced filtration makes the source a mixed Hodge structure as soon as it induces a Hodge structure on each graded piece; on $\Gr^{W}_{-2}=\ZZ(1)$ and $\Gr^{W}_{0}=\ZZ$ this is automatic, and on $\Gr^{W}_{-1}=\Lambda$ it holds because $\Lambda$ is a sub-Hodge structure.
\end{proof}

\begin{notation}\label{not_Q_lambda}
Let $Q'_{\Lambda}$ be the subgroup of $P/U$ consisting of the classes modulo $U$ of the elements of $P$ whose reductive part lies in $G_{T}$ and whose unipotent blocks $(y,\underline x)$ satisfy $(y,\underline xg^{-1})\in\Lambda$; explicitly, for a ring $C$, $Q'_{\Lambda}(C)$ consists of the classes modulo $U(C)$ of the matrices
\[
  \begin{pmatrix}
    \mu(g)1_n&(v^{\top}\beta_{i,\ZZ}\,g)_{1\leq i\leq n}&\ast\\
    0&g&\alpha_{\ZZ}v\\
    0&0&1
  \end{pmatrix},
  \qquad v\in\M_{2d,1}(C),\ g\in G_{T}(C),
\]
where $G_{T}\subseteq\GSp(\Psi)$ is the group of the pure Shimura datum of $T$ fixed at the beginning of the subsection, which preserves $\Lambda$ by hypothesis. Thus $Q'_{\Lambda}$ is parametrized by the pairs $(v,g)$, and by \cite[Proposition 3.4]{gaoziyang} it is the group of the mixed Shimura datum of the abelian subscheme $\calA'$; it does not involve $\ttt$.

Let $Q_{\Lambda,\ttt}$ be the subgroup of the preimage of $Q'_{\Lambda}$ in $P$ consisting of the elements that stabilize the image of $\Theta_{\Lambda,\chi}$ for every $\chi\in\Xi$. Its reductive quotient is $G_{T}$.
\end{notation}

Besides being a condition on the Ribet datum, the Ribet property is what allows the quadratic corner $\tfrac12v^{\top}\beta_{\chi,\ZZ}\alpha_{\ZZ}v$ below to exist at all: the symmetry of $\beta_{\chi,\ZZ}\alpha_{\ZZ}$ is the symmetry of a $2$-cocycle, and it is what lets that cocycle be trivialized by a quadratic form.

The functor of points of $Q_{\Lambda,\ttt}$ can now be computed explicitly: the conditions imposed by the tensors are conditions on the fibral coordinate $\underline z$ alone, one for each character in $\Xi$. For $\ttt=1$, so $\Xi=\ZZ^{n}$, and $n=1$ the following is \cite[(5.1.3)]{BE20}, where it is stated without proof; the factor $\tfrac12$ in the upper right entry is the reason for working over $\ZZ[\tfrac12]$, as in the remark on the unipotent radical in loc.\ cit.

\begin{proposition}\label{prop_Q_points}
Let $(\Lambda,\ttt)$ be a linear Ribet datum over $T$. Then, for every $\ZZ[\tfrac12]$-algebra $C$, the group $Q_{\Lambda,\ttt}(C)$ consists of the matrices
\[
  \begin{pmatrix}
    \mu(g)1_n&(v^{\top}\beta_{i,\ZZ}\,g)_{1\leq i\leq n}&\underline z\\
    0&g&\alpha_{\ZZ}v\\
    0&0&1
  \end{pmatrix}
\]
with $v\in\M_{2d,1}(C)$, $g\in G_{T}(C)$ and $\underline z\in C^{n}$ subject to
\[
  \chi\underline z=\tfrac12 v^{\top}\beta_{\chi,\ZZ}\alpha_{\ZZ}v
  \qquad\text{for every }\chi\in\Xi .
\]
In particular $Q_{\Lambda,\ttt}\cap U=U'$.
\end{proposition}
\begin{proof}
Let $\tilde g$ be an element of the preimage of $Q'_{\Lambda}(C)$ in $P(C)$, so that
\[
  \tilde g=\begin{pmatrix}\mu(g)1_n&\underline x&\underline z\\0&g&y\\0&0&1\end{pmatrix},
  \qquad y=\alpha_{\ZZ}v,\quad \underline x=(v^{\top}\beta_{i,\ZZ}\,g)_{i},
\]
for some $v\in\M_{2d,1}(C)$ and $g\in G_{T}(C)$, and
\[
  \tilde g^{-1}=\begin{pmatrix}
    \mu(g)^{-1}1_n&-\mu(g)^{-1}\underline xg^{-1}&\mu(g)^{-1}(\underline xg^{-1}y-\underline z)\\
    0&g^{-1}&-g^{-1}y\\
    0&0&1
  \end{pmatrix}.
\]
Fix $\chi\in\Xi$ and write $x:=\chi\underline x=v^{\top}\beta_{\chi,\ZZ}g$ and $z:=\chi\underline z$ for the $\chi$-components, which are the entries of the image of $\tilde g$ in the automorphisms of $M_\chi$. That image acts on $(M_\chi)_R$ by left multiplication and on $(M_\chi^{\vee}(1))_R$ by right multiplication by $\mu(g)\tilde g^{-1}$, so
\begin{itemize}
  \item $\tilde g\cdot(k,w_1,m)=(\mu(g)k+xw_1+mz,\ gw_1+my,\ m)$,
  \item $\tilde g\cdot(-m,w_2,-k)=\bigl(-m,\ mxg^{-1}+\mu(g)w_2g^{-1},\ z'\bigr)$ with
  \[
    z'=-m(xg^{-1}y-z)-\mu(g)w_2g^{-1}y-\mu(g)k .
  \]
\end{itemize}
Hence $\tilde g$ stabilizes the image of $\Theta_{\Lambda,\chi}$ if and only if for every $(k,w,m)$, with $w_1:=\alpha_{\ZZ}w$ and $w_2:=w^{\top}\beta_{\chi,\ZZ}$,
\begin{enumerate}
  \item $(gw_1+my,\ mxg^{-1}+\mu(g)w_2g^{-1})$ again lies in the image of $\Lambda_C$ under $(\alpha_{\ZZ},\beta_{\chi,\ZZ})$, and
  \item $xw_1+2mz=(mxg^{-1}+\mu(g)w_2g^{-1})y$.
\end{enumerate}
Condition (1) holds automatically: $\Lambda_C$ is stable under $G_{T}(C)$, so $g\cdot(w_1,w_2)$ comes from $\Lambda_C$, and $(y,xg^{-1})=(\alpha_{\ZZ}v,v^{\top}\beta_{\chi,\ZZ})$ does so by construction. For condition (2), write $(gw_1,\mu(g)w_2g^{-1})=(\alpha_{\ZZ}u,u^{\top}\beta_{\chi,\ZZ})$, which is possible by the same stability. Then
\[
  xw_1=v^{\top}\beta_{\chi,\ZZ}\alpha_{\ZZ}u,
  \qquad
  (\mu(g)w_2g^{-1})y=u^{\top}\beta_{\chi,\ZZ}\alpha_{\ZZ}v ,
\]
and these agree because $\beta_{\chi,\ZZ}\alpha_{\ZZ}$ is symmetric. So (2) holds for every $(k,w,m)$ if and only if it holds for $m=1$, $w=0$, that is, if and only if
\[
  \chi\underline z=z=\tfrac12 xg^{-1}y=\tfrac12 v^{\top}\beta_{\chi,\ZZ}\alpha_{\ZZ}v .
\]
Letting $\chi$ run through $\Xi$ gives the description in the statement, and its last assertion is the case $v=0$, $g=1$.
\end{proof}

\begin{remark}\label{rem_Q_group}
That the set of Proposition~\ref{prop_Q_points} is a group is not evident from its shape, and it is here that the Ribet conditions are used a second time. Given $v,w\in\M_{2d,1}(C)$ and $g,h\in G_{T}(C)$, the product of the two corresponding matrices has upper blocks
\[
  \Bigl(\mu(gh)1_n,\ \bigl(\mu(g)w^{\top}\beta_{i,\ZZ}h+v^{\top}\beta_{i,\ZZ}gh\bigr)_i,\ \underline z''\Bigr),
\]
where
\[
  \chi\underline z''=\tfrac12\mu(g)w^{\top}\beta_{\chi,\ZZ}\alpha_{\ZZ}w
  +v^{\top}\beta_{\chi,\ZZ}g\alpha_{\ZZ}w
  +\tfrac12v^{\top}\beta_{\chi,\ZZ}\alpha_{\ZZ}v ,
\]
and lower blocks $(0,gh,g\alpha_{\ZZ}w+\alpha_{\ZZ}v)$. Since $\Lambda_C$ is $G_{T}(C)$-stable there is $y\in\M_{2d,1}(C)$ with $g\alpha_{\ZZ}w=\alpha_{\ZZ}y$ and $\mu(g)w^{\top}\beta_{i,\ZZ}=y^{\top}\beta_{i,\ZZ}g$ for all $i$; substituting and using the symmetry of $\beta_{\chi,\ZZ}\alpha_{\ZZ}$ turns the displayed data into the data attached to $x:=y+v$ and $gh$, which again satisfies the equations of Proposition~\ref{prop_Q_points}. Closure under inverses is checked in the same way.
\end{remark}

\begin{remark}\label{rem_Q_extremes}
The two extremes of Subsection~\ref{subsubsec_ribet_type} are visible in Proposition~\ref{prop_Q_points}. For $\ttt=\gm^{n}$ one has $\Xi=0$, no equation is imposed on $\underline z$, and $Q_{\Lambda,\gm^{n}}$ is the full preimage in $P$ of the group of $\calA'$; the corresponding subvariety is the preimage $\pi^{-1}(\calA')$. For $\ttt=1$ one has $\Xi=\ZZ^{n}$, every coordinate of $\underline z$ is pinned to $\tfrac12v^{\top}\beta_{i,\ZZ}\alpha_{\ZZ}v$, and $Q_{\Lambda,1}\cap U=0$; the corresponding subvariety is the Ribet multi-section itself. In general the toric directions of $\ttt$ are free and the directions indexed by $\Xi$ are pinned, one scalar equation for each character.
\end{remark}

So far no point of $\HH_{g,T}$ has been used. To pass from these data to the geometric families they define we now fix one: we choose $\tau\in\HH_{g,T}$, we let
\[
  \tilde\tau:=(\tau,0,0,0)\in\HH_g\times\M_{1,g}(\CC)^{n}\times\M_{g,1}(\CC)\times\CC^{n}\ =\ D ,
\]
the identification being the bijection of Lemma~\ref{lem_D_coordinates}, and we endow $M$ with the mixed Hodge structure that $\tilde\tau$ defines. We write $A:=A_\tau$ for the fiber of $\Ag$ at $\tau$, the abelian variety of Example~\ref{ex_ppav}, so that
\[
  A=\M_{g,1}(\CC)/(1_g\ \ {-\tau})\M_{2g,1}(\ZZ),
  \qquad
  A^{\vee}=\M_{1,g}(\CC)/\M_{1,2g}(\ZZ(1))\begin{pmatrix}\begin{smallmatrix}\tau\\1_g\end{smallmatrix}\end{pmatrix} .
\]
Nothing below depends on which $\tau\in\HH_{g,T}$ is taken, the whole construction being $G_{T}(\RR)^{+}$-equivariant.

\begin{remark}\label{rem_tensor_meaning}
The summand $M_\chi^{\vee}(1)$ carries the dual Hodge filtration, and for a single $\chi$ the situation is that of \cite[\S5]{BE20} applied to the pushout of $M$ along $\chi$. The geometric meaning of the tensor is the following. Its image picks out a point of the Poincar\'e torsor of $A$ whose projection to $A\times A^{\vee}$ lies in the image of $(\alpha,\beta_\chi)$; under the uniformization of Proposition~\ref{prop_uniformization}, the point $\tilde\tau$ then corresponds to a point of the Ribet multi-section of Proposition-Definition~\ref{def_ribet_section} in the $\chi$-direction. The choice $\tilde\tau=(\tau,0,0,0)$ places that point over the origin of $A'$, at the unit element of the fibral $\gm^{n}$. So $\Theta_{\Lambda,\chi}$ \emph{is} the Ribet multi-section in the $\chi$-direction, written in linear algebra.
\end{remark}

We now list the mixed Shimura data involved, as in \cite[\S5]{BE20}.

\begin{notation}\label{not_MSD_list}
Let $(\Lambda,\ttt)$ be a linear Ribet datum over $T$, and let $(\calA',\ttt)$ be the Ribet datum it corresponds to under Proposition~\ref{prop_datum_bijection}.
\begin{itemize}
  \item \emph{The base.} The pure Shimura datum is the $(G_{T},\HH_{g,T})$ of $T$, fixed at the beginning of the subsection; recall that $\HH_{g,T}=G_{T}(\RR)^{+}\cdot\tau$, the point $\tau$ having just been chosen in $\HH_{g,T}$.
  \item \emph{The ambient family.} Let $P_{T}$ be the preimage of $G_{T}$ under $P\to P/P^{u}=\GSp(\Psi)$ and $E_{T}:=P_{T}(\RR)^{+}U(\CC)\cdot\tilde\tau$; the mixed Shimura datum is $(P_{T},E_{T})$ and the corresponding subvariety is $\pgn|_{T}$.
  \item \emph{The family of abelian varieties.} Let $P'_{T}:=(P_{T}\cdot U)/U$ be the image of $P_{T}$ in $P/U$ and $E'_{T}:=P'_{T}(\RR)^{+}\cdot\tau'$; the mixed Shimura datum is $(P'_{T},E'_{T})$ and the corresponding subvariety is $\calB_{g,n}|_{T}$.
  \item \emph{The abelian subscheme.} Let $Q'_{\Lambda}:=(Q_{\Lambda,\ttt}\cdot U)/U$, the image of $Q_{\Lambda,\ttt}$ in $P/U$, and $D'_{\Lambda}:=Q'_{\Lambda}(\RR)^{+}\cdot\tau'$, where $\tau'$ is the image of $\tilde\tau$ in $\HH_g\times\M_{1,g}(\CC)^{n}\times\M_{g,1}(\CC)$; the mixed Shimura datum is $(Q'_{\Lambda},D'_{\Lambda})$ and the corresponding subvariety is the abelian subscheme $\calA'\subseteq\calB_{g,n}|_{T}$. It depends only on $\Lambda$, not on $\ttt$.
  \item \emph{The subvariety of Ribet type.} Put $D_{\Lambda,\ttt}:=Q_{\Lambda,\ttt}(\RR)^{+}U'(\CC)\cdot\tilde\tau\subseteq D$; the mixed Shimura datum is $(Q_{\Lambda,\ttt},D_{\Lambda,\ttt})$.\qedhere
\end{itemize}
\end{notation}

The last item is the construction promised at the start of the subsection: a linear Ribet datum $(\Lambda,\ttt)$ over $T$, hence by Proposition~\ref{prop_datum_bijection} a Ribet datum $(\calA',\ttt)$ for $\pgn|_{T}$, produces the mixed Shimura datum $(Q_{\Lambda,\ttt},D_{\Lambda,\ttt})$, which we call the \textbf{mixed Shimura datum of Ribet type} attached to it. That the geometric object it defines is the subvariety of Ribet type $R(\calA',\ttt)$ of Section~\ref{section4} is the content of Subsection~\ref{subsec_comparison}.

Torsion translations are visible on the level of these data, and this is what will allow us, in Subsection~\ref{subsec_finish}, to dispose of them.

\begin{lemma}[dictionary]\label{lem_hecke_dictionary}
Write $\pgn|_{T}=\Gamma_P\backslash E_{T}$ with $\Gamma_P=P^{u}(\ZZ)\rtimes\Gamma_{T}$, where $\Gamma_{T}\subseteq G_{T}(\QQ)$ is the congruence subgroup giving the level. For $\nu\in P^{u}(\QQ)$ put $\Gamma_{\nu}:=\Gamma_P\cap\nu^{-1}\Gamma_P\nu$ and consider the Hecke correspondence
\[
  \pgn|_{T}=\Gamma_P\backslash E_{T}
  \ \xleftarrow{\ \varrho\ }\ \Gamma_{\nu}\backslash E_{T}
  \ \xrightarrow{\ \cdot\,\nu\ }\ \nu\Gamma_{\nu}\nu^{-1}\backslash E_{T}
  \ \xrightarrow{\ \varrho'\ }\ \Gamma_P\backslash E_{T}=\pgn|_{T} .
\]
Then:
\begin{enumerate}[label=\textup{(\arabic*)}]
  \item $\Gamma_{\nu}$ and $\nu\Gamma_{\nu}\nu^{-1}$ are congruence subgroups of finite index in $\Gamma_P$, so all three quotients are mixed Shimura varieties and $\varrho$, $\varrho'$ are finite Shimura morphisms;
  \item for an irreducible closed $Z\subseteq\pgn|_{T}$, the irreducible components of $\varrho'\bigl(\nu\cdot\varrho^{-1}(Z)\bigr)$ are the torsion translates of $Z$ in the sense of Definition~\ref{def_torsion_translation} attached to the class of $\nu$, under the identifications
  \[
    (P/U)^{u}(\QQ)/(P/U)^{u}(\ZZ)\ \longleftrightarrow\ \calB_{g,n}|_{T,\mathrm{tors}},
    \qquad
    U(\QQ)/U(\ZZ)\ \longleftrightarrow\ \mu_{\infty}^{n}\subseteq\gm^{n},
  \]
  compatible with the exact sequences of Lemma~\ref{lem_Pu_law}\textup{(3)} and Proposition~\ref{prop_heisenberg_exact}; conversely every torsion translate of $Z$ arises this way;
  \item a special subvariety with datum $(Q,D_Q)$ is carried to the special subvariety with datum $(\nu Q\nu^{-1},\nu D_Q)$.
\end{enumerate}
\end{lemma}
\begin{proof}
(1) Conjugation by a rational point of $P^{u}$ carries congruence subgroups to congruence subgroups: for the unipotent part it multiplies denominators by a bounded factor, and for the Levi part $\nu^{-1}\gamma\nu=\gamma\cdot(\gamma^{-1}\nu^{-1}\gamma\nu)$ with the second factor in $P^{u}(\QQ)$ with denominators bounded along a congruence subgroup. So $\Gamma_{\nu}$ contains a principal congruence subgroup, and the quotients are mixed Shimura varieties at a deeper level; see \cite[Ch.~3]{PinkThesis} for Hecke correspondences of mixed Shimura data. Images and preimages of special subvarieties along $\varrho,\varrho'$ are special, and dimensions are preserved, these morphisms being finite and surjective.

(2) Unwinding $P^{u}(\ZZ)\backslash D=\pgn|_{\HH_g}$ of Proposition~\ref{prop_uniformization}, the operation ``take a component of the preimage under a finite covering, translate, project down'' is Definition~\ref{def_torsion_translation} for a suitable level $N$; the two identifications are the statements that $\M_{n,2g}(\QQ)/\M_{n,2g}(\ZZ)$ and $\M_{2g,1}(\QQ)/\M_{2g,1}(\ZZ)$ are the torsion of $(\Agv)^{[n]}$ and of $\Ag$, and that $\QQ^{n}/\ZZ^{n}\cong\mu_{\infty}^{n}$ through $\underline z\mapsto e^{2\pi i\underline z}$. The lift of a translation of $\calB_{g,n}$ to $\pt_{\Ag,n}$ is unique up to $\gm^{n}$ by Proposition~\ref{prop_heisenberg_exact}, which matches the ambiguity of $\nu$ modulo $U(\ZZ)$.

(3) Immediate from the description of the correspondence on the uniformizing spaces.
\end{proof}

\begin{proposition}\label{prop_torsion_conjugation}
A torsion translate of a special subvariety of $\pgn$ is again special, and conversely; for $(Q,D_Q)$ and $\nu\in P^{u}(\QQ)$ the translate has datum $(\nu Q\nu^{-1},\nu D_Q)$. In particular, for $\nu\in U(\QQ)$ the translation is by a torsion section of the fibral $\gm^{n}$, namely by the image of $\nu$ in $U(\QQ)/U(\ZZ)\cong\mu_{\infty}^{n}$.
\end{proposition}
\begin{proof}
This is Lemma~\ref{lem_hecke_dictionary}: specialty is preserved because the correspondence consists of finite Shimura morphisms between mixed Shimura varieties, by part (1), and the effect on data is part (3). The converse follows since torsion translation is invertible, Remark~\ref{rem_translation_design}(2), the inverse being the translation attached to $\nu^{-1}$.
\end{proof}

\begin{remark}\label{rem_translation_level}
The coverings appearing in Definition~\ref{def_torsion_translation}, where a point of a $\gm^{n}$-torsor and an $N^{2}$-th root of unity have to be extracted, are arbitrary finite surjective ones and are of no consequence: by Lemma~\ref{lem_hecke_dictionary} the operation they compute is the Hecke correspondence attached to $\nu$, whose coverings $\varrho,\varrho'$ are congruence, hence Shimura morphisms. It is along the latter that specialty and codimension are transported; see Remark~\ref{rem_coverings}.

The correspondence of data is exact, not merely one of existence. By Proposition~\ref{prop_heisenberg_exact} the lifts of a torsion section $c$ of $\calB_{g,n}$ to $\calH_N$ form a $\gm^{n}$-torsor, of which the torsion ones form a coset of $\mu_{\infty}^{n}$; under Lemma~\ref{lem_hecke_dictionary} the component of $\nu$ in $(P/U)^{u}(\QQ)/(P/U)^{u}(\ZZ)$ is $c$, and its component in $U(\QQ)/U(\ZZ)\cong\mu_{\infty}^{n}$ is the choice of torsion lift. So every torsion translation available geometrically is a Hecke move, and conversely. That the lift be torsion is essential: a lift differing from a torsion one by an element of $\gm^{n}$ which is not a root of unity translates by a non-torsion point of the fibral $\gm^{n}$, which is not a Hecke move and does not preserve specialty --- this is the phenomenon of \cite{BE20} that also appears in Theorem~\ref{thm_ribet_type_fibers}, where the components of a fiber are translates but need not be torsion translates.
\end{remark}

\subsection{Comparison with the Algebraic Construction}\label{subsec_comparison}

We can now carry out the second step: the geometric objects defined by the two sides of Proposition~\ref{prop_datum_bijection} agree.

\begin{proposition}\label{prop_comparison}
Let $(\Lambda,\ttt)$ be a linear Ribet datum over $T$ and $(\calA',\ttt)$ the corresponding Ribet datum, as in Notation~\ref{not_MSD_list}. Then $(Q'_{\Lambda})^{u}(\ZZ)\backslash D'_{\Lambda}$ is the abelian subscheme $\calA'$ of $\calB_{g,n}|_{T}$, and
\[
  Q_{\Lambda,\ttt}^{u}(\ZZ)\backslash D_{\Lambda,\ttt}=R(\calA',\ttt),
\]
the subvariety of Ribet type of Definition~\ref{def_ribet_type} attached to the Ribet datum $(\calA',\ttt)$ over $T$.
\end{proposition}
\begin{proof}
The first assertion is Corollary~\ref{cor_uniformization} together with the description of $Q'_{\Lambda}$ in Notation~\ref{not_Q_lambda}: a $G_{T}$-stable saturated sublattice of $\M_{2g,1}(\ZZ)\oplus\M_{1,2g}(\ZZ(1))^{n}$ is the same thing as an abelian subscheme of $\calB_{g,n}|_{T}$, by \cite[Proposition 3.4]{gaoziyang}.

For the second, let $(\tilde Q_{\Lambda,\ttt})^{u}(2\ZZ)$ be the group of matrices
\[
  \begin{pmatrix}
    1_n&(v^{\top}\beta_{i,\ZZ})_i&\underline z\\
    0&1&\alpha_{\ZZ}v\\
    0&0&1
  \end{pmatrix},
  \qquad v\in\M_{2d,1}(2\ZZ),\quad \chi\underline z=\tfrac12v^{\top}\beta_{\chi,\ZZ}\alpha_{\ZZ}v\ \ (\chi\in\Xi),
\]
and let $(\tilde Q'_{\Lambda})^{u}(2\ZZ)$ be its image in $(P/U)(\ZZ)$; the factor $2$ makes the entries integral. Then $(\tilde Q_{\Lambda,\ttt})^{u}(2\ZZ)$ is a subgroup of both $(P_{T})^{u}(\ZZ)$ and $(Q'_{\Lambda})^{u}(\ZZ)$, and we obtain the diagram
\[
\begin{tikzcd}[column sep=10ex]
(\tilde Q_{\Lambda,\ttt})^{u}(2\ZZ)\backslash D_{\Lambda,\ttt} \arrow[d,"\cong"'] \arrow[r]
& Q_{\Lambda,\ttt}^{u}(\ZZ)\backslash D_{\Lambda,\ttt} \arrow[d] \arrow[r]
& (P_{T})^{u}(\ZZ)\backslash E_{T} \arrow[d]\\
(\tilde Q'_{\Lambda})^{u}(2\ZZ)\backslash D'_{\Lambda} \arrow[r,"{[2]}"]
& (Q'_{\Lambda})^{u}(\ZZ)\backslash D'_{\Lambda} \arrow[r]
& (P'_{T})^{u}(\ZZ)\backslash E'_{T}\ \rlap{,}
\end{tikzcd}
\]
whose left vertical arrow is an isomorphism because $Q_{\Lambda,\ttt}\cap U=U'$ has been divided out on the left and because the equations of Proposition~\ref{prop_Q_points} determine the remaining coordinates of $\underline z$. Comparing it with the two cartesian squares of Subsection~\ref{subsubsec_ribet_type}, whose composite defines $R(\calA',\ttt)$ from the Ribet multi-section $\tilde r$ of $(\calA',\ttt)$, we see that both objects are, over $[2]_{\calA'}$, the locus in the $\gm^{n}$-torsor $\pt_{\Ag,n}|_{\calA'}$ where the coordinates indexed by $\Xi$ take the values prescribed by a trivialization of $\pt_{\Ag,n}|_{\calA',\gm^{n}/\ttt}$, and where the coordinates in $\ttt$ are free. Two such loci differ by the section of $\gm^{n}/\ttt$ over $\calA'$ by which the two trivializations differ, that is, by a morphism $\calA'\to\gm^{n}/\ttt$. Such a morphism is constant on the fibers of $\calA'\to T$ by Lemma~\ref{lem_noconst}, and by the choice of $\tilde\tau$ in Remark~\ref{rem_tensor_meaning} it takes the value $1$ above the origin of $\calA'$; hence it is trivial and the two loci coincide. The trivialization in question is therefore the Ribet multi-section, by the uniqueness in Proposition-Definition~\ref{def_ribet_section}(2).
\end{proof}

\begin{theorem}\label{thm_ribet_type_special}
Let $R\subseteq\pgn$ be a subvariety of Ribet type in the sense of Definition~\ref{def_ribet_type}, for $\calA=\Ag$ over a base $T$ which is a special subvariety of $\ag$. Then $R$ is a special subvariety of $\pgn$.
\end{theorem}
\begin{proof}
By Definition~\ref{def_ribet_type} we may write $R$ as a torsion translate of $R(\calA',\ttt)$ for a Ribet datum $(\calA',\ttt)$ over $T$. Apply the setup of Subsection~\ref{subsec_partial_rank_MSD} to this $T$: let $(G_{T},\HH_{g,T})$ be its pure Shimura datum and fix $\tau\in\HH_{g,T}$. By Proposition~\ref{prop_datum_bijection} the Ribet datum $(\calA',\ttt)$ corresponds to a linear Ribet datum $(\Lambda,\ttt)$ over $T$. Hence Proposition~\ref{prop_Q_points} applies, $(Q_{\Lambda,\ttt},D_{\Lambda,\ttt})$ is a mixed Shimura subdatum of $(P_{T},E_{T})$, and $R(\calA',\ttt)$ is the special subvariety it defines by Proposition~\ref{prop_comparison}. Finally $R$, being a torsion translate of it, is special by Proposition~\ref{prop_torsion_conjugation}.
\end{proof}

\section{Special Subvarieties of \texorpdfstring{$\pgn$}{P\^{}x\_\{g,n\}} Are of Ribet Type}\label{section6}

We now prove the converse of Theorem~\ref{thm_ribet_type_special}: every special subvariety of $\pgn$ arises from the construction of Section~\ref{section4}. Recall the two projections
\[
  \Pi:\pgn\to\ag,\qquad \pi:\pgn\to\calB_{g,n} .
\]
Let $R\subseteq\pgn$ be a special subvariety, with mixed Shimura datum $(Q,D_Q)$. Images of special subvarieties under Shimura morphisms are special, so $T:=\Pi(R)$ is a special subvariety of $\ag$ and $W:=\pi(R)$ is a special subvariety of $\calB_{g,n}$. We restrict everything to $T$ and write $(G_{T},\HH_{g,T})$ for its pure Shimura datum, so that $(Q,D_Q)$ is a mixed Shimura subdatum of the datum $(P_{T},E_{T})$ of $\pgn|_{T}$ in Notation~\ref{not_MSD_list}.

The proof is divided into three steps.
\begin{itemize}
  \item First we replace $R$ by a torsion translate for which $\pi(R)$ is honestly an abelian subscheme $\calA'$ of $\calB_{g,n}|_{T}$, rather than a torsion translate of one. This is Subsection~\ref{subsec_moving} and is independent of $n$.
  \item Then, in Subsection~\ref{subsec_converse}, we read off the toric part: the subspace $U'=Q\cap U_{\QQ}$ corresponds to a subtorus $\ttt\subseteq\gm^{n}$, and dividing by $\ttt$ leaves a subvariety finite over $\calA'$. Finiteness forces the Ribet conditions, so $(\calA',\ttt)$ is a Ribet datum. This does not yet identify $R$ with $R(\calA',\ttt)$: being finite over $\calA'$ is much weaker than being the particular finite cover cut out by the Ribet multi-section. To get the identification we compute the group $Q$ and find that it is $Q_{\Lambda,\ttt}$ up to conjugation by a rational point of $U$.
  \item That conjugation is a torsion translation in the toric direction, by Proposition~\ref{prop_torsion_conjugation}, and removing it in Subsection~\ref{subsec_finish} identifies $R$ with a torsion translate of $R(\calA',\ttt)$.
\end{itemize}

\subsection{Step 1: Clear Torsion Translations in the Abelian Part}\label{subsec_moving}

\begin{proposition}\label{prop_step1}
There is a torsion translate $R_1$ of $R$, in the sense of Definition~\ref{def_torsion_translation}, such that $\pi(R_1)$ is an abelian subscheme of $\calB_{g,n}|_{T}$ over $T$. It is again a special subvariety of $\pgn|_{T}$, and $R$ is in turn a torsion translate of it.
\end{proposition}
\begin{proof}
The special subvariety $W=\pi(R)$ of $\calB_{g,n}$ has image $T$ in $\ag$, and by \cite[Proposition 3.4]{gaoziyang}, applied to each of the $n+1$ factors of $\calB_{g,n}=\Ag\times_{\ag}(\Agv)^{[n]}$, it is of the form
\[
  W=\calA'+t
\]
for an abelian subscheme $\calA'\subseteq\calB_{g,n}|_{T}$ over $T$ and a torsion multi-section $t$ of $\calB_{g,n}|_{T}$; passing to a level as in Remark~\ref{rem_coverings} we may take $t$ to be a torsion section, killed by $N$ say. Choose a torsion section $c$ of $\calB_{g,n}|_{T}$, again after passing to a level, with $Nc=-t$; it is killed by $N^{2}$, so by Lemma~\ref{lem_torsion_lift} it lifts to a torsion element $\tau\in\calH_N$, and we let $R_1:=R+\tau$ be the corresponding torsion translate of Definition~\ref{def_torsion_translation}. By Proposition~\ref{prop_fibral_transport} the composite of $\tau$ with $\varpi_N$ lies over the map $\calB_{g,n}\to\calB_{g,n}$, $b\mapsto N(b+c)$, so that
\[
  \pi(R_1)=[N]\bigl([N]^{-1}W+c\bigr)=W+Nc=\calA'+t-t=\calA' ,
\]
an abelian subscheme over $T$. That $R_1$ is again special is Proposition~\ref{prop_torsion_conjugation}: the translation performed is a torsion translation, so by Lemma~\ref{lem_hecke_dictionary} it is the Hecke correspondence attached to a $\nu\in P^{u}(\QQ)$, and $R_1$ has mixed Shimura datum $(\nu Q\nu^{-1},\nu D_Q)$. In particular the finite coverings used in the construction of $\tau$ carry no specialty statement; that is Remark~\ref{rem_translation_level}. Finally $R$ is a torsion translate of $R_1$ by Remark~\ref{rem_translation_design}(2).
\end{proof}

\begin{remark}\label{rem_order_of_translations}
In the treatment of \cite[\S4]{BE20} the two directions were cleared one after the other, first the $(\Agv)^{[n]}$-direction and then the $\Ag$-direction, and the two orders differ by a translation in the fibral $\gm^{n}$, namely by the Weil pairing of the two torsion sections involved. Definition~\ref{def_torsion_translation} clears both at once, through a single element of the geometric Heisenberg group $\calH_N$, so that no such discrepancy arises here; on the group side this is Lemma~\ref{lem_Pu_law}(2), the translations of $\calB_{g,n}$ forming an abelian group whatever $n$ is. The fibral ambiguity is exactly what Subsection~\ref{subsec_finish} removes.
\end{remark}

By Proposition~\ref{prop_step1} we may replace $R$ by $R_1$: the latter is again special, so it has a mixed Shimura datum to work with, and being of Ribet type is by Definition~\ref{def_ribet_type} stable under torsion translation, so the conclusion for $R_1$ gives the conclusion for $R$. We therefore assume for the rest of this section that
\[
  \pi(R)=\calA'
\]
is an abelian subscheme of $\calB_{g,n}|_{T}$ over $T$. By \cite[Proposition 3.4]{gaoziyang} again, $\calA'$ is given by a $G_{T}$-stable saturated sublattice $\Lambda$ of $\M_{2g,1}(\ZZ)\oplus\M_{1,2g}(\ZZ(1))^{n}$, with matrices $\alpha_{\ZZ}$ and $\beta_{i,\ZZ}$ as in Subsection~\ref{subsec_partial_rank_MSD}, and its mixed Shimura datum is the $(Q'_{\Lambda},D'_{\Lambda})$ of Notation~\ref{not_MSD_list}.

\subsection{Step 2: Identify the Ribet Property}\label{subsec_converse}

Put
\[
  U':=Q\cap U_{\QQ},
\]
a $\QQ$-subspace of $U_{\QQ}\cong\QQ^{n}$, and let $\ttt\subseteq\gm^{n}$ be the subtorus corresponding to it as in Subsection~\ref{subsec_partial_rank_MSD}, with $\Xi=\chl(\gm^{n}/\ttt)$ the group of characters trivial on $\ttt$. Recall from Subsection~\ref{subsubsec_ribet_sections} the pushout
\[
  \pt_{\Ag,n}\big|_{\calA',\gm^{n}/\ttt}=(q_\ttt)_{*}\bigl(\pt_{\Ag,n}|_{\calA'}\bigr),
\]
a $\gm^{n}/\ttt$-torsor over $\calA'$, and write $R_\ttt$ for the image of $R$ in it.

\begin{lemma}\label{lem_toric_part_of_R}
The subtorus $\ttt$ is the largest subtorus of $\gm^{n}$ whose fibral action on $\pgn$ preserves $R$, and $R_\ttt\to\calA'$ is finite and surjective.
\end{lemma}
\begin{proof}
By Corollary~\ref{cor_uniformization} the action of $U(\CC)$ on $D$ is, modulo $U(\ZZ)$, the fibral action of $\gm^{n}$. A subtorus preserves $R$ if and only if the corresponding subspace of $U_{\QQ}$ is contained in $Q$, whence the first assertion. For the second, $R_\ttt$ is the image of $R$ under the Shimura morphism given by the pushout along $q_\ttt$, hence a special subvariety, and the group of its mixed Shimura datum is $Q/U'$, whose intersection with $(U/U')_{\QQ}$ is trivial by the definition of $U'$. By \cite[Facts 2.6(a)]{pink04} the fibers of $R_\ttt$ over its image $\calA'$ are orbits under the unipotent part of that intersection, hence finite; and the map is surjective because $\pi(R)=\calA'$.
\end{proof}

The following proposition is the crucial observation of this section: the Ribet property is forced by finiteness.

\begin{proposition}\label{prop_ribet_from_finiteness}
The pair $(\calA',\ttt)$ is a Ribet datum in the sense of Definition~\ref{def_ribet_datum}.
\end{proposition}
\begin{proof}
Let $\chi\in\Xi$ and write $f:R_\ttt\to\calA'$ for the finite surjective morphism of Lemma~\ref{lem_toric_part_of_R}. Pushing out along $\chi$ turns $\pt_{\Ag,n}|_{\calA',\gm^{n}/\ttt}$ into the $\gm$-torsor of $(\alpha,\beta_\chi)^{*}\pl_{\Ag}$ by \eqref{eq_pushout_along_chi}, and $R_\ttt$ into a closed subvariety of it, still finite and surjective over $\calA'$; we denote it again by $R_\ttt$. Consider the cartesian square
\[
\begin{tikzcd}[column sep=10ex]
\bigl((\alpha,\beta_\chi)^{*}\pl_{\Ag}\bigr)^{\times}\times_{\calA'}R_\ttt \arrow[d] \arrow[r] \arrow[dr,phantom,"\ulcorner",very near start]
& \bigl((\alpha,\beta_\chi)^{*}\pl_{\Ag}\bigr)^{\times} \arrow[d]\\
R_\ttt \arrow[r,"f"] \arrow[u,"\iota",dotted,bend left]
& \calA' \rlap{,}
\end{tikzcd}
\]
in which $\iota$ is the tautological section given by $(\id_{R_\ttt},R_\ttt\hookrightarrow(\cdots)^{\times})$. The existence of a section forces the $\gm$-torsor on the upper left to be trivial, so the line bundle $f^{*}\bigl((\alpha,\beta_\chi)^{*}\pl_{\Ag}\bigr)$ on $R_\ttt$ is trivial.

Let $x$ be a geometric point of $T$. Then $f_x:(R_\ttt)_x\to\calA'_x$ is finite and surjective and $f_x^{*}\bigl((\alpha,\beta_\chi)^{*}\pl_{\Ag}\bigr)_x$ is trivial, so by the projection formula \cite[Proposition 2.5(c)]{fulton}
\begin{align*}
  \deg(f_x)\cdot c_1\bigl(((\alpha,\beta_\chi)^{*}\pl_{\Ag})_x\bigr)\cap[\calA'_x]
  &=c_1\bigl(((\alpha,\beta_\chi)^{*}\pl_{\Ag})_x\bigr)\cap f_{x,*}[(R_\ttt)_x]\\
  &=f_{x,*}\Bigl(c_1\bigl(f_x^{*}((\alpha,\beta_\chi)^{*}\pl_{\Ag})_x\bigr)\cap[(R_\ttt)_x]\Bigr)=0 .
\end{align*}
The N\'eron--Severi group of an abelian variety is torsion free, so $((\alpha,\beta_\chi)^{*}\pl_{\Ag})_x$ lies in $\pic^{0}_{\calA'_x/x}(\kappa(x))$ for every geometric point $x$ of $T$; by the moduli interpretation of $\pic^{0}_{\calA'/T}$ this means
\[
  (\alpha,\beta_\chi)^{*}\pl_{\Ag}\in\pic^{0}_{\calA'/T}(T) .
\]
By Theorem~\ref{thm_ribet_property} the pair $(\alpha,\beta_\chi)$ satisfies the Ribet property. As $\chi\in\Xi$ was arbitrary, $(\calA',\ttt)$ is a Ribet datum.
\end{proof}

At this point we know that $R$ is finite over $\calA'$ after dividing by $\ttt$, and that $(\calA',\ttt)$ is a Ribet datum, hence that the subvariety $R(\calA',\ttt)$ of Definition~\ref{def_ribet_type} exists. This is not yet enough to identify the two: $R(\calA',\ttt)$ is a very particular finite cover of $\calA'$, the one cut out by the Ribet multi-section, whereas so far $R$ is an arbitrary one. The identification comes from computing the group $Q$, which the Ribet property just established now makes possible.

\begin{lemma}\label{lem_Q_computation}
There is $\underline z_0\in U(\QQ)$, unique modulo $U'(\QQ)$, such that for every $\QQ$-algebra $C$ the group $Q(C)$ consists of the matrices
\[
  \begin{pmatrix}
    \mu(g)1_n&(v^{\top}\beta_{i,\ZZ}\,g)_{1\leq i\leq n}&\underline z\\
    0&g&\alpha_{\ZZ}v\\
    0&0&1
  \end{pmatrix}
\]
with $v\in\M_{2d,1}(C)$ and $g\in G_{T}(C)$, subject to
\[
  \chi\underline z=\tfrac12v^{\top}\beta_{\chi,\ZZ}\alpha_{\ZZ}v+(\mu(g)-1)\,\chi\underline z_0
  \qquad\text{for every }\chi\in\Xi .
\]
Equivalently, $Q=\nu\,Q_{\Lambda,\ttt}\,\nu^{-1}$ for the element $\nu\in U(\QQ)$ with coordinate $-\underline z_0$.
\end{lemma}
\begin{proof}
Since $\pi(R)=\calA'$, the composite $Q\hookrightarrow P\to P/U$ has image $Q'_{\Lambda}$, and its kernel is $Q\cap U=U'$. By Notation~\ref{not_Q_lambda} the elements of $Q(C)$ are then exactly the matrices of the displayed shape, with $v\in\M_{2d,1}(C)$ and $g\in G_{T}(C)$ arbitrary and with $\underline z$ determined modulo $U'(R)$, the sequence $1\to U'\to Q\to Q'_{\Lambda}\to1$ being exact on $R$-points because $U'$ is a vector group; that is, for each $\chi\in\Xi$ the quantity $\chi\underline z$ is a function
\[
  f_\chi:G_{T}\times_{\QQ}\M_{2d,1}\longrightarrow\mathbb{A}^{1}_{\QQ},
\]
a morphism of schemes over $\QQ$, obtained by composing $Q'_{\Lambda}\cong Q/U'\hookrightarrow(P/U')$ with the $\chi$-component of the upper right corner. The group law of $Q$ gives, for $v,w\in\M_{2d,1}(C)$ and $g,h\in G_{T}(C)$,
\[
  \mu(g)f_\chi(w,h)+v^{\top}\beta_{\chi,\ZZ}\,g\,\alpha_{\ZZ}w+f_\chi(v,g)=f_\chi(x,gh),
\]
where $x\in\M_{2d,1}(C)$ satisfies
\[
  \alpha_{\ZZ}x=g\alpha_{\ZZ}w+\alpha_{\ZZ}v,
  \qquad
  x^{\top}\beta_{i,\ZZ}=\mu(g)w^{\top}\beta_{i,\ZZ}g^{-1}+v^{\top}\beta_{i,\ZZ}\quad(1\leq i\leq n).
\]
Specializing gives the three identities
\begin{enumerate}
  \item $\mu(g)f_\chi(0,h)+f_\chi(0,g)=f_\chi(0,gh)$;
  \item $f_\chi(0,h)+f_\chi(v,1)=f_\chi(v,h)$;
  \item $f_\chi(w,1)+v^{\top}\beta_{\chi,\ZZ}\alpha_{\ZZ}w+f_\chi(v,1)=f_\chi(w+v,1)$.
\end{enumerate}
As $G_{T}$ contains the scalars, substituting $g=2\cdot\id$ in (1) yields
\[
  4f_\chi(0,h)+f_\chi(0,2\cdot\id)=f_\chi(0,2h)=f_\chi(h\cdot2)=\mu(h)f_\chi(0,2\cdot\id)+f_\chi(0,h),
\]
so with $z_\chi:=f_\chi(0,2\cdot\id)/3\in\QQ$ we get $f_\chi(0,h)=(\mu(h)-1)z_\chi$. Next put $\varphi_\chi(v):=f_\chi(v,1)-\tfrac12v^{\top}\beta_{\chi,\ZZ}\alpha_{\ZZ}v$. By Proposition~\ref{prop_ribet_from_finiteness} and Proposition~\ref{prop_ribet_symmetric} the matrix $\beta_{\chi,\ZZ}\alpha_{\ZZ}$ is symmetric, so (3) is equivalent to the additivity $\varphi_\chi(v)+\varphi_\chi(w)=\varphi_\chi(v+w)$, whence $\varphi_\chi(v)=m_\chi v$ for some $m_\chi\in\M_{1,2d}(\QQ)$. Combining with (2),
\begin{equation}\label{eq_form_of_f}
  f_\chi(v,g)=\tfrac12v^{\top}\beta_{\chi,\ZZ}\alpha_{\ZZ}v+m_\chi v+(\mu(g)-1)z_\chi .
\end{equation}
It remains to see that $m_\chi=0$. Substituting \eqref{eq_form_of_f} into the functional equation and using the symmetry of $\beta_{\chi,\ZZ}\alpha_{\ZZ}$ once more leaves
\[
  m_\chi(x-v-w)=0\qquad\text{for all }v,w\in\M_{2d,1}(C)\text{ and }g\in G_{T}(C).
\]
Take $g=2\cdot\id$, so that $\mu(g)=4$ and $g^{-1}=\tfrac12$. The two equations determining $x$ become
\[
  \alpha_{\ZZ}x=2\alpha_{\ZZ}w+\alpha_{\ZZ}v=\alpha_{\ZZ}(v+2w),
  \qquad
  x^{\top}\beta_{i,\ZZ}=2w^{\top}\beta_{i,\ZZ}+v^{\top}\beta_{i,\ZZ}=(v+2w)^{\top}\beta_{i,\ZZ}
\]
for every $i$, and $v\mapsto\bigl(\alpha_{\ZZ}v,(v^{\top}\beta_{i,\ZZ})_i\bigr)$ is injective, being the inclusion of $\Lambda$; hence $x=v+2w$. So $x-v-w=w$ and $m_\chi w=0$ for every $w$, that is, $m_\chi=0$.

Finally, $\chi\mapsto z_\chi$ is additive, being so in each of the identities above, and vanishes on the characters trivial on $\gm^{n}/\ttt$; so $z_\chi=\chi\underline z_0$ for some $\underline z_0\in U(\QQ)$, unique modulo $U'(\QQ)$. Comparing with Proposition~\ref{prop_Q_points} and computing
\[
  u(\underline z_0)\begin{pmatrix}
    \mu(g)1_n&\underline x&\underline z\\ 0&g&y\\ 0&0&1
  \end{pmatrix}u(\underline z_0)^{-1}
  =\begin{pmatrix}
    \mu(g)1_n&\underline x&\underline z+(1-\mu(g))\underline z_0\\ 0&g&y\\ 0&0&1
  \end{pmatrix},
\]
where $u(\underline z_0)\in U(\QQ)$ is the element with coordinate $\underline z_0$, gives the last assertion.
\end{proof}

\begin{remark}\label{rem_group_theoretic_route}
There is a purely group-theoretic route to the results of this section, which bypasses the geometry entirely: one classifies the special algebraic subgroups of $P$ directly, and the Ribet conditions then appear as the symmetry of the bilinear form induced on the weight $-1$ part with values in the weight $-2$ part, exactly as in Proposition~\ref{prop_ribet_symmetric}. This is carried out in the thesis \cite{Lopuhaa14}, for the ambient group in general and, for $g=1$, down to the special subvarieties of the Poincar\'e torsor itself; see also \cite{EdixhovenPrivate}. We have preferred the geometric route because it does not presuppose the coordinates of Subsection~\ref{subsec_partial_rank_MSD}.
\end{remark}

\subsection{Step 3: Clear Torsion Translations in the Toric Part}\label{subsec_finish}

\begin{proposition}\label{prop_step3}
Assume, as we may by Subsection~\ref{subsec_moving}, that $\pi(R)=\calA'$ is an abelian subscheme over $T$. Then $(\calA',\ttt)$ is a Ribet datum over $T$ and $R$ is a torsion translate of $R(\calA',\ttt)$; in particular $R$ is of Ribet type.
\end{proposition}
\begin{proof}
The first assertion is Proposition~\ref{prop_ribet_from_finiteness}. By Lemma~\ref{lem_Q_computation} we have $Q=\nu Q_{\Lambda,\ttt}\nu^{-1}$ for an element $\nu\in U(\QQ)$, and $D_Q=\nu D_{\Lambda,\ttt}$, both data having the same unipotent intersection $U'$ with $U$. By Proposition~\ref{prop_torsion_conjugation} the subvariety attached to $(Q,D_Q)$ is therefore the translate of the one attached to $(Q_{\Lambda,\ttt},D_{\Lambda,\ttt})$ by the torsion section of the fibral $\gm^{n}$ given by the class of $\nu$ in $U(\QQ)/U(\ZZ)$, and by Proposition~\ref{prop_comparison} the latter subvariety is $R(\calA',\ttt)$.
\end{proof}

\begin{theorem}\label{thm_classification}
Let $R\subseteq\pgn$ be an irreducible closed subvariety. Then $R$ is a special subvariety of $\pgn$, viewed as a mixed Shimura variety, if and only if $R$ is a subvariety of Ribet type in the sense of Definition~\ref{def_ribet_type}, for $\calA=\Ag$ over a special subvariety of $\ag$.
\end{theorem}
\begin{proof}
If $R$ is of Ribet type over a special subvariety of $\ag$, it is special by Theorem~\ref{thm_ribet_type_special}. Conversely, let $R$ be special, with $T=\Pi(R)$, a special subvariety of $\ag$. By Proposition~\ref{prop_step1} some torsion translate $R_1$ of $R$ has $\pi(R_1)$ an abelian subscheme of $\calB_{g,n}|_{T}$ over $T$, and $R$ is a torsion translate of $R_1$. By Proposition~\ref{prop_step3} the subvariety $R_1$ is a torsion translate of $R(\calA',\ttt)$ for a Ribet datum $(\calA',\ttt)$ over $T$. Torsion translation being transitive up to composition, $R$ is a torsion translate of $R(\calA',\ttt)$, hence of Ribet type over the special base $T$.
\end{proof}

\section{The Tales of Zilber--Pink and Relative Manin--Mumford}\label{section7}

We now apply the classification theorem of Section~\ref{section6}. The issue is to pass from the universal mixed Shimura variety $\pgn$ to an arbitrary semi-abelian variety, or to an arbitrary semi-abelian scheme. We obtain two consequences. First, in Subsection~\ref{subsec_ZPimpliesZPSA}, we show that the Zilber--Pink conjecture for mixed Shimura varieties implies its semi-abelian counterpart, correcting \cite[Theorem 5.7]{Pink05}. Second, in Subsection~\ref{subsec_reform_RMM}, we give a corrected formulation of the relative Manin--Mumford conjecture for semi-abelian schemes, since \cite[Conjecture 6.1]{Pink05} is false as stated, and deduce this formulation from Zilber--Pink.

Both arguments have the same basic form. The given semi-abelian object is realized as a fiber, or as a pullback, of the universal family $\pgn\to(\Agv)^{[n]}$. One then needs two comparison results between the family and its fibers. Theorem~\ref{thm_ribet_type_fibers} describes the restriction of a special subvariety to a fiber: it is a finite union of translates of one algebraic subgroup. Conversely, Lemma~\ref{lem_extension} spreads an algebraic subgroup of a fiber, and any torsion translate of it, to a special subvariety of the family with the same codimension.

\subsection{The Zilber--Pink Conjecture for Mixed Shimura Varieties}\label{subsec_ZP}

All mixed Shimura varieties in this article are connected mixed Shimura varieties; we omit the adjective from the notation. Thus, for a mixed Shimura variety $X$ over $\CC$, a \textbf{special subvariety} of $X$ means an irreducible component of a mixed Shimura subvariety of $X$.

The Zilber--Pink conjecture for mixed Shimura varieties concerns unlikely intersections between an irreducible closed subvariety and special subvarieties of a mixed Shimura variety. It implies the Mordell--Lang, Andr\'e--Oort, and Andr\'e--Pink conjectures \cite{pink04}. We recall the formulation \cite[Conjecture 1.3]{Pink05}, which is equivalent to \cite[Conjectures 1.1 and 1.2]{Pink05}. A subvariety $Z\subseteq X$ is \textbf{Hodge generic} if it is not contained in any proper special subvariety of $X$.

\begin{conjecture}[{\cite[Conjecture 1.3]{Pink05}}]\label{conj_ZP}
Consider a mixed Shimura variety $X$ over $\CC$ and a Hodge generic irreducible closed subvariety $Z$. Then the union
\[
  \bigcup_{\substack{\text{special subvariety } R\subseteq X\\ \operatorname{codim}_X R>\dim Z}}R\cap Z
\]
is not Zariski dense in $Z$.
\end{conjecture}

For semi-abelian varieties, the corresponding special subvarieties are the torsion translates of algebraic subgroups. This gives the semi-abelian Zilber--Pink conjecture \cite[Conjecture 5.1]{Pink05}, recalled below as Conjecture~\ref{conj_ZP_SA}. Pink also proposed a relative version, \cite[Conjecture 6.1]{Pink05}, recalled as Conjecture~\ref{conj_RMM_pink}.

One expects Conjecture~\ref{conj_ZP} to imply the semi-abelian conjecture and its relative analogue. Pink claimed these deductions in \cite{Pink05}. As observed in \cite{BE20}, however, the arguments in \cite[Theorem 5.7, Theorem 6.3]{Pink05} are not correct; moreover, \cite[Conjecture 6.1]{Pink05} itself is contradicted by Ribet sections.

The correction is based on Theorem~\ref{thm_classification}. We prove that Conjecture~\ref{conj_ZP} implies \cite[Conjecture 5.1]{Pink05}, and hence \cite[Conjecture 5.2]{Pink05}. We then formulate a corrected relative conjecture and prove that it also follows from Conjecture~\ref{conj_ZP}.

\subsection{The Zilber--Pink Conjecture for Semi-abelian Varieties}\label{subsec_ZP_semiab}

We first recall the semi-abelian form of the conjecture and fix the elementary group-theoretic notation used in its proof.
\begin{notation}\label{not_codim_union}
For a commutative algebraic group $G$ over $\CC$ we set $G^{[>d]}$ to be the union of all algebraic subgroups of $G$ of codimension $>d$, that is,
\[
  G^{[>d]}:=\bigcup_{\substack{\text{algebraic subgroup } H\subseteq G\\ \operatorname{codim}_G H>d}}H .
\]
More generally, for an algebraic family of commutative algebraic groups $\calG\to S$ over a variety $S$ over $\CC$ we set $\calG^{[>d]}:=\bigcup_{s\in S}\calG_s^{[>d]}$.
\end{notation}

Pink discussed the following analogues of the Manin--Mumford and Mordell--Lang conjectures in \cite{Pink05}.

\begin{conjecture}[Zilber--Pink for semi-abelian varieties, {\cite[Conjecture 5.1]{Pink05}}]\label{conj_ZP_SA}
Consider a semi-abelian variety $G$ over $\CC$ and an irreducible closed subvariety $Z$ with $\dim Z=d$. Assume that $Z$ is \textbf{not contained in any proper algebraic subgroup} of $G$. Then $Z\cap G^{[>d]}$ is not Zariski dense in $Z$.
\end{conjecture}

By \cite[Theorem 5.3, Theorem 5.5]{Pink05}, Conjecture~\ref{conj_ZP_SA} is equivalent to the following.

\begin{conjecture}[{\cite[Conjecture 5.2]{Pink05}}]\label{conj_ZP_ML}
Consider a semi-abelian variety $G$ over $\CC$, a subgroup $\Gamma\subseteq G(\CC)$ of finite $\ZZ$-rank, and an irreducible closed subvariety $Z$ with $\dim Z=d$. Assume that $Z$ is \textbf{not contained in any translate of a proper algebraic subgroup} of $G$. Then $Z\cap(\Gamma+G^{[>d]})$ is not Zariski dense in $Z$.
\end{conjecture}

We use the following elementary facts about algebraic subgroups of semi-abelian varieties.

\begin{lemma}\label{lem_subgroup_structure}
Let $G$ be a semi-abelian variety over $\CC$ and let $H\subseteq G$ be an algebraic subgroup. Then:
\begin{enumerate}[label=\textup{(\arabic*)}]
  \item the identity component $H^{\circ}$ is a semi-abelian subvariety of $G$;
  \item $H/H^{\circ}$ is finite and $H$ is a finite union of translates of $H^{\circ}$ by torsion points of $G$;
  \item the set of torsion points of $G$ is exactly $G^{[>\dim G-1]}$.
\end{enumerate}
\end{lemma}
\begin{proof}
(1) $H^{\circ}$ is a connected commutative algebraic group, so by Chevalley's theorem it is an extension of an abelian variety by a connected linear group $L$. Every homomorphism from a connected linear group to an abelian variety is trivial, so $L$ is contained in the toric part $\gm^{n}$ of $G$; being connected it is a subtorus. Hence $H^{\circ}$ has no unipotent part and is semi-abelian.

(2) An algebraic group has finitely many connected components, so $H=\bigcup_{j}(H^{\circ}+t_j)$ for finitely many $t_j\in H(\CC)$. Let $m:=|H/H^{\circ}|$, so that $mt_j\in H^{\circ}$. A semi-abelian variety over an algebraically closed field is divisible, so $mt_j=mu_j$ for some $u_j\in H^{\circ}$, and then $t_j-u_j$ is killed by $m$ while $H^{\circ}+t_j=H^{\circ}+(t_j-u_j)$.

(3) By definition $G^{[>\dim G-1]}$ is the union of the algebraic subgroups of $G$ of dimension $0$, that is, of the finite subgroups of $G$; and a point of $G$ lies in a finite subgroup if and only if it is torsion.
\end{proof}

\begin{remark}\label{rem_torsion_in_codim}
For a closed subvariety $Z\subsetneq G$ one has $d=\dim Z<\dim G$, so by Lemma~\ref{lem_subgroup_structure}(3)
\[
  \left\{\text{torsion points of }G\right\}=G^{[>\dim G-1]}\ \subseteq\ G^{[>d]} .
\]
Thus semi-abelian Zilber--Pink implies Manin--Mumford.
\end{remark}
The hypotheses in Conjectures~\ref{conj_ZP_SA} and~\ref{conj_ZP_ML} differ slightly. The Mordell--Lang form excludes translates of proper algebraic subgroups, whereas Conjecture~\ref{conj_ZP_SA} excludes only proper algebraic subgroups. The latter condition is weaker: a coset $H+c$ of a proper subgroup may itself be contained in no proper algebraic subgroup, equivalently $H+\overline{\langle c\rangle}=G$.

\subsection{ZP for Mixed Shimura Varieties Implies ZP for Semi-abelian Varieties}\label{subsec_ZPimpliesZPSA}

The deduction from Conjecture~\ref{conj_ZP} requires a moduli-theoretic realization of a given semi-abelian variety. We use the following universal families.

\begin{definition}\label{def_shimura_family}
Let $S_1\subseteq(\Agv)^{[n]}$ be a special subvariety. The \textbf{Shimura family of semi-abelian varieties} over $S_1$ is the restriction
\[
  \pgn|_{S_1}:=\pgn\times_{(\Agv)^{[n]}}S_1\ \longrightarrow\ S_1 .
\]
It is again a mixed Shimura variety, being the preimage of $S_1$ under the Shimura morphism $\pgn\to(\Agv)^{[n]}$; in particular its special subvarieties are exactly the special subvarieties of $\pgn$ contained in it, and $\dim\pgn|_{S_1}=\dim S_1+g+n$.
\end{definition}

\begin{remark}\label{rem_minimal_family}
Every semi-abelian variety over $\CC$ occurs in a Shimura family, and in a smallest one. Explicitly, let $G$ be a semi-abelian variety over $\CC$ of toric rank $n$ and put $g:=\dim G-n$. Here and below we suppress a suitable level structure and a suitable polarization type in the definition of $\pgn$. With these choices made, Proposition~\ref{prop_universal_semiabelian} gives a point $\underline b\in(\Agv)^{[n]}(\CC)$ and an isomorphism $G\cong(\pgn)_{\underline b}$. Taking for $S_1$ the \textbf{special closure} of $\underline b$, that is, the smallest special subvariety of $(\Agv)^{[n]}$ containing $\underline b$, we obtain a Shimura family $\pgn|_{S_1}$ containing $G$ as the fiber at a Hodge generic point.
\end{remark}

\begin{remark}\label{rem_coverings}
We record the convention on finite coverings used below. Specialty is not invariant under arbitrary finite base change: a finite-index subgroup of an arithmetic group need not be congruence, so an arbitrary finite covering of a Shimura variety need not be a Shimura variety. The coverings used in Section~\ref{section4}, in the definitions of torsion translation and Ribet type, are internal to the construction: the resulting subvarieties live over the original base. By contrast, the coverings in Lemmas~\ref{lem_hodge_locus} and~\ref{lem_extension} are obtained by raising the level. They are finite Shimura morphisms; the covering spaces are Shimura varieties, images and preimages of special subvarieties are special, and codimensions are preserved.
\end{remark}

As in Remark~\ref{rem_ribet_type_extremes}, \emph{special} is reserved for the mixed Shimura setting, and subvarieties of a semi-abelian scheme are called \emph{of Ribet type}, in the sense of Definition~\ref{def_ribet_type_general}; a single semi-abelian variety is the case $S=\spec\CC$.

The next lemma spreads an abelian subvariety of one fiber to an abelian subscheme over the special closure of the moduli point.

\begin{lemma}\label{lem_hodge_locus}
Let $s\in\ag(\CC)$, let $A:=(\Ag)_s$, and let $S_0\subseteq\ag$ be the special closure of $s$. For every abelian subvariety $X\subseteq A$ there exist a finite surjective Shimura morphism $S_0'\to S_0$, obtained by raising the level, and an abelian subscheme $\calX\subseteq\Ag\times_{S_0}S_0'$ such that $S_0'$ is again a special subvariety of a Siegel moduli variety, and the fiber of $\calX$ over a point of $S_0'$ above $s$ is $X$.
\end{lemma}
\begin{proof}
Let $V_{\ZZ}:=\mathrm{H}_1(A,\ZZ)$ and $V_{\QQ}:=\mathrm{H}_1(A,\QQ)$, a polarizable $\QQ$-Hodge structure of type $\{(-1,0),(0,-1)\}$. Let $\underline V_{\ZZ}$ and $\underline V_{\QQ}$ be the restrictions to $S_0$ of the corresponding universal VHS, and let $\pi\in\End(V_{\QQ})$ be the projector onto $W_{\QQ}:=\mathrm{H}_1(X,\QQ)$. Since $W_{\QQ}$ is a sub-Hodge structure of $V_{\QQ}$, the tensor $\pi$ is a Hodge class of $\End(\underline V_{\QQ})$ at $s$.

Since $s$ is Hodge generic in $S_0$, it is Hodge generic for $\underline V_{\QQ}$ as well; equivalently, $\mathrm{MT}(V_{\QQ})$ is the generic Mumford--Tate group of $\underline V_{\QQ}$, by \cite[2.9]{Moonen98}. Hence the Hodge tensor $\pi$ extends over the universal cover of $S_0$ to a flat section of $\End(\underline V_{\QQ})$ which is a Hodge class at every point, by \cite[\S2.3]{MoonenOort}. The monodromy acts on these Hodge classes through a finite group by \cite[Lemma 4]{Andre92} and its proof, so after a finite covering of $S_0$ the tensor $\pi$ extends to such a flat section.

This covering may be taken to be congruence. Choose $k\geq1$ with $k\pi\in\End(V_{\ZZ})$. The monodromy orbit of $k\pi$ is finite; an element of the principal congruence subgroup of level $N$ fixes $k\pi$ modulo $N$, and hence fixes $k\pi$ once $N$ exceeds the entries occurring in that finite orbit. The covering of $S_0$ attached to that congruence subgroup is a special subvariety of the Siegel moduli variety of level $N$, finite and surjective over $S_0$, and we take it for $S_0'$.

The image of the extended projector is a sub-VHS $\underline W_{\QQ}\subseteq\underline V_{\QQ}$ with fiber $W_{\QQ}$ above $s$ and with the same Hodge type. Then $\underline W_{\QQ}\cap\underline V_{\ZZ}$ is a polarizable sub-$\ZZ$-VHS of that type, and therefore corresponds to an abelian subscheme $\calX\subseteq\Ag\times_{S_0}S_0'$ whose fiber over a point above $s$ is $X$.
\end{proof}

We also need the corresponding semi-abelian statement: algebraic subgroups of a fiber, together with their torsion translates, spread to special subvarieties of the family without changing codimension.

\begin{lemma}\label{lem_extension}
Let $\underline b\in(\Agv)^{[n]}(\CC)$, let $S_1$ be the special closure of $\underline b$ in $(\Agv)^{[n]}$, and put $E:=(\pgn)_{\underline b}$. Let $H\subseteq E$ be an algebraic subgroup. Then there exist a finite surjective Shimura morphism $S_1'\to S_1$, obtained by raising the level, and a point $\underline b'\in S_1'(\CC)$ above $\underline b$ with the following properties. The variety $S_1'$ is again a special subvariety of a space of the form $(\Agv)^{[n]}$, the restriction $\pgn|_{S_1'}:=\pgn\times_{(\Agv)^{[n]}}S_1'$ is a mixed Shimura variety, the fiber of $\pgn|_{S_1'}$ at $\underline b'$ is $E$, and:
\begin{enumerate}
  \item there is a special subvariety $T\subseteq\pgn|_{S_1'}$ satisfying
  \[
    T_{\underline b'}=H^{\circ},\qquad
    \operatorname{codim}_{\pgn|_{S_1'}}T=\operatorname{codim}_EH^{\circ};
  \]
  \item writing $H=\bigcup_{j=1}^{m}(H^{\circ}+t_j)$ with $t_j$ torsion as in Lemma~\ref{lem_subgroup_structure}(2), there are special subvarieties $T_1,\dots,T_m\subseteq\pgn|_{S_1'}$, each of codimension $\operatorname{codim}_EH$ in $\pgn|_{S_1'}$, with $H^{\circ}+t_j\subseteq(T_j)_{\underline b'}$ for every $j$.
\end{enumerate}
\end{lemma}
\begin{proof}
We prove (1). Write $s$ for the image of $\underline b$ in $\ag$, $A:=(\Ag)_s$, and let $\ttt\subseteq\gm^{n}$ be the toric part of $H^{\circ}$ and $\iota:X\hookrightarrow A$ its abelian quotient, so that $H^{\circ}$ is an extension of $X$ by $\ttt$ inside $E$. Pushing that extension out along $q_\ttt:\gm^{n}\to\gm^{n}/\ttt$ splits it, so by Theorem~\ref{thm_weil_barsotti}
\begin{equation}\label{eq_split_condition}
  \iota^{\vee}b_\chi=0\qquad\text{for every }\chi\in\chl(\gm^{n}/\ttt),
\end{equation}
exactly as in \eqref{eq_flat_subgroup_class}.

Let $S_0$ be the special closure of $s$ in $\ag$. By Lemma~\ref{lem_hodge_locus} there exist a finite surjective morphism $S_0'\to S_0$, a point $s'\in S_0'(\CC)$ above $s$, and an abelian subscheme $\calX\subseteq\Ag\times_{S_0}S_0'$ whose fiber at $s'$ is $X$. We write again $\iota:\calX\hookrightarrow\Ag\times_{S_0}S_0'$ for its inclusion, which extends $\iota$, and $\iota^{\vee}$ for the dual homomorphism $\Agv\times_{S_0}S_0'\to\calX^{\vee}$.

Let $\underline b'$ be a point of $(\Agv)^{[n]}\times_{S_0}S_0'$ above $\underline b$ and over $s'$, and let $S_1'$ be its special closure. The image of $S_1'$ in $(\Agv)^{[n]}$ is a special subvariety containing $\underline b$, hence contains $S_1$. On the other hand, $S_1'$ is contained in the preimage of $S_1$, a special subvariety of $(\Agv)^{[n]}\times_{S_0}S_0'$ containing $\underline b'$. Thus $S_1'\to S_1$ is finite and surjective, and the fiber of $\pgn|_{S_1'}$ at $\underline b'$ is $E$.

Let $\calZ\subseteq(\Agv)^{[n]}\times_{S_0}S_0'$ be the subgroup scheme cut out by the conditions $\iota^{\vee}c_\chi=0$, $\chi\in\chl(\gm^{n}/\ttt)$, as in the proof of Proposition~\ref{prop_flat_subgroup}. Its identity component is an abelian subscheme and its component group is finite. Hence every irreducible component of $\calZ$ is a torsion translate of an abelian subscheme over $S_0'$, and is therefore special: the space $(\Agv)^{[n]}$ has no toric part, so its special subvarieties are exactly the torsion translates of abelian subschemes over special subvarieties of $\ag$, by \cite[Proposition 1.1]{gaoziyang} applied to each factor, and the same holds over the finite covering $S_0'$.

By \eqref{eq_split_condition} the point $\underline b'$ lies on $\calZ$. The component of $\calZ$ through $\underline b'$ is therefore a special subvariety containing $\underline b'$, and so contains $S_1'$. Thus
\[
  \iota^{\vee}c_\chi=0\qquad\text{for every }\chi\in\chl(\gm^{n}/\ttt)\text{ and every point }\underline c\text{ of }S_1' .
\]
The preceding display says that the extension $(q_\ttt)_{*}\bigl(\pgn|_{\calX}\bigr)$ of $\calX$ by $\gm^{n}/\ttt$ is split over the whole of $S_1'$. Its splitting is unique, since two splittings differ by an element of $\Hom_{S_1'}(\calX,\gm^{n}/\ttt)=0$ (Lemma~\ref{lem_noconst}). Let $\calE\subseteq\pgn|_{S_1'}$ be the preimage of the image of this splitting. Then $\calE$ is a subgroup scheme of $\pgn|_{S_1'}$, flat over $S_1'$, and an extension of $\calX$ by $\ttt$. By the same uniqueness, $\calE_{\underline b'}$ is the unique subgroup of $E$ with toric part $\ttt$ and abelian quotient $X$, namely $H^{\circ}$.

By Proposition~\ref{prop_flat_subgroup} the subgroup scheme $\calE$ is of Ribet type, hence a special subvariety of $\pgn|_{S_1'}$ by Theorem~\ref{thm_classification}. Finally
\[
  \operatorname{codim}_{\pgn|_{S_1'}}\calE=\bigl(\dim S_1'+g+n\bigr)-\bigl(\dim S_1'+\dim H^{\circ}\bigr)=\operatorname{codim}_E H^{\circ},
\]
which is the codimension comparison asserted in (1); it is here that the flatness of $\calE$ over the whole of $S_1'$, rather than merely the existence of $\calE_{\underline b'}$, is used.

For (2) it suffices, given a torsion $t\in E(\CC)$, to produce a special subvariety of $\pgn|_{S_1'}$ of the same dimension as $\calE$ whose fiber at $\underline b'$ contains $H^{\circ}+t$; applying this to $t_1,\dots,t_m$ then gives the assertion, the codimensions of $H$ and of $H^{\circ}$ in $E$ being equal.

Let $t\in E(\CC)$ be a torsion point and let $a_0\in A$ be its image. Then $a_0$ is torsion, say $N^{2}a_0=0$. Since $\Ag[N^{2}]\to\ag$ is finite \'etale, $a_0$ extends, after enlarging the covering $S_0'$ and hence $S_1'$, to a torsion section $(a_0,0)$ of $\calB_{g,n}[N^{2}]$ over $S_1'$. By Lemma~\ref{lem_torsion_lift} this section lifts to a torsion element $\tau\in\calH_N$.

The torsion translate $\calE+\tau$ is again of Ribet type, hence special, and has the same dimension as $\calE$, since torsion translation is defined through finite morphisms and an isomorphism. By Proposition~\ref{prop_fibral_transport}(2),(3), the fiber $(\calE+\tau)_{\underline b'}$ is the union, over the relevant preimages $\underline c$ of $\underline b'$ under multiplication by $N$, of the images of $(\calE_N)_{\underline c}$ under a homomorphism followed by translation by the point $u=\tau(0)$, which lies over $a_0$; here $\calE_N$ is the irreducible component of $\varpi_N^{-1}(\calE)$ used in the torsion translation.

The torsion lifts of $(a_0,0)$ in $\calH_N$ form a torsor under the roots of unity of $\gm^{n}$. The torsion points of $E$ above $a_0$ form a nonempty coset of the same group: nonemptiness follows from the divisibility of $\gm^{n}(\CC)$, which allows a point above $a_0$ to be corrected by a root of unity so as to be killed by the order of $a_0$. Since $\gm^{n}$ is central in $\calH_N$, the assignment $\tau\mapsto u=\tau(0)$ is equivariant for the two actions of $\mu_{\infty}^{n}$, hence surjective onto the torsion points above $a_0$. Thus $\tau$ may be chosen so that $u=t$, and then $(\calE+\tau)_{\underline b'}\supseteq H^{\circ}+t$.
\end{proof}

We now prove the semi-abelian Zilber--Pink conjecture from Conjecture~\ref{conj_ZP}.

\begin{theorem}[Zilber--Pink for mixed Shimura varieties implies Zilber--Pink for semi-abelian varieties]\label{thm_ZP_implies_ZPSA}
Conjecture~\ref{conj_ZP} implies Conjecture~\ref{conj_ZP_SA}, and hence Conjecture~\ref{conj_ZP_ML}.
\end{theorem}
\begin{proof}
By \cite[Theorem 5.3, Theorem 5.5]{Pink05}, it is enough to prove Conjecture~\ref{conj_ZP_SA}. Let $G$ and $Z\subseteq G$ be as in that conjecture, and put $d:=\dim Z$. Thus $Z$ is not contained in any proper algebraic subgroup of $G$.

\smallskip
\noindent\textit{Step 1: realization in a Shimura family.} Let $n$ be the toric rank of $G$ and put $g:=\dim G-n$. By Remark~\ref{rem_minimal_family}, there exist a point $\underline b\in(\Agv)^{[n]}(\CC)$ and an isomorphism $G\cong(\pgn)_{\underline b}$. Let $S_1$ be the special closure of $\underline b$; then $\underline b$ is Hodge generic in $S_1$. Put $X:=\pgn|_{S_1}$, and regard $Z$ as a subvariety of the fiber $G\subseteq X$.

\smallskip
\noindent\textit{Step 2: the special closure of $Z$.} Let $R_0\subseteq X$ be the special closure of $Z$. Then $Z$ is Hodge generic in the mixed Shimura variety $R_0$, and Conjecture~\ref{conj_ZP}, applied inside $R_0$, implies that
\[
  \bigcup_{\substack{\text{special }Y\subseteq R_0\\ \operatorname{codim}_{R_0}Y>d}}\!\!Y\cap Z
\]
is not Zariski dense in $Z$. It therefore suffices to show the inclusion
\[
  Z\cap G^{[>d]}\ \subseteq\!\!\bigcup_{\substack{\text{special }Y\subseteq R_0\\ \operatorname{codim}_{R_0}Y>d}}\!\!Y\cap Z ,
\]
that is, that every point of $Z\cap G^{[>d]}$ lies on a special subvariety of $R_0$ of codimension $>d$.

Since $\pr(R_0)$ is a special subvariety of $S_1$ containing $\underline b$, and $\underline b$ is Hodge generic in $S_1$, we have $\pr(R_0)=S_1$. By Theorem~\ref{thm_classification} the subvariety $R_0$ is of Ribet type, so Theorem~\ref{thm_ribet_type_fibers} provides an algebraic subgroup $H\subseteq G$ of which every irreducible component of $(R_0)_{\underline b}$ is a translate, all of them of dimension $\dim H$. Let $F$ be the component containing the irreducible $Z$, and write $F=H+c_1$. By Corollary~\ref{cor_fiber_dim},
\begin{equation}\label{eq_R0_dim}
  \dim R_0=\dim S_1+\dim H .
\end{equation}

\smallskip
\noindent\textit{Step 3: the codimension estimate.} Let $p\in Z\cap G^{[>d]}$. Then $p$ lies in an algebraic subgroup $H_1\subseteq G$ with $\operatorname{codim}_GH_1>d$. By Lemma~\ref{lem_subgroup_structure}(2), we may write $p\in H_1^{\circ}+t$ with $t$ torsion, and $\operatorname{codim}_GH_1^{\circ}=\operatorname{codim}_GH_1>d$.

Lemma~\ref{lem_extension}(2) realizes this torsion translate in a Shimura family. After a finite surjective morphism $S_1'\to S_1$, there exist a point $\underline b'$ above $\underline b$ and a special subvariety $T'\subseteq\pgn|_{S_1'}$ of codimension $\operatorname{codim}_GH_1^{\circ}$ such that $H_1^{\circ}+t\subseteq T'_{\underline b'}$. Let $T$ be the image of $T'$ in $X$. Since $S_1'\to S_1$ induces a finite surjective Shimura morphism $\pgn|_{S_1'}\to X$, the image $T$ is special and has the same codimension in $X$; hence
\[
  H_1^{\circ}+t\subseteq T_{\underline b},\qquad \operatorname{codim}_XT=\operatorname{codim}_GH_1^{\circ}>d .
\]

We determine the fiber of $T$ at $\underline b$. Since $\pr(T)$ is special and contains the Hodge generic point $\underline b$, it is equal to $S_1$. Hence Corollary~\ref{cor_fiber_dim}, applied to $T$, gives
\[
  \dim T_{\underline b}=\dim T-\dim S_1=\dim H_1^{\circ}.
\]
By Theorem~\ref{thm_classification} and Theorem~\ref{thm_ribet_type_fibers}, the irreducible components of $T_{\underline b}$ are translates of one algebraic subgroup. The component containing the irreducible coset $H_1^{\circ}+t$ has the same dimension as that coset, hence is equal to it. Therefore $T_{\underline b}$ is a finite union of translates of $H_1^{\circ}$.

Let $W$ be the irreducible component of $T\cap R_0$ containing $p$. It is a special subvariety of $R_0$, since irreducible components of intersections of special subvarieties are special. We prove that $\operatorname{codim}_{R_0}W>d$.

As above, $\pr(W)$ is special and contains $\underline b$, so $\pr(W)=S_1$. Corollary~\ref{cor_fiber_dim}, applied to $W$ and to $R_0$, gives
\[
  \operatorname{codim}_{R_0}W
  =(\dim S_1+\dim H)-(\dim S_1+\dim W_{\underline b})
  =\dim H-\dim W_{\underline b}.
\]
Let $U$ be the irreducible component of $W_{\underline b}$ containing $p$. Since $W$ is special, Theorem~\ref{thm_classification} and Theorem~\ref{thm_ribet_type_fibers} imply that all irreducible components of $W_{\underline b}$ have the same dimension. Thus $\dim W_{\underline b}=\dim U$, and it remains to bound $\dim U$ in the fiber $G$.

The fiber $(R_0)_{\underline b}$ is a finite union of translates of $H$, one of which is $F=H+c_1$. This coset contains $Z$, and therefore lies in no proper algebraic subgroup of $G$. The fiber $T_{\underline b}$ is a finite union of translates of $H_1^\circ$, and the translate $H_1^\circ+t$ is one of them; moreover $(H_1^\circ+t)\cap F$ contains $p$.

Put $M:=H+H_1^\circ$. Since $(H_1^\circ+t)\cap(H+c_1)$ is nonempty, we have $c_1-t\in M$, whence $H+c_1\subseteq M+t$. Thus $F$ is contained in the algebraic subgroup generated by $M$ and $t$, namely $M+\langle t\rangle$. The point $t$ is torsion, so this subgroup contains $M$ with finite index.

Since $F$ is contained in no proper algebraic subgroup of $G$, we get $M+\langle t\rangle=G$. As $G$ is connected, this forces $M=G$.

Now $U\subseteq T_{\underline b}\cap(R_0)_{\underline b}$. This intersection is a finite union of intersections of translates of $H$ with translates of $H_1^\circ$; every nonempty such intersection is a coset of $H\cap H_1^\circ$. Hence $\dim U\leq\dim(H\cap H_1^\circ)$. Therefore
\[
\begin{aligned}
  \operatorname{codim}_{R_0}W
    &=\dim H-\dim U\\
    &\geq\dim H-\dim(H\cap H_1^\circ)\\
    &=\dim(H+H_1^\circ)-\dim H_1^\circ\\
    &=\operatorname{codim}_GH_1^{\circ}>d .
\end{aligned}
\]
This proves the required inclusion, and therefore the theorem.
\end{proof}

\begin{remark}\label{rem_pink_route}
The proof follows the general strategy of \cite{Pink05}: one first deduces Conjecture~\ref{conj_ZP_SA} from Conjecture~\ref{conj_ZP}, and then obtains Conjecture~\ref{conj_ZP_ML} by \cite[Theorem 5.3]{Pink05}. We also use three ingredients appearing in \cite[Theorem 5.7]{Pink05}: torsion translates of semi-abelian subvarieties of a fiber are realized by special subvarieties of the family (Lemma~\ref{lem_extension}); codimension is preserved in this realization; and the fibers of special subvarieties over their images are equidimensional (Corollary~\ref{cor_fiber_dim}).

The point requiring correction is the converse fibral assertion. An irreducible component of the intersection of a special subvariety with a fiber is a translate of a semi-abelian subvariety by Theorem~\ref{thm_ribet_type_fibers}, but it need not be a torsion translate \cite{BE20}. Hence the hypothesis of Conjecture~\ref{conj_ZP_SA} does not imply that the subvariety is Hodge generic in the ambient family. In the proof above this is replaced by passing to the special closure $R_0$ of $Z$, where $Z$ is Hodge generic by definition, and by the group-theoretic codimension estimate in Step 3.
\end{remark}

\begin{remark}\label{rem_edixhoven}
Theorem~\ref{thm_ZP_implies_ZPSA} was obtained independently by Edixhoven \cite{EdixhovenPrivate}, as announced in \cite[Remark 1.3]{BE20}, by a direct group-theoretic computation generalizing \cite{BE20}; this proof does not use the geometry of the special subvarieties of $\pgn$. Moreover, \cite[Proposition 3.9]{EdixhovenPrivate} gives a finer fibral statement than Theorem~\ref{thm_ribet_type_fibers}. The proof above instead derives the required group-theoretic input from the classification of special subvarieties.
\end{remark}

\subsection{Reformulation of the Relative Manin--Mumford Conjecture}\label{subsec_reform_RMM}

Throughout this subsection $S$ is an irreducible variety over $\CC$. Pink proposed the following relative form of Conjecture~\ref{conj_ZP_SA}.

\begin{conjecture}[{\cite[Conjecture 6.1]{Pink05}}]\label{conj_RMM_pink}
Consider an algebraic family of semi-abelian varieties $\calG\to S$ over $\CC$ and an irreducible closed subvariety $Z\subseteq\calG$ of dimension $d$ that is not contained in any proper closed subgroup scheme (with a possibly smaller base) of $\calG\to S$. Then $Z\cap\calG^{[>d]}$ is not Zariski dense in $Z$.
\end{conjecture}

If $\calG\to S$ has relative dimension $g+n$, with toric rank $n$, then for every $d<g+n$ the set $\calG^{[>d]}$ contains all points of $\calG$ which are torsion in their own fiber, by Remark~\ref{rem_torsion_in_codim}. Denoting this set of fiberwise torsion points by $\calG_{\mathrm{tor}}$, Conjecture~\ref{conj_RMM_pink} implies the following special case.

\begin{conjecture}[{\cite[Conjecture 6.2]{Pink05}}]\label{conj_RMM_pink_cor}
Consider an algebraic family of semi-abelian varieties $\calG\to S$ of relative dimension $g+n$ over an irreducible variety over $\CC$, and an irreducible closed subvariety $Z\subseteq\calG$. Assume that $Z$ is not contained in any proper closed subgroup scheme of $\calG\to S$ and that $Z(\CC)\cap\calG_{\mathrm{tor}}$ is Zariski dense in $Z$. Then $\dim Z\geq g+n$.
\end{conjecture}

Conjecture~\ref{conj_RMM_pink_cor} is the relative Manin--Mumford conjecture. For families of abelian varieties it has been proved in \cite{GaoHab_RMM}. In positive toric rank, however, Conjecture~\ref{conj_RMM_pink} is false; a counterexample is constructed in \cite{BE20}, after \cite[Theorem 2.4]{BE20}.

The example in \cite{BE20} is the Poincar\'e torsor over a CM elliptic curve $E$ over $\CC$. Let $\lambda:E\stackrel{\sim}{\to}E^{\vee}$ be the canonical principal polarization associated with $\calO_E(0_E)$, choose an endomorphism $\phi\in\End_\CC(E)$ with $\phi^{\dagger}\neq\phi$, and set $\alpha:=\lambda\circ(\phi-\phi^{\dagger})$. In the language of Section~\ref{section4}, this is the Ribet datum $(\Gamma(\alpha),1)$ of Example~\ref{ex_ribet_graph} with $n=1$, and the corresponding subvariety of Ribet type is the Ribet section $t_\phi$ over the graph of $\alpha$. The two properties established in \cite{BE20} are:
\begin{itemize}
  \item the image of $t_\phi$ is not contained in any proper subgroup scheme of $\pt_E$, and
  \item $t_\phi(x,\alpha(x))$ is torsion in the fiber of $\pt_E\to E$ for every torsion point $x$ of $E$.
\end{itemize}
Taking $Z=t_\phi(\Gamma(\alpha))$ and $d=1$, the set $\pt_E{}^{[>1]}$ is the set of points of $\pt_E$ that are torsion in their fiber, and $Z\cap\pt_E{}^{[>1]}$ is dense in $Z$. This contradicts Conjecture~\ref{conj_RMM_pink}. On the other hand, $t_\phi$ is a special subvariety of $\pt_E$ as a mixed Shimura variety.

The counterexample shows that, in the semi-abelian setting, torsion translates of subgroup schemes do not exhaust the exceptional subvarieties. By the classification theorem, the additional exceptional subvarieties are precisely those of Ribet type. They must therefore be included in the formulation.

The notion of a subvariety of Ribet type is the one of Definition~\ref{def_ribet_type_general}, and the statements below invoke it after a finite base change of $S$. Whenever we speak of a subvariety of Ribet type of $\calG\times_SS'$ for a finite surjective morphism $S'\to S$, it is understood that $S'$ is chosen so that the toric part of $\calG$ becomes split over it, since only then is the notion defined; such an $S'$ always exists, the character sheaf of the toric part being an isotrivial \'etale local system on $S$.

We now state the revised conjectures, beginning with the analogue of \cite[Conjecture 6.1]{Pink05}.

\begin{conjecture}\label{conj_re_RMM}
Let $\calG\to S$ be a semi-abelian scheme over $\CC$ with $S$ irreducible, and let $Z\subseteq\calG$ be an irreducible closed subvariety of dimension $d$ such that, for every finite surjective morphism $S'\to S$, no irreducible component of $Z\times_SS'$ is contained in a proper subvariety of Ribet type of $\calG\times_SS'$. Then $Z\cap\calG^{[>d]}$ is not Zariski dense in $Z$.
\end{conjecture}

\begin{remark}\label{rem_hypothesis_weakens}
Replacing Pink's hypothesis ``$Z$ is not contained in a proper closed subgroup scheme, with a possibly smaller base'' by ``no irreducible component of $Z\times_SS'$ is contained in a proper subvariety of Ribet type of $\calG\times_SS'$, for any finite surjective $S'\to S$'' strengthens the hypothesis, by Proposition~\ref{prop_flat_subgroup} applied over the smaller base, and therefore weakens the conjecture. This removes the counterexample of \cite{BE20}: in that example $Z$ is itself a proper subvariety of Ribet type of $\pt_E$.
\end{remark}

The revised relative Manin--Mumford statement is the following immediate consequence of Conjecture~\ref{conj_re_RMM}, by Remark~\ref{rem_torsion_in_codim}.

\begin{conjecture}[Relative Manin--Mumford, revised]\label{conj_re_RMM_cor}
Let $\calG\to S$ be a semi-abelian scheme of relative dimension $g+n$ with toric rank $n$, over an irreducible variety $S$ over $\CC$. Let $Z\subseteq\calG$ be an irreducible closed subvariety, and let $\calG_{\mathrm{tor}}$ denote the set of points of $\calG$ that are torsion in their own fiber. If $Z(\CC)\cap\calG_{\mathrm{tor}}$ is Zariski dense in $Z$, then either $\dim Z\geq g+n$, or there are a finite surjective morphism $S'\to S$ and a \emph{proper} subvariety of Ribet type of $\calG\times_SS'$ containing an irreducible component of $Z\times_SS'$.
\end{conjecture}

\begin{proposition}[Zilber--Pink implies the revised relative Manin--Mumford]\label{prop_ZP_implies_re_RMM}
Conjecture~\ref{conj_ZP} implies Conjecture~\ref{conj_re_RMM}, and hence Conjecture~\ref{conj_re_RMM_cor}.
\end{proposition}
\begin{proof}
Conjecture~\ref{conj_re_RMM} is invariant under finite surjective base change of $S$. We may therefore assume that the toric part of $\calG$ is split and that a polarization of the abelian quotient and a level structure have been chosen, as in Remark~\ref{rem_minimal_family}. We then have a classifying morphism $S\to(\Agv)^{[n]}$ and a morphism $f:\calG\to\pgn$ over it. Let $S_1\subseteq(\Agv)^{[n]}$ be the special closure of the image of $S$. Then $f$ factors through the Shimura family $\pgn|_{S_1}$ of Definition~\ref{def_shimura_family}. Let $Z'\subseteq\pgn|_{S_1}$ be the Zariski closure of $f(Z)$; the morphism $f|_Z:Z\to Z'$ is dominant.

\smallskip
\noindent\textit{Step 1: $Z'$ is Hodge generic in $\pgn|_{S_1}$.} Suppose that $Z'\subseteq R$ for a special subvariety $R\subseteq\pgn|_{S_1}$. We show that $R=\pgn|_{S_1}$. By Theorem~\ref{thm_classification}, the subvariety $R$ is of Ribet type. Hence the irreducible components of $f^{-1}(R)$ are subvarieties of Ribet type of $\calG$, the construction of Section~\ref{section4} being compatible with base change. Since $f^{-1}(R)$ contains $Z$, and $Z$ is by hypothesis contained in no proper subvariety of Ribet type of $\calG$, the case $S'=S$ of the hypothesis, one irreducible component of $f^{-1}(R)$ must be all of $\calG$. Thus $f(\calG)\subseteq R$.

It follows that $\pr(R)$ contains the image of $S$ in $(\Agv)^{[n]}$. Since $\pr(R)$ is special, it contains the special closure $S_1$ of that image; hence $\pr(R)=S_1$. Moreover, for any $s\in S(\CC)$, with image $\underline b\in S_1(\CC)$, the fiber $R_{\underline b}$ contains the image of $\calG_s$, which is the whole fiber $(\pgn)_{\underline b}$ because $f$ is fiberwise an isomorphism. By Theorem~\ref{thm_classification} and Corollary~\ref{cor_fiber_dim}, all nonempty fibers of $R\to S_1$ have the same dimension. Thus $\dim R=\dim S_1+g+n=\dim\pgn|_{S_1}$, and since $R$ is closed and irreducible in the irreducible variety $\pgn|_{S_1}$, we have $R=\pgn|_{S_1}$.

\smallskip
\noindent\textit{Step 2: applying Conjecture~\ref{conj_ZP}.} Put $d':=\dim Z'$, so $d'\leq d$. By Step 1 we may apply Conjecture~\ref{conj_ZP} to $Z'\subseteq\pgn|_{S_1}$: the set
\[
  W':=\!\!\bigcup_{\substack{\text{special }R\subseteq\pgn|_{S_1}\\ \operatorname{codim}_{\pgn|_{S_1}}R>d}}\!\!R\cap Z'
\]
is not Zariski dense in $Z'$, since the union over the special subvarieties of codimension $>d$ is contained in the union over those of codimension $>d'$. Let $\overline{W'}\subsetneq Z'$ be its Zariski closure, a proper closed subset of $Z'$.

Then $Z\cap f^{-1}\bigl(\overline{W'}\bigr)$ is a proper closed subset of $Z$. Indeed, it is closed, and it cannot be all of $Z$: otherwise $f(Z)\subseteq\overline{W'}$, contradicting the density of $f(Z)$ in $Z'$ and the properness of $\overline{W'}\subsetneq Z'$.

\smallskip
\noindent\textit{Step 3: points in $\calG^{[>d]}$.} Let
\[
  \Omega:=\!\!\bigcup_{\substack{\text{special }R\subseteq\pgn|_{S_1}\\ \operatorname{codim}_{\pgn|_{S_1}}R>d}}\!\!R .
\]
Fix $s\in S(\CC)$ and let $s'$ be its image in $S_1$. We claim that
\[
  f\bigl(\calG_s^{[>d]}\bigr)\ \subseteq\Omega .
\]
Let $S_2\subseteq S_1$ be the special closure of $s'$, and let $H\subseteq\calG_s\cong(\pgn)_{s'}$ be an algebraic subgroup with $\operatorname{codim}_{\calG_s}H>d$. By Lemma~\ref{lem_extension}(2), after a finite Shimura covering $S_2'\to S_2$ there are finitely many special subvarieties of $\pgn|_{S_2'}$, each of codimension $\operatorname{codim}_{\calG_s}H$, whose fibers above the chosen point contain the components of $H$. Their images in $\pgn|_{S_2}$ are special subvarieties of the same codimension, and these images contain the components of $H$. If $R$ is any one of them, viewed also as a special subvariety of $\pgn|_{S_1}$, then
\[
  \operatorname{codim}_{\pgn|_{S_1}}R=\operatorname{codim}_{\pgn|_{S_1}}\bigl(\pgn|_{S_2}\bigr)+\operatorname{codim}_{\pgn|_{S_2}}R\ \geq\ \operatorname{codim}_{\calG_s}H\ >\ d ,
\]
This proves the claim. Taking the union over $s\in S(\CC)$ gives $f\bigl(\calG^{[>d]}\bigr)\subseteq\Omega$. Therefore,
\[
  Z\cap\calG^{[>d]}\ \subseteq\ Z\cap f^{-1}\bigl(\Omega\cap Z'\bigr)
  =Z\cap f^{-1}(W')
  \subseteq Z\cap f^{-1}\bigl(\overline{W'}\bigr),
\]
which by Step 2 is not Zariski dense in $Z$. This is the conclusion of Conjecture~\ref{conj_re_RMM}.
\end{proof}

\begin{remark}
At present few cases of the reformulated relative Manin--Mumford conjecture are known. For abelian schemes, that is for $n=0$, it is the theorem of Gao and Habegger \cite{GaoHab_RMM}, itself the outcome of a long line of work beginning with \cite{MZ10} and passing through \cite{CMZ18}. In positive toric rank the only case we know of is that of semi-abelian surfaces, treated by Bertrand, Masser, Pillay and Zannier \cite{Known_for_semiab_surf}: families over a curve whose fibers are extensions of elliptic curves by $\gm$, equivalently the universal Poincar\'e torsor over a modular curve. That is also the setting in which the counterexample of \cite{BE20} lives, so already the first case in positive toric rank requires the subvarieties of Ribet type.

Let us also record why the semi-abelian case is worth having. Torsion in generalized Jacobians --- extensions of a Jacobian by copies of $\gm$ and of $\Ga$ --- governs the solvability of the Pell--Abel equation $A^{2}-DB^{2}=1$ in polynomials when $D$ has multiple roots, and, through a criterion of Risch, the integrability in elementary terms of the members of a pencil of algebraic differentials. Masser and Zannier obtained finiteness statements for both problems by proving the relevant relative Manin--Mumford statements \cite{MZPell}; see \cite{ZannierICM,Zannier8ECM}. Extensions by $\Ga$ are not of mixed Shimura type and are not covered by the results of this article; see \cite{Schmidt19}.
\end{remark}

\bibliographystyle{amsalpha}
\bibliography{mybibliography}

\end{document}